 \documentclass{elsart}

\usepackage{amssymb}
\usepackage{color}
\usepackage{graphicx}
\usepackage{subfigure}
\usepackage{amsmath}
\usepackage{amsfonts}
\usepackage{epstopdf}

\renewcommand{\qed}{\hfill \nobreak \ifvmode \relax \else
      \ifdim\lastskip<1.5em \hskip-\lastskip
      \hskip1.5em plus0em minus0.5em \fi \nobreak
      \vrule height0.75em width0.5em depth0.25em\fi}

\def\pd#1#2{\frac{\partial #1}{\partial #2}}
\def\dfr#1#2{\displaystyle{\frac{#1}{#2}}}

\renewcommand{\vec}[1]{\mbox{\boldmath \small $#1$}}
\def\dd#1#2{\frac{\partial #1}{\partial #2}}

 \usepackage[ 
              unicode=true,
             pdfstartview=FitH,
             CJKbookmarks=true,
             bookmarksnumbered=true,
             bookmarksopen=true,
             colorlinks=true, 
             citecolor=magenta,
             linkcolor=blue,
             linktocpage=true,
             ]{hyperref}       

\usepackage{multirow}
\allowdisplaybreaks
    \numberwithin{equation}{section}

    \newtheorem{Example}{Example} [section]

        \def\DG{RKDG methods}
        \def\CDG{Runge-Kutta CDG methods}

  \newtheorem{Remark}{Remark}[section]

\begin{document}

\begin{frontmatter}
\title{Central Runge-Kutta discontinuous Galerkin methods
for the special relativistic hydrodynamics}


 \author{Jian Zhao}
\ead{everease@163.com}
\address{HEDPS, CAPT \& LMAM, School of Mathematical Sciences, Peking University,
Beijing 100871, P.R. China}
\author[label2]{Huazhong Tang}
\thanks[label2]{Corresponding author. Tel:~+86-10-62757018;
Fax:~+86-10-62751801.}
\ead{hztang@math.pku.edu.cn}
\address{HEDPS, CAPT \& LMAM, School of Mathematical Sciences, Peking University,
Beijing 100871, P.R. China; School of Mathematics and Computational Science,
 Xiangtan University, Hunan Province, Xiangtan 411105, P.R. China}
 \date{\today{}}

\maketitle

\begin{abstract}
	This paper developes Runge-Kutta $P^K$-based central discontinuous Galerkin (CDG) methods with WENO limiter to the one- and two-dimensional special relativistic hydrodynamical (RHD) equations, $K=1,2,3$.
Different from the non-central DG methods, the \CDG{} have to find
two approximate solutions defined on mutually dual meshes.
For each mesh, the CDG
approximate solutions on its dual mesh are used to calculate the flux values in the cell and on the
cell boundary so that the approximate solutions on mutually dual meshes are coupled with
each other, and the use of numerical flux may be avoided.
The WENO limiter  is adaptively implemented via two steps:
 the ``troubled'' cells are first identified by using a modified TVB minmod function,
and then   the  WENO technique is used to
  locally reconstruct new  polynomials of degree $(2K+1)$  replacing the CDG solutions inside the ``troubled'' cells
  by  the cell average values of the CDG solutions in the neighboring cells as well as the original cell averages of the ``troubled'' cells.
Because the WENO limiter is only employed for finite ``troubled'' cells,
the computational cost can be as little as possible.
The accuracy of the CDG without the numerical dissipation is analyzed and
calculation of the flux integrals over the cells is also addressed.
Several test problems in one and two dimensions are solved by using our  \CDG{} with WENO limiter. The computations demonstrate that our methods are stable, accurate, and robust in solving
complex RHD problems.
\end{abstract}

\end{frontmatter}



\renewcommand\baselinestretch{1.1}   

\section{Introduction}

Relativistic fluid widely appears in nuclear physics, astrophysics, plasma physics,
and other fields. For example,  in the physical phenomena
such as the formation of neutron stars and black holes and the high-speed jet,
the local fluid velocity may be close to the speed of light,
at this time the relativistic effect can not be neglected
and the relativistic fluid dynamics (RHD) is needed.
Because the  RHD equations are more complicated,
their theoretical analysis is impractical so that conversely
numerical simulation has become a primary and powerful way
to study and understand the physical mechanisms in the RHDs.

The pioneering numerical work may date back to
the finite difference code via artificial viscosity for the
spherically symmetric general RHD equations in the Lagrangian coordinate
\cite{May-White1966,May-White1967}.
Wilson first
attempted to solve multi-dimensional RHD equations in the
Eulerian coordinate by using the finite difference method
with the artificial viscosity technique  \cite{Wilson:1972}.
Since 1990s, the numerical study of the RHDs began to attract
considerable attention, and various modern shock-capturing methods
with an exact or approximate Riemann solver
have been developed for the RHD equations,
the readers are referred to the early review articles
\cite{MartiRHDRiview,DEAOdyck:2004}.
Some examples on existing methods, which are extensions of Godunov type shock capturing methods,
are the upwind schemes based on local linearization
 \cite{EulderinkMel:1995,FalleKom:1996},
the two shock spproximation solvers \cite{Balsara:1994,DaiWood:1997,Mignoneetal:2005},
 flux-vector splitting scheme \cite{DonatFluxSplitRHD},
 HLL (Harten-Lax-van Leer) schemes \cite{Schneider:1993,Duncan:1994},
 HLLC (Harten-Lax-van Leer-Contact) scheme  \cite{MignoneHLLCRHD},
 non-oscillatory essentially (ENO) schemes \cite{DolezalWongRHDENO,ZannaBucciantini:2002},
 and kinetic schemes \cite{YangBeam1997, KunikKinetic2004} and so on.
Recently the second author and his co-workers developed
adaptive moving mesh method  \cite{HeAdaptiveRHD}, derived
the second-order accurate generalized Riemann problem (GRP) methods
for the one- and two-dimensional RHD equations  \cite{YangGRPRHD1D,YangGRPRHD2D},
and  the finite volume local evolution Galerkin scheme for two-dimensional RHD equations
\cite{WuEGRHD}.
Later, the third-order accurate GRP scheme in \cite{WuGRP3} was extended to
the one-dimensional RHD equations \cite{WuGRPRHD3}, and the direct Eulerian GRP scheme
was developed for the spherically symmetric general relativistic hydrodynamics \cite{Wu-Tand-SISC2016}.
The physical-constraints-preserving (PCP)
schemes were also studied for the special RHD equations recently.
The high-order accurate PCP finite difference weighted essentially non-oscillatory
(WENO) schemes and discontinuous Galerkin (DG) methods were proposed in \cite{Wu-Tand-JCP2015,Wu-Tang-RHD2016,QinShuYang-JCP2016}.
Moreover,  the set of admissible states  and the PCP schemes of
the ideal relativistic magnetohydrodynamics was studied
for the first time in \cite{Wu-Tang-RMHD2016}, where  the importance of divergence-free fields
was revealed in achieving PCP methods especially.

The DG methods
have been rapidly developed  in recent decades and  has become a kind of important methods
in computational fluid dynamics. They are easy to achieve high order accuracy, suitable for parallel computing, and adapt to complex domain boundary.
The DG method was first developed by Reed and Hill \cite{ReedFirstDG} to solve  steady-state scalar linear hyperbolic equation but it had not been widely used.
A major development of the DG method was carried out  in a series of papers \cite{CockburnDG1,CockburnDG2,CockburnDG3,CockburnDG4,CockburnDG5}, where
the DG spatial approximation was combined with explicit Runge-Kutta time discretization to develop
Runge-Kutta DG (RKDG) methods  and  a general framework of DG methods was established
for the nonlinear equation or system.
After that, the RKDG methods began to get a wide range of research and application,
 such as the Euler equations \cite{TangRKDG,momentlimiter,Remade2003},
Maxwell equations \cite{DGMaxwell}, nonlinear Dirac equations \cite{ShaoTangDiracDG} etc.
Moreover, the DG methods have also been used to solve other partial differential
equations, such as convection-diffusion type equation or system \cite{BassiDGNS,CockLDG} and
Hamilton-Jacobi equation  \cite{hu1999DGHamilton,AnaDGHJ,CDGHamiltonLi} etc.
The readers are referred to the review article \cite{CockDGReview}.
The \CDG{} \cite{LiuCentralDG}
were developed by combing \DG{} and central scheme \cite{LiuCentralSch} and  found two approximate solutions
defined on  mutually dual meshes. Although
 two approximate solutions are redundant, the numerical flux may be avoided
due to  the use of the solution on the dual mesh to calculate the
flux at the cell interface. It is one of the advantages of the central scheme.
Because the \CDG{} can be considered as a variant of \DG{}, they keep  many advantages of \DG{}, such as
compact stencil and  parallel implementation etc.
Moreover, the \CDG{} allow a larger CFL number than \DG{}
and reduce numerical oscillations for some problems. Up to now, the \CDG{} have also been used to solve the Euler equations \cite{LiuCentralDG} and the ideal magneto-hydrodynamical equations \cite{CentralMHD,arbiCDGMHD} and so on.

A deficiency of the RKDG methods is that when the strong discontinuity
appears in the solution,
the numerical oscillations should be suppressed
after each Runge-Kutta inner stage or after some complete Runge-Kutta
steps by using the nonlinear limiter, which  is a commonly used technique
of the modern shock-capturing methods for hyperbolic conservation laws.
The commonly used limiter  is  the minmod limiter, which
 limits the slope of solution such that the values of limited solution in the cell
falls in the certain interval determined by the cell average values of neighboring cells.
The  minmod limiter has good robustness but becomes
 only first-order accurate near extreme points.
 The modified TVB minmod limiter is given in \cite{CockburnDG2} and applied to the \DG{}.
It does not limit the solution near extreme points by
choosing a parameter $M$, thus  the accuracy  of \DG{} is not destroyed near the extreme point.
In general, for nonlinear equation, the parameter is dependent on the problem, and
the accuracy of $P^K$-based \DG{} for $K\ge3$ may still be destroyed because
 more than three of the higher order moments will be set to zero in the  modified TVB minmod limiter.
%
Besides those commonly used limiters,
some other limiters are porposed, such as
the moment based limiters \cite{momentlimiter}  and its improvement
 \cite{Krivodonova2007} etc. Those limiters may suppress numerical oscillations
  near the discontinuity, however, the accuracy of \DG{}
 may be reduced in the some region.

In the modern shock-capturing methods, the ENO  and WENO  methods
 are more robust than the slope limiters especially for high order schemes
 and have been widely used,see the review article
 \cite{Shu-SIEV2009}.
  An attempt was made to use them as limiters for the DG methods \cite{QiuShuWENOlimiter,QiuShuWENOlimiterUn,QiuShuWENOlimiter3D}.
The WENO limiter  first identifies  the ``troubled'' cells  by using a modified TVB minmod function,
and then   new  polynomials inside the ``troubled'' cells are locally reconstructed to
replace the DG solutions by using the  WENO technique based on the cell average values of
the DG solutions in the neighboring cells as well as the original cell averages of the ``troubled'' cells.
It is only employed for finite ``troubled'' cells,
so the computational cost can be as little as possible.

This paper proposes the Runge-Kutta $P^K$-based CDG methods with WENO limiter for the one- and two-dimensional special RHD equations, $K=1,2,3$. It
is organized as follows. Section
\ref{chap:rhdrkdg} introduces the system of special RHD equations.
 Section \ref{sec:cdgdis} proposes Runge-Kutta $P^K$-based CDG methods with WENO limiter.
 Section \ref{sec:cdgana} gives some discussions of the \CDG{}.
  Section \ref{sec:rhdcdgnum} gives several numerical examples to verify the accuracy
robustness, and effectiveness of the proposed \CDG{}.
Concluding remarks are presented in Section \ref{Section-conclusion}.

\section{Special RHD equations} 
\label{chap:rhdrkdg}
This section introduces the governing equations of the special relativistic hydrodynamics (RHD).
Similar to the non-relativistic case,
the special RHD equations may be established
by the laws of local baryon number conservation
and energy-momentum conservation \cite{Landau:1959} and
  cast into the
 following  covariant form 
\begin{equation}
\label{eq:RHDcov-form}
\begin{cases}
\partial_\alpha (\rho u^\alpha) = 0,   \\
\partial_\alpha \Big(\rho h u^\alpha u^\beta + p g^{\alpha \beta} \Big) = 0,
\end{cases}
\end{equation}
where the Greek indices
$\alpha$ and $\beta$ run from 0 to 3, $\partial_\alpha =
\partial_{x^\alpha}$ stands for the covariant derivative,
$g^{\alpha \beta}$ denotes the metric tensor and is restricted to the Minkowski tensor
throughout the paper, i.e.
 $\big( g^{\alpha \beta}\big)_{4\times4} =
\mbox{diag}\{-1,1,1,1\}$,
$\rho$, $u^\alpha$ and $p$
denote the rest-mass density,  four-velocity vector,
and pressure, respectively,
and $h$  is the specific enthalpy defined by
\begin{equation}\label{eq:RHDenh}
h = 1 + e + \dfr{p}{\rho}, 
\end{equation}
here $e$ denotes the specific internal energy.
For the sake of convenience, units in which the speed of light is equal to one will be used so that
  $x^\alpha=(t, x_1,x_2,x_3)^T$ and $u^\alpha=\gamma(1,  v_1,  v_2,
v_3)^T$, where $\gamma=1/\sqrt{1-v^2}$ is the   Lorentz factor and
$v:=\sqrt{v_1^2+v_2^2+v_3^2}$ is the size of fluid velocity.

In order to close the above system  \eqref{eq:RHDcov-form},
an equation of state (EOS) for the thermodynamical variables
\begin{equation}\label{RHDEOS0}
p =p(\rho,e),
\end{equation} is needed.
For example,  the EOS for an ideal gas can be expressed in the
$\Gamma$-law form
\begin{equation}\label{RHDgamma-law}
p = (\Gamma-1)\rho e,
\end{equation}
where $\Gamma$ is the adiabatic index, taken as 5/3  for the mildly relativistic case and
4/3 for the ultra-relativistic case.

The  covariant form of special RHD equations \eqref{eq:RHDcov-form}
 is usually written into
a time-dependent system of conservation laws in the laboratory frame as follows
\begin{equation}\label{eqn:RHDconeqn}
\displaystyle\frac{\partial \vec{U}}{\partial t} +
\sum^d_{i=1}\frac{\partial \vec{F}_i(\vec{U})}{\partial x_i}=0,
\end{equation}
where $\vec{U}$ is  conservative variable vector
and  $\vec{F}_i$  denotes and the flux vector in the $x_i$ direction, $i=1,\cdots,d$.
For example, in the case of  $d=3$,  the detailed expressions of $\vec{U}$ and $\vec{F}_i$ are
\begin{align}
\begin{aligned}
\vec{U} =& \Big(D, m_1, m_2, m_3, E\Big)^T,\\
\vec{F}_1 =& \Big(Dv_1, m_1 v_1 + p,  m_2 v_1,  m_3 v_1, m_1\Big)^T,\\
\vec{F}_2 =& \Big(Dv_2, m_1 v_2,  m_2 v_2+p,  m_3 v_2,
m_2\Big)^T,\\
\vec{F}_3 =& \Big(Dv_3, m_1 v_3,  m_2 v_3,  m_3 v_3+p, m_3\Big)^T,
\end{aligned}
\label{eqn:coneqn-04}
\end{align}
here $D=\rho \gamma$, $m_i=\rho h \gamma^2 v_i$ ,
and $E=\rho h\gamma^2-p$  denote the mass,
$x_i$-momentum, and energy densities relative to the laboratory frame, respectively.

The formal structure of \eqref{eqn:RHDconeqn}
is  identical to that of the three-dimensional
non-relativistic Euler equations.
The momentum equations in \eqref{eqn:RHDconeqn}  are only with a Lorentz-contracted momentum density replacing
$\rho v_i$ in the non-relativistic Euler equations.
When the fluid velocity is small ($v\ll 1=c$) and the velocity of the internal (microscopic) motion of the fluid particles is small, the RHD equations \eqref{eqn:RHDconeqn} reduce to the non-relativistic Euler equations.
 The system  \eqref{eqn:RHDconeqn} also satisfies the properties of the rotational invariance and the homogeneity  as well as the hyperbolicity in time when \eqref{RHDgamma-law} is used, see \cite{ZhaoHeTang2012}. However,
in comparison to the non-relativistic Euler equations,
 a strong coupling between the hydrodynamic equations is introduced
 and  additional numerical difficulties are posed due to
the relations between the laboratory quantities (the mass density $D$, the momentum density $m_i$,
and the energy density $E$) and the quantities in the local rest frame (the mass density $\rho$,
and the fluid velocity $ v_i$, the internal energy density $e$).
Especially,
the flux $\vec F_i$ in   \eqref{eqn:RHDconeqn}
can not be formulated in an explicit form of the conservative vector $\vec U$
and  the physical constraints $E\geq D$, $\rho>0$, $p>0$, and $v<1$
have to be fulfilled.
Thus, in practical computations of the system \eqref{eqn:RHDconeqn}
by using the shock-capturing methods, the primitive variable vector $\vec V=(\rho,v_1,\cdots,v_d,p)^T$
has to be first recovered from the known conservative vector $\vec U=(D,m_1,\cdots,m_d,E)^T$
at each time step by numerically solving  a nonlinear pressure equation such as
\begin{equation}\label{pressure}
E + p = D\gamma + \displaystyle\frac{\Gamma}{\Gamma-1} p\gamma^2,
\end{equation}
where  $\gamma=(1-|\vec m|^2/(E+p)^2)^{-1/2}$.
Any standard root-finding algorithm, e.g. Newton's iteration, may be used to
solve \eqref{pressure} to get the pressure, and then $\gamma$, $\rho$,
$e$, $h$, and $v_i$ in order, the readers are referred to \cite{ZhaoTang2013}
for the choice of initial guess.




\section{\CDG{}}
\label{sec:cdgdis}
\begin{figure}[!htbp]
	\centering{}
	\includegraphics[width=10cm]{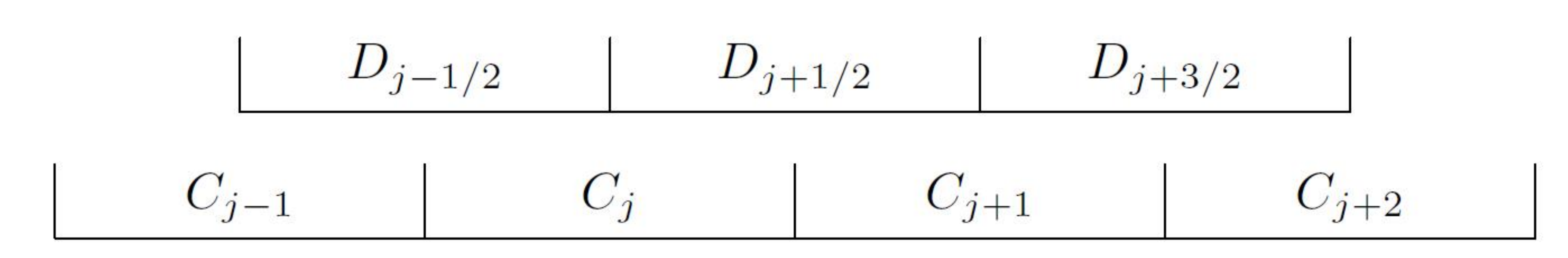}
	\caption{Schematic diagram of mutually dual meshes for the 1D \CDG{}.}\label{fig:Cdg1DMesh}
\end{figure}

This section gives the Runge-Kutta central DG methods for the hyperbolic conservation laws.
%
For the sake of simplicity,  one-dimensional scalar equation
\begin{equation}\pd{u(x,t)}{t}+\pd{f(u(x,t))}{x}=0,\quad
x\in\Omega ,\label{eq:scalar1D}\end{equation}
is taken as  an example  to introduce the  \CDG{} \cite{LiuCentralDG}.
Similar to the non-central \DG{},  the \CDG{} also employ
the discontinuous Galerkin finite element in
the spatial discretization and the explicit Runge-Kutta method for the time discretization.
Their difference between them is that \CDG{} need two  mutually dual meshes, see
the one-dimensional schematic diagram in Fig. \ref{fig:Cdg1DMesh} for the mesh
 $\{C_{j}=(x_{j-\frac12},x_{j+\frac12}),\forall j\in\mathbb Z\}$ and
 its dual mesh $\{D_{j+\frac 12}=(x_j,x_{j+1}),\forall j\in \mathbb Z\}$.

The aim of 
\CDG{} is to find two approximate solutions $u^C_h(x,t)$ and $u^D_h(x,t)$
such that at any time $t\in(0,T]$, they belong to the following finite spaces respectively
\begin{align*}
&\mathcal{V}^C:=\left\{v(x)\in L^1(\Omega)|   v(x)\in
\mathbb  P^{K}(C_{j}),\mbox{} ~x\in C_{j}\subset\Omega,
\forall j \right\},\\
&\mathcal{V}^D:=\left\{w(x)\in L^1(\Omega)|~  w(x)\in
\mathbb  P^{K}(D_{j+\frac12}),\mbox{ } ~x\in D_{j+\frac12}\subset\Omega,
\forall j \right\}.\end{align*}
Consider the CDG scheme for the approximate solution $u^C_h$.
Multiplying \eqref{eq:scalar1D} by the test function $v(x)\in \mathbb  P^{K}(C_{j})$
and  integrating it over the cell $C_{j}$ by parts
gives
\begin{equation}\frac{d}{dt}\int_{C_{j}} u v dx=\int_{C_{j}}
f(u)\dd{v}{x}dx -f\big(u(x_{j+\frac{1}{2}},t)\big) v (x_{j+\frac{1}{2}})+
f\big(u(x_{j-\frac{1}{2}},t)\big) v (x_{j-\frac{1}{2}}).
\label{eq:RHDCDG1}
\end{equation}
If replacing $u$ at the left- and right-hand sides of \eqref{eq:RHDCDG1} with the approximate solution $u^C_h$
and  $u^D_h$, respectively,
then one has
\begin{align}
\nonumber \frac{d}{dt}\int_{C_{j}}u^C_h  v
dx=&\frac{1}{\tau_{max}}\int_{C_{j}
} (u^D_h-u^C_h) v(x) dx+\int_{C_{j}} f(u^D_h) \dd{v}{x}
dx\\
-&f\big(u^D_h(x_{j+\frac{1}{2}},t)\big)
v(x_{j+\frac{1}{2}})+f\big(u^D_h(x_{j-\frac{1}{2}},t)\big) v(x_{j-\frac{1}{2}}),
\label{eq:weakuch}\end{align}
where the first term at the right-hand side of \eqref{eq:weakuch}
denotes the numerical dissipation term borrowing from the central scheme \cite{LiuCentralSch},
 and $\tau_{max}$ denotes the maximum time step size allowed by the CFL condition.
Because the approximate solution $u^D_h$ is continuous at the boundary of cell $C_{j}$,
the fluxes $f (u(x_{j\pm \frac12},t) )$  may be directly evaluated and thus numerical flux is not required
in the CDG methods.

If
using $\phi_{j}^{(l)}(x)$, $l=0,...K$, to denote
 a basis of the space $\mathbb  P^{K}(C_{j})$,  then $u^C_h$ may be expressed
 as
 $$u^C_h(x,t)=\sum_{l=0}^K u^{C,(l)}_{j}(t)\phi_{j}^{(l)}(x)=:u^C_j(x,t),\quad
~x\in C_{j}.$$
Replacing  $v(x)$ in \eqref{eq:weakuch} with the basis function $\phi_{j}^{(\ell)}(x)$
and using the numerical quadrature with $q$ points to calculate
the integral with flux $f(u)$ gives the semi-discrete scheme
for $u^C_h$ as follows
\begin{align}
\nonumber
\sum_{l=0}^K\big(\int_{C_{j}} \phi_j^{(\ell)}(x)  \phi_j^{(l)}(x)
dx\big) \frac{d u^{C,(l)}_j}{dt}=\frac{1}{\tau_{max}}\int_{C_{j}
} (u^D_h-u^C_h) \phi_j^{(\ell)}(x) dx\\
\nonumber
-f\big(u^D_h(x_{j+\frac{1}{2}},t)\big)
\phi_j^{(\ell)}(x_{j+\frac{1}{2}})+f\big(u^D_h(x_{j-\frac{1}{2}},t)\big)
\phi_j^{(\ell)}(x_{j-\frac{1}{2}})
\\ +h_{j} \sum_{m=1}^q \omega_m^Cf\big(u^D_h(x^C_m,t)\big) \dd{}{x}\phi_j^{(\ell)}(x_m^C), \quad \ell=0,1,\cdots,K,\label{eq:semiuch}\end{align}
where $h_j=x_{j+\frac{1}{2}}-x_{j-\frac{1}{2}}$, and $x^C_m$ and
$\omega^C_m$  denote the point and weight for the numerical integration over the cell $C_j$, $m=1,\cdots,q$.
It needs to be pointed out that the first term at the right-hand side of \eqref{eq:semiuch}
is an integral of piecewise polynomial and may be exactly calculated.

The semi-discrete scheme for the approximate solution $u^D_h$ may be similarly derived.
If choosing a basis of the space $\mathbb  P^{K}(D_{j+\frac{1}{2}})$ as
$\{\varphi_{j+\frac{1}{2}}^{(l)}(x),~ l=0,...K\}$, then the   semi-discrete
scheme for $u^D_h$ is given as follows
\begin{align}
\nonumber
\sum_{l=0}^K\big(\int_{D_{j+\frac{1}{2}}} \varphi_{j+\frac{1}{2}}^{(\ell)}(x)  \varphi_{j+\frac{1}{2}}^{(l)}(x)
dx\big)
\frac{d u^{D,(l)}_{j+\frac{1}{2}}}{d t}=\frac{1}{\tau_{max}}\int_{D_{j+\frac{1}{2}}}
(u^C_h-u^D_h) \varphi_{j+\frac{1}{2}}^{(\ell)}(x) dx\\
\nonumber
-f\big(u^C_h(x_{j+1},t)\big)
\varphi_{j+\frac12}^{(\ell)}(x_{j+1})+f\big(u^C_h(x_{j},t)\big)
\varphi_{j+\frac12}^{(\ell)}(x_{j})
\\+h_{j+\frac{1}{2}} \sum_{m=1}^q
\omega_m^Df\big(u^C_h(x^D_m,t)\big) \dd{}{x} \varphi_{j+\frac12}^{(\ell)}(x_m^D), \quad \ell=0,1,\cdots,K.\label{eq:semiudh}\end{align}

\begin{Remark}\label{whycdgdis}
	If $f(u)=u$, then  
	two approximate solutions satisfy \cite{LiuCDGError}
	\begin{equation}
	\frac{1}{2}\frac{d}{dt}\int_\Omega\big((u_h^C)^2+(u_h^D)^2\big)dx=-\frac{1}{\tau_{max}}\int_\Omega
	(u_h^C-u_h^D)^2dx\le 0.
	\end{equation}
It can be known,  that is why  the term $\frac{1}{\tau_{max}}\int_{C_{j} } (u^D_h-u^C_h) v(x)dx$ in \eqref{eq:weakuch} is called as numerical dissipation.
\end{Remark}

\begin{Remark}\label{re2gauss}	
It is worth noting that	  the flux within the cell $C_j$ in
\eqref{eq:weakuch} is evaluated by using $u^D_h$,
but the approximate solution $u^D_h$ is not continuous at $x_j$,
thus before using the numerical integration to evaluate
the integral of flux over $C_j$ in \eqref{eq:weakuch}, one has to split
it into two parts 
	\begin{equation}\label{eq:flux-integral} \int_{x_{j-\frac{1}{2}}}^{x_{j+\frac{1}{2}}} f(u^D_h)\dd{v}{x}dx=
	\int_{x_{j-\frac{1}{2}}}^{x_j} f(u^D_{j-\frac{1}{2}})\dd{v}{x}
	dx+\int_{x_j}^{x_{j+\frac{1}{2}}} f(u^D_{j+\frac{1}{2}})\dd{v}{x}
	dx,	\end{equation}
and then use Gaussian quadrature with $K+1$ points
to calculate two integrals at the right-hand side of the above equation.
The flux integral in \eqref{eq:semiudh} should be similarly treated.
Section \ref{sec:cdgnodiff} will give a further discussion on the evaluation of
such flux integral.
 \end{Remark}

Both semi-discrete CDG schemes \eqref{eq:semiuch} and \eqref{eq:semiudh} may be cast into the following abstract form
$$\dd{\vec{U}}{t}=\vec L(\vec{U}),$$
which is a nonlinear system of ordinary differential equation of $\vec U$ with respect to $t$,
and thus  the time derivatives    may be further approximated
to give the fully-discrete CDG methods may be derived for the degrees of freedom or the moments
by using
the third-order accurate TVD (total variation
diminishing) Runge-Kutta method \cite{ShuOsher1988}
\begin{align}
\begin{aligned}
\vec{U}^{(1)}=&\vec{U}^n+\Delta t_n \vec L(\vec{U}^n),\\
\vec{U}^{(2)}=&\frac{3}{4}\vec{U}^n+\frac{1}{4}\big(\vec{U}^{(1)}+\Delta
t_n \vec L(\vec{U}^{(1)})\big),\\
\vec{U}^{n+1}=&\frac{1}{3}\vec{U}^n+\frac{2}{3}\big(\vec{U}^{(2)}+\Delta
t_n \vec L(\vec{U}^{(2)})\big),
\end{aligned}\label{eq:RK3}
\end{align}
or the fourth-order accurate non-TVD Runge-Kutta method
\begin{align}
\begin{aligned}
&\vec{U}^{(1)}=\vec{U}^n+\frac{1}{2}\Delta t_n \vec L(\vec{U}^n),\\
&\vec{U}^{(2)}=\vec{U}^n+\frac{1}{2}\Delta
t_n \vec L(\vec{U}^{(1)}),\\
&\vec{U}^{(3)}=\vec{U}^n+\Delta
t_n \vec L(\vec{U}^{(2)}),\\
&\vec{U}^{n+1}=\frac{1}{3}\big(\vec{U}^{(1)}+2\vec{U}^{(2)}+3\vec{U}^{(3)}-\vec{U}^n+\frac{1}{2}\Delta
t_n \vec L(\vec{U}^{(3)})\big).\end{aligned}\label{eq:RK4}
\end{align}
and so on.

As mentioned above, $\tau_{max}$ is determined by using the CFL condition.
After determining $\tau_{max}=\tau_n$ at $t=t_n$, the practical time stepsize $\Delta t_n$ should satisfy
$0< \Delta t_n \le \tau_n$. If denoting $\theta=\Delta t_n/\tau_n$, then $\theta\in (0,1]$.
For hyperbolic equations, $\Delta t_n$ may usually be taken as $\tau_n$, that is,
$\theta=1$.

Although the \CDG{} are only introduced   for one-dimensional scalar conservation law,
their extension to one-dimensional RHD equations and two-dimensional rectangular mesh
is easy and direct. Similar to the non-central \DG{}, the limiting procedure
is necessary for the \CDG{} when the solution contains strong discontinuity.
The WENO limiting procedure in Section 3.3 of \cite{ZhaoTang2013} may be directly and independently
applied to the solutions $u^C_h$ and $u^D_h$ of \CDG{} by the following two steps:
\begin{itemize}
  \item identify the ``troubled'' cells in the meshes $\{C_j\}$ and $\{D_{j+ 1/2}\}$,
  namely, those cells which might need the limiting procedure,
  \item replace the CDG solution polynomials $u^C_h$ and $u^D_h$ in those   ``troubled'' cells with WENO reconstructed polynomials of degree $(2K+1)$, denoted by  $u^{C,WENO}_h$ and $u^{D,WENO}_h$,
      which maintain the original cell averages (conservation) and the accuracy, but have less numerical oscillation.
\end{itemize}
In order to save the length of paper,  those details are omitted here.

\section{Some discussions of \CDG{}}
\label{sec:cdgana}
In comparison to the non-central \DG{}, the \CDG{} has an additional numerical dissipation term.
Remark \ref{whycdgdis} has shown that such dissipation term is important for the $L^2$ stability of CDG methods.
This section discusses the accuracy of \CDG{} without the numerical dissipation
in order to understand that the impact of numerical dissipation term on the accuracy of CDG methods.
Furthermore, this section will also discuss the calculation of flux integral over the cell
 mentioned in Remark \ref{re2gauss}. 

\subsection{Accuracy of CDG methods without numerical dissipation}
\label{sec:p1cdgnodiff}
The  semi-discrete version of \CDG{} without numerical dissipation
can be written as follows
\begin{align}\begin{aligned}
\frac{d}{dt}\int_{C_{j}} u^C_h  v
dx=& \int_{C_{j}} f(u^D_h) \dd{v}{x}
dx
-f\big(u^D_h(x_{j+\frac{1}{2}},t)\big) v(x_{j+\frac{1}{2}})
\\ &
+f\big(u^D_h(x_{j-\frac{1}{2}},t)\big) v(x_{j-\frac{1}{2}}), \quad \forall v(x)\in
\mathbb P^{K}(C_{j}),
\\
\frac{d}{dt}\int_{D_{j+\frac{1}{2}}} u^D_h  w
dx=&\int_{D_{j+\frac{1}{2}}} f(u^C_h) \dd{w}{x}
dx
-f\big(u^C_h(x_{j+1},t)\big)
w(x_{j+1})
\\ &+f\big(u^C_h(x_{j},t)\big) w(x_{j}), \quad \forall w(x)\in
\mathbb P^{K}(D_{j+\frac{1}{2}}).
\end{aligned}\label{eq:seminodiff}
\end{align}

The accuracy of $P^1$-based CDG methods without numerical dissipation
is first discussed here by using the Fourier method similar to \cite{zhang2003analysis,ShuAnaDGSV}.
Use \eqref{eq:seminodiff} to solve the scalar equation
\begin{equation}u_t+u_x=0,\quad x\in[0,2\pi],
\label{eq:oneDlinear}
\end{equation}
subject to  the initial condition $u(x,0)=\sin(x)$.

\begin{figure}[!htbp]
	\centering{}
	\includegraphics[width=0.3\textwidth]{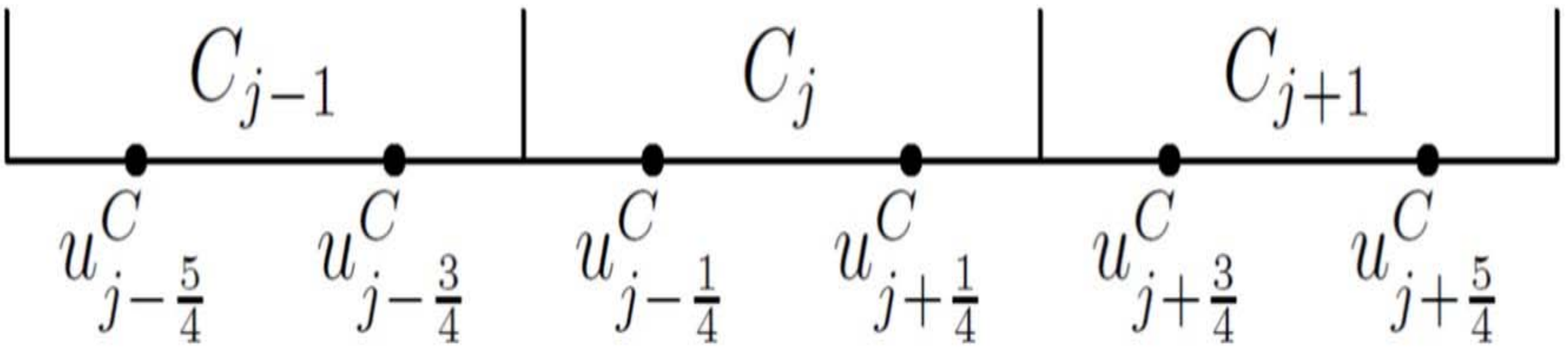}\\[-0.5cm]
	\caption{Schematic diagram of degrees of freedom.}\label{fig:cdgfdm}
\end{figure}

For the CDG solution on the mesh $\{C_j\}$.
As shown in Fig.~\ref{fig:cdgfdm},
for the sake of convenience, the degrees of freedom are chosen as
the function values at $2N$ points distributed with equal distance
$$u^C_{j-\frac{1}{4}},~ u^C_{j+\frac{1}{4}},~ j=1,\cdots, N,$$
instead of all order moments $\{u^{C,(l)}_{j}(t)\}$.
%

Within the cell $C_j$, the solution $u^C_j(x,t)$ may be expressed as
$$u^C_{j}(x,t)=u^C_{j-\frac14}(t)\phi_{j-\frac14}(x)+u^C_{j+\frac 14}(t)\phi_{j+\frac 14}(x),$$
where $\phi_{j\pm \frac14}(x)$ denote the basis functions,
given by
$$\phi_{j-\frac 14}(x)=\frac{1}{2}-\frac{2(x-x_j)}{h},\quad
\phi_{j+\frac 14}(x)=\frac{1}{2}+\frac{2(x-x_j)}{h}.$$
Similarly, within the cell $D_{j+\frac 12}$, the solution $u^D_{j+\frac 12}(x,t)$ may be written as follows
$$u^D_{j+\frac12}(x,t)=u^D_{j+\frac 14}(t)\varphi_{j+\frac 14}(x)+u^D_{j+\frac 34}(t)\varphi_{j+\frac 34}(x).$$
If replacing $v(x)$ in \eqref{eq:seminodiff} with $\phi_{j-\frac 14}(x)$ and $\phi_{j+\frac 14}(x)$,
and $w(x)$ with $\varphi_{j+\frac 14}(x)$ and $\varphi_{j+\frac 34}(x)$, respectively,
 and performing the mass matrix inversion, 
 then one has
\begin{align}
\begin{aligned}
\frac{d u^C_{j-\frac{1}{4}}}{dt}
=& \frac{1}{4h} \big( 5u^D_{j-\frac{3}{4}}-u^D_{j-\frac{1}{4}}-5
u^D_{j+\frac{1}{4}}+u^D_{j+\frac{3}{4}}\big),\\
\frac{d u^C_{j+\frac{1}{4}}}{dt}
=&\frac{1}{4h} \big(-u^D_{j-\frac{3}{4}}+5 u^D_{j-\frac{1}{4}}+ u^D_{j+\frac{1}{4}}-5 u^D_{j+\frac{3}{4}}\big),
\\
\frac{d u^D_{j+\frac{1}{4}}}{dt}
=& \frac{1}{4h} \big( 5u^C_{j-\frac{1}{4}}-u^C_{j+\frac{1}{4}}-5 u^C_{j+\frac{3}{4}}+u^C_{j+\frac{5}{4}}\big),\\
\frac{d u^D_{j+\frac{3}{4}}}{dt}
=&\frac{1}{4h} \big(-u^C_{j-\frac{1}{4}}+5 u^C_{j+\frac{1}{4}}+ u^C_{j+\frac{3}{4}}-5 u^C_{j+\frac{5}{4}}\big),
\end{aligned}\nonumber 
\end{align}
which may be rewritten as follows
\begin{equation}
\frac{d}{dt}\begin{pmatrix}
u^C_{j-\frac 14}\\
u^C_{j+\frac14}\\
u^D_{j+\frac14}\\
u^D_{j+\frac34}\end{pmatrix}= \vec A \begin{pmatrix}
u^C_{j-\frac 54}\\
u^C_{j-\frac34}\\
u^D_{j-\frac34}\\
u^D_{j-\frac14}
\end{pmatrix}
+\vec B\begin{pmatrix}
u^C_{j-\frac 14}\\
u^C_{j+\frac14}\\
u^D_{j+\frac14}\\
u^D_{j+\frac34}
\end{pmatrix}
+\vec C\begin{pmatrix}
u^C_{j+\frac34}\\
u^C_{j+\frac54}\\
u^D_{j+\frac54}\\
u^D_{j+\frac74}
\end{pmatrix},
\label{eq:evolvedeg}
\end{equation}
where three coefficient matrices are respectively given by
$$\vec A= \begin{pmatrix} 0 & 0 & \frac{5}{4h} &
-\frac{1}{4h}\\
0 & 0 & -\frac{1}{4h} &
\frac{5}{4h}\\
0 & 0 & 0 & 0\\
0 & 0 & 0 & 0\\
\end{pmatrix},
~\vec B= \begin{pmatrix} 0 & 0 &  -\frac{5}{4h} & \frac{1}{4h}\\
0 &  0 & \frac{1}{4h} & -\frac{5}{4h}\\
\frac{5}{4h} & -\frac{1}{4h}  & 0 & 0\\
-\frac{1}{4h} & \frac{5}{4h} & 0 & 0
\end{pmatrix},
~\vec C=\begin{pmatrix}  0 & 0 & 0 & 0\\
0 & 0 & 0 & 0\\
-\frac{5}{4h} &  \frac{1}{4h} & 0 & 0\\
\frac{1}{4h} &  -\frac{5}{4h} & 0 & 0
\end{pmatrix}.
$$
Because the solution of \eqref{eq:oneDlinear} is periodic and the mesh is uniform,
the solution of \eqref{eq:evolvedeg} may be assumed to be of the following form 
\begin{equation}\label{eq:solexp}
\begin{pmatrix}
u^C_{j-\frac 14}(t)\\
u^C_{j+\frac14}(t)\\
u^D_{j+\frac14}(t)\\
u^D_{j+\frac34}(t)\end{pmatrix}=\begin{pmatrix}
\hat{u}^C_{q,-\frac14}(t)\\
\hat{u}^C_{q,\frac14}(t)\\
\hat{u}^D_{q,\frac14}(t)\\
\hat{u}^D_{q,\frac34}(t)\end{pmatrix} e^{iqx_j},
\end{equation}
where $i=\sqrt{-1}$ is the imaginary unit.
Substituting \eqref{eq:solexp} into \eqref{eq:evolvedeg} gives
\begin{equation}
\frac{d}{dt}\begin{pmatrix}
\hat{u}^C_{q,-\frac 14}\\
\hat{u}^C_{q,\frac14}\\
\hat{u}^D_{q,\frac14}\\
\hat{u}^D_{q,\frac34}\end{pmatrix}
=\vec G(q,h)\begin{pmatrix}
\hat{u}^C_{q,-\frac14}\\
\hat{u}^C_{q,\frac14}\\
\hat{u}^D_{q,\frac14}\\
\hat{u}^D_{q,\frac34}\end{pmatrix}
,\end{equation}
where $\vec G(q,h)$ denotes the amplification matrix and is defined by
\begin{equation}
\label{eq:matG}
\vec G(q,h)=\vec Ae^{-i\alpha}+\vec B+\vec Ce^{i\alpha},
\end{equation}
here $\alpha=qh$ and $h$ denotes the spatial stepsize.
Four eigenvalues of $\vec G$ are 
$$\lambda_{1,2}=\pm\frac{1}{h}\sqrt{2\cos\alpha-2},\quad \lambda_{3,4}=\pm\frac{3}{2h}\sqrt{2\cos\alpha-2},$$
and corresponding right eigenvectors may be taken as follows
$$\vec r_1=\begin{pmatrix} \frac{e^{-\alpha i}-1}{\sqrt{2\cos \alpha-2}}\\
\frac{e^{-\alpha i}-1}{\sqrt{2\cos\alpha-2}}\\
1\\
1\\
\end{pmatrix},~\vec r_2=\begin{pmatrix} -\frac{e^{-\alpha
		i}-1}{\sqrt{2\cos \alpha-2}}\\
-\frac{e^{-\alpha i}-1}{\sqrt{2\cos \alpha-2}}\\
1\\
1\\
\end{pmatrix},
~\vec r_3=\begin{pmatrix} -\frac{e^{-\alpha i}-1}{\sqrt{2\cos \alpha-2}}\\
\frac{e^{-\alpha i}-1}{\sqrt{2\cos \alpha-2}}\\
-1\\
1\\
\end{pmatrix},~\vec r_4=\begin{pmatrix} \frac{e^{-\alpha i}-1}{\sqrt{2\cos\alpha-2}}\\
-\frac{e^{-\alpha i}-1}{\sqrt{2\cos\alpha-2}}\\
-1\\
1\\
\end{pmatrix}.
$$
Thus the solution of  Eq.~\eqref{eq:evolvedeg} may be expressed as follows
\begin{align*}
& \big(u^C_{j-\frac14}(t),~u^C_{j+\frac14}(t),~u^D_{j+\frac14}(t),~u^D_{j+\frac34}(t)\big)^T\\
&=c_1 e^{iqx_j+\lambda_1 t} \vec r_1+c_2 e^{iqx_j+\lambda_2 t} \vec r_2+ c_3
e^{iqx_j+\lambda_3 t}\vec r_3+c_4 e^{iqx_j+\lambda_4 t}\vec r_4,
\end{align*}
where $c_i$, $1\le i\le 4$, are four undetermined coefficients.

Let us discuss the accuracy of methods. Take $q=1$ and define
\begin{equation}\label{eq:CDGACini}
u^C_{j-\frac{1}{4}}(0)=e^{ix_{j-\frac{1}{4}}},\ \
u^C_{j+\frac{1}{4}}(0)=e^{ix_{j+\frac{1}{4}}},\ \
u^D_{j+\frac{1}{4}}(0)=e^{ix_{j+\frac{1}{4}}}, \ \
u^D_{j+\frac{3}{4}}(0)=e^{ix_{j+\frac{3}{4}}},
\end{equation}
then their imaginary parts satisfy the initial condition
$u(x,0)=\sin(x)$.
It should be pointed out that,
 the initial degrees of freedom in the DG methods  are generally derived
 by using  the $L^2$ projection to the initial condition, but
 the  approach for setting initial value \cite{ShuAnaDGSV} is used here
 and does not effect the final results on accuracy. 

The undetermined coefficients $c_i$, $1\le i\le 4$, may be determined by
\eqref{eq:CDGACini} as follows
\begin{align*}
c_{1,2}&= \pm \frac{1}{2}\frac{\sqrt{2\cos h-2}}{e^{-ih}-1}\cos
\frac{h}{4} + \frac{1}{4}(e^{\frac{3}{4}ih}+e^{\frac{1}{4}ih}), \\
c_{3,4}&= \pm \frac{i}{2}\frac{\sqrt{2\cos h-2}}{e^{-ih}-1}\sin \frac{h}{4} + \frac{1}{4}(e^{\frac{3}{4}ih}-e^{\frac{1}{4}ih}),
\end{align*}
which may give the expression of $u^C_{j-\frac{1}{4}}(t)$.
Using the Taylor expansion to the imaginary part of $u^C_{j-\frac{1}{4}}(t)$
with respect to $h$ gives
\begin{equation}\label{eq:P1NDCDGRes}
Im\{u^C_{j-\frac{1}{4}}(t)\}=\sin(x_{j-\frac{1}{4}}-t)
+\frac{h}{4}\big(\cos(x_{j-\frac{1}{4}}-t)-\cos(x_{j-\frac{1}{4}}-\frac{3}{2}t)\big)+O(h^2).
\end{equation}
The solutions of $P^1$-based methods with numerical dissipation satisfy \cite{LiuCDGError}
$$ Im\{u^C_{j-\frac{1}{4}}(t)\}=\sin(x_{j-\frac{1}{4}}-t)
+ \sigma \sin(x_{j-\frac{1}{4}}-t) h^2+O(h^3),$$
where $\sigma$ is a constant only depending on $\tau_{max}/h$.
Comparing them gives their obvious difference.
The similar differences may be given by using the Taylor expansion
to the expression of $u^C_{j+\frac{1}{4}}(t),u^D_{j+\frac{1}{4}}(t)$, and
$u^D_{j+\frac{3}{4}}(t)$.
From the above analysis, it is seen that the $P^1$-based \CDG{} without
numerical dissipation term are only first-order accurate in space.

\begin{figure}[htbp]
	\centering{}
	\includegraphics[width=0.55\textwidth]{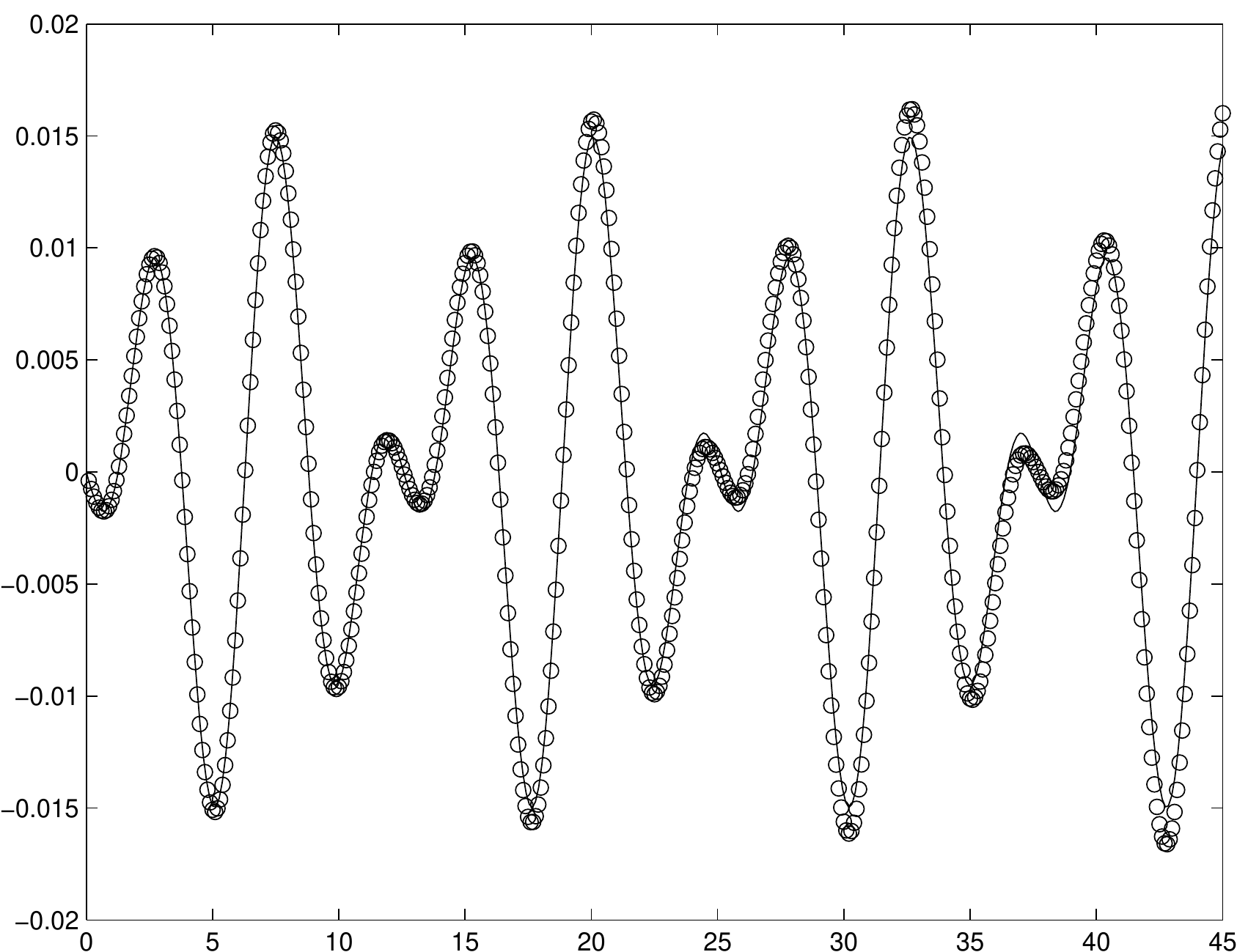}\\
	\caption{The symbol ``$\circ$'' denotes the deference between the numerical and exact solutions of \eqref{eq:oneDlinear} at point $x_{j-\frac{1}{4}}$, $u^C_h(x_{j-\frac{1}{4}},t)-u(x_{j-\frac{1}{4}},t)$,
while the solid line denote the result obtained by \eqref{eq:P1NDCDGRes}, $ h\big(\cos(x_{j-\frac{1}{4}}-t)-\cos(x_{j-\frac{1}{4}}-\frac{3}{2}t)\big)/4$,
where $j=50$. Numerical solution is obtained by using the $P^1$-based \CDG{} without
numerical dissipation but with $200$ cells. }\label{fig:cdginferror}
\end{figure}

In the following, we use the $P^1$-based method to solve Eq. \eqref{eq:oneDlinear}
in order to numerically demonstrate
\eqref{eq:P1NDCDGRes}.
In order to reduce the errors arising from the time discretization,
the fourth-order accurate Runge-Kutta method \eqref{eq:RK4} is used with the time stepsize
$\Delta t=0.01h$.
Fig.~\ref{fig:cdginferror} shows the time evolution of
error at the point $x_{j-\frac{1}{4}}$.
Except for a few moments,
numerical results are highly consistent with  the theoretical result  given by \eqref{eq:P1NDCDGRes}.
 Moreover, Table~\ref{tab:acnodf} presents $l^1$ and $l^\infty$ errors
of solution at $t=15$, as well as  the results estimated in theory.
It is seen that the  numerical results are in good agreement with the theoretical analysis.

\begin{table}
	\centering{}
	\caption{$l^1$ and $l^\infty$ errors and orders at $t=15$ obtained by the $P^1$-based \CDG{} without
numerical dissipation for Eq.~\eqref{eq:oneDlinear}.  The  fourth-order Runge-Kutta method is employed and $h$ denotes the spatial stepsize. }
	\begin{tabular}{|c|c|c|c|c|c|c|c|c|}
		\hline
		&\multicolumn{4}{|c|}{Numerical results}&\multicolumn{4}{|c|}{Theoretical results} \\
		\hline
		$h$&$l^1$ error &  order & $l^
		\infty$ error & order &  $l^1$ error & order & $l^\infty$ error & order\\
		\hline
		$2\pi/40$&1.90e-01& &9.57e-02& &1.80e-01& &9.13e-02& \\
		\hline
		$2\pi/80$  &9.13e-02& 1.05&4.64e-02& 1.04&8.98e-02& 1.00&4.53e-02& 1.01\\
		\hline
		$2\pi/160$&4.51e-02& 1.02&2.28e-02& 1.03&4.49e-02& 1.00&2.25e-02& 1.01\\
		\hline
		$2\pi/320$&2.25e-02& 1.00&1.13e-02& 1.01&2.25e-02& 1.00&1.12e-02& 1.00\\
		\hline
		$2\pi/640$&1.12e-02& 1.00&5.62e-03& 1.01&1.12e-02& 1.00&5.61e-03& 1.00\\
		\hline
	\end{tabular}\label{tab:acnodf}
\end{table}

It is difficult to use the above Fourier method to accuracy of
$P^2$- and $P^3$-based \CDG{}   without numerical dissipation.
For this reason, the numerical experiments are provided to replace the above Fourier method.
Table~\ref{tab:acnodfP23} lists errors and orders of solutions
obtained by using the $P^2$- and $P^3$-based \CDG{}   without numerical dissipation.
It is seen that the convergence rate of $P^2$-based \CDG{}   without numerical dissipation
is essentially consistent with the predicated value 3,
but the convergence rate of $P^3$-based \CDG{}   without numerical dissipation is about
3, lesser than  the predicated value 4.

\begin{table}
	\centering{}
	\caption{Same as Table \ref{tab:acnodf} except for $P^2$- and $P^3$- based \CDG{}.
}
	\begin{tabular}{|c|c|c|c|c|c|c|c|c|}
		\hline
		&\multicolumn{4}{|c|}{$P^2$}&\multicolumn{4}{|c|}{$P^3$} \\
		\hline
		$h$&$l^1$ error & order & $l^
		\infty$ error & order &  $l^1$ error & order & $l^\infty$ error & order\\
		\hline
		$2\pi/40$&6.67e-05& &2.62e-05& &3.24e-06& &1.83e-06& \\
		\hline
		$2\pi/80$&6.80e-06& 3.29&2.04e-06& 3.68&4.01e-07& 3.01&2.24e-07& 3.03\\
		\hline
		$2\pi/160$&8.92e-07& 2.93&3.73e-07& 2.46&5.01e-08& 3.00&2.77e-08& 3.01\\
		\hline
		$2\pi/320$&9.24e-08& 3.27&3.22e-08& 3.53&6.26e-09& 3.00&3.45e-09& 3.01\\
		\hline
		$2\pi/640$&1.20e-08& 2.94&4.79e-09& 2.75&7.82e-10& 3.00&4.31e-10& 3.00\\
		\hline
	\end{tabular}\label{tab:acnodfP23}
\end{table}

\subsection{Discussion on the flux integrals over the cell}
\label{sec:cdgnodiff}

\begin{figure}[!htbp]
	\centering{}
	\includegraphics[width=0.3\textwidth]{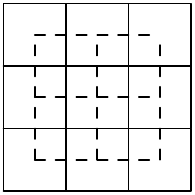}
	\caption{Schematic diagram of two-dimensional dual mesh for \CDG{}. }\label{fig:Cdg2DMesh}
\end{figure}

As mentioned in Remark~\ref{re2gauss}, because the DG solution $u^D_h$
is discontinuous at the point $x_j$ which is an internal point of
the cell $C_{j}=(x_{j-\frac{1}{2}},x_{j+\frac{1}{2}})$,
the integral of ``flux'' $f(u^D_h) \dd{v}{x} $ over $C_j$ 
becomes \eqref{eq:flux-integral}. 
If the Gaussian quadrature is used to evaluate such flux integral,
then  the numerical integration point number is twice
the non-central \DG{}.
When \CDG{} are used to solve two-dimensional conservation laws on the
dual meshes displayed
in Fig.~\ref{fig:Cdg2DMesh}, the the numerical integration point number becomes
four times that of the non-central \DG{}.

In order to reduce the computational cost of
 numerical integration,
 an attempt may be considered that only approximate
solution on the dual mesh is used to evaluate the flux on the cell boundary for the DG approximations
 on the mesh, thus the integrals of DG solution over the dual cell may be avoided and
the cost of numerical integration is hopefully reduced.
Specifically,
if using $$\int_{C_{j}} f(u^C_h) \dd{v}{x}
dx$$
to replace
$$\int_{C_{j}} f(u^D_h) \dd{v}{x}
dx,$$
then then semi-discrete CDG methods may be expressed as follows
\begin{align}\begin{aligned}
\frac{d}{dt} \int_{C_{j}} u^C_h & vdx=\frac{1}{\tau_{max}}\int_{C_{j}}(u^D_h-u^C_h) v(x)  dx+\int_{C_{j}} f(u^C_h)\dd{v}{x} dx\\
-&f\big(u^D_h(x_{j+\frac{1}{2}},t)\big)
v(x_{j+\frac{1}{2}})+f\big(u^D_h(x_{j-\frac{1}{2}},t)\big)
v(x_{j-\frac{1}{2}}), \ \forall v(x)\in \mathbb
P^{K}(C_{j}),\\
\frac{d}{dt}\int_{D_{j+\frac{1}{2}}} u^D_h & w
dx=\frac{1}{\tau_{max}}\int_{D_{j+\frac{1}{2}}} (u^C_h-u^D_h) w(x)
dx+\int_{D_{j+\frac{1}{2}}} f(u^D_h)\dd{w}{x} dx\\
-&f\big(u^C_h(x_{j+1},t)\big) w(x_{j+1})+f\big(u^C_h(x_{j},t)\big) w(x_{j}), \ \forall    w(x)\in \mathbb P^{K}(D_{j+\frac{1}{2}}).
\end{aligned}\label{eq:semiuchnew}
\end{align}
In the following, the \CDG{} based on \eqref{eq:semiuchnew} is called as the new \CDG{},
otherwise the old \CDG{}.
%
A natural problem is whether such change does effect the stability?
Consider the $P^1$-based methods.
If using the similar way to that in Section \ref{sec:p1cdgnodiff} and some
simple algebraic operations, and applying the scheme \eqref{eq:semiuchnew}  to Eq. \eqref{eq:oneDlinear},
then the evolution equation of the degrees of freedom  may be derived as follows
\begin{equation}
\frac{d}{dt}\begin{pmatrix}
u^C_{j-\frac 14}\\
u^C_{j+\frac 14}\\
u^D_{j+\frac 14}\\
u^D_{j+\frac 34}\end{pmatrix}=\vec A \begin{pmatrix}
u^C_{j-\frac 54}\\
u^C_{j-\frac 34}\\
u^D_{j-\frac 34}\\
u^D_{j-\frac 14}
\end{pmatrix}
+\vec B\begin{pmatrix}
u^C_{j-\frac 14}\\
u^C_{j+\frac 14}\\
u^D_{j+\frac 14}\\
u^D_{j+\frac 34}
\end{pmatrix}
+\vec C\begin{pmatrix}
u^C_{j+\frac 34}\\
u^C_{j+\frac 54}\\
u^D_{j+\frac 54}\\
u^D_{j+\frac 74}
\end{pmatrix},
\label{eq:uabc}
\end{equation}
where
$$\vec A= \begin{pmatrix} 0 & 0 & \frac{1}{16\tau_{max}}+\frac{5}{4h} &
\frac{13}{16\tau_{max}}+\frac{5}{4h}\\
0 & 0 & \frac{-1}{16\tau_{max}}-\frac{1}{4h} &
\frac{3}{16\tau_{max}}-\frac{1}{4h}\\
0 & 0 & 0 & 0\\
0 & 0 & 0 & 0\\
\end{pmatrix},
$$
$$\vec B= \begin{pmatrix} -\frac{1}{\tau_{max}}-\frac{3}{2h} & -\frac{3}{2h} & \frac{3}{16\tau_{max}}+\frac{1}{4h} &
\frac{-1}{16\tau_{max}}+\frac{1}{4h}\\
\frac{3}{2h} & -\frac{1}{\tau_{max}}+\frac{3}{2h} & \frac{13}{16\tau_{max}}-\frac{5}{4h} &
\frac{1}{16\tau_{max}}-\frac{5}{4h}\\
\frac{1}{16\tau_{max}}+\frac{5}{4h} &
\frac{13}{16\tau_{max}}+\frac{5}{4h}  & -\frac{1}{\tau_{max}}-\frac{3}{2h} & -\frac{3}{2h}\\
\frac{-1}{16\tau_{max}}-\frac{1}{4h} &
\frac{3}{16\tau_{max}}-\frac{1}{4h} & \frac{3}{2h} & -\frac{1}{\tau_{max}}+\frac{3}{2h}
\end{pmatrix},
$$
and
$$\vec C=\begin{pmatrix}  0 & 0 & 0 & 0\\
0 & 0 & 0 & 0\\
\frac{3}{16\tau_{max}}+\frac{1}{4h} &  \frac{-1}{16\tau_{max}}+\frac{1}{4h} & 0 & 0\\
\frac{13}{16\tau_{max}}-\frac{5}{4h} &  \frac{1}{16\tau_{max}}-\frac{5}{4h} & 0 & 0
\end{pmatrix}.
$$
Because the mesh is uniform and the periodic condition is considered
here, the solution is still assumed to be 
\begin{equation}
\nonumber
\begin{pmatrix}
u^C_{j-\frac 14}(t)\\
u^C_{j+\frac 14}(t)\\
u^D_{j+\frac 14}(t)\\
u^D_{j+\frac 34}(t)\end{pmatrix}=\begin{pmatrix}
\hat{u}^C_{q,-\frac 14}(t)\\
\hat{u}^C_{q,\frac 14}(t)\\
\hat{u}^D_{q,\frac 14}(t)\\
\hat{u}^D_{q,\frac 34}(t)\end{pmatrix} e^{iqx_j}.
\end{equation}
If denoting $\vec
u(t)=\big(\hat{u}^C_{q,-\frac 14}(t),\hat{u}^C_{q,\frac 14}(t),\hat{u}^D_{q,\frac 14}(t),\hat{u}^D_{q,\frac 34}(t)\big)^T$,
then Eq.~\eqref{eq:uabc} reduces to
\begin{equation}\label{eq:odenew}
\frac{d}{dt} \vec u(t)=\vec G \vec u(t),
\end{equation}
where the definition of amplification matrix $\vec G(q,h)$ is the same as that in \eqref{eq:matG},
and its four eigenvalues are
\begin{align}\begin{aligned}
\lambda_1=\frac{1}{\tau_{max}}\Big(-1+\frac{1}{8}\big(e^{\frac{\alpha
		i}{2}}( 4\mu-1)-e^{\frac{-\alpha
		i}{2}}( 4\mu+1)+\sqrt{a}\big)\Big),\\
\lambda_2=\frac{1}{\tau_{max}}\Big(-1+\frac{1}{8}\big(e^{\frac{\alpha
		i}{2}}( 4\mu-1)-e^{\frac{-\alpha
		i}{2}}( 4\mu+1)-\sqrt{a}\big)\Big),\\
\lambda_3=\frac{1}{\tau_{max}}\Big(-1-\frac{1}{8}\big(e^{\frac{\alpha
		i}{2}}( 4\mu-1)-e^{\frac{-\alpha
		i}{2}}( 4\mu+1)+\sqrt{b}\big)\Big),\\
\lambda_4=\frac{1}{\tau_{max}}\Big(-1-\frac{1}{8}\big(e^{\frac{\alpha
		i}{2}}( 4\mu-1)-e^{\frac{-\alpha
		i}{2}}( 4\mu+1)-\sqrt{b}\big)\Big),
\end{aligned}\label{eq:newshemeeig}
\end{align}
here $\mu=\tau_{max}/h$ denotes the CFL number, and $a$ and $b$ are given by
\begin{align*}
a=42-32\mu^2+(16\mu^2+24 \mu
-3)e^{\alpha i}+(16\mu^2-24\mu-3) e^{-\alpha
	i}+96 \mu (e^{\frac{1}{2}\alpha i}-e^{-\frac{1}{2}\alpha
	i}),\\
b=42-32\mu^2+(16\mu^2+24 \mu
-3)e^{\alpha i}+(16\mu^2-24\mu-3) e^{-\alpha
	i}-96 \mu (e^{\frac{1}{2}\alpha i}-e^{-\frac{1}{2}\alpha i}).
\end{align*}

Because the above expressions  of eigenvalues are more complicated,
the special case of small $\alpha$ is only considered here.
Using the Taylor expansions to \eqref{eq:newshemeeig} with respect to $\alpha$
gives 
\begin{align}
\nonumber
\lambda_1=&\frac{1}{\tau_{max}}\big(-\frac{1}{2}+2\mu\alpha
i+(\frac{4}{3}\mu^2+\frac{1}{16})\alpha^2\big)+O(\alpha^3),\\
\nonumber
\lambda_2=&\frac{1}{\tau_{max}}(-2-\mu\alpha  i-\frac{4}{3}\mu^2\alpha^2)+O(\alpha^3),\\
\nonumber
\lambda_3=&\frac{1}{\tau_{max}}(-\frac{3}{2}-\frac{1}{16}\alpha^2)+O(\alpha^3),\\
\nonumber
\lambda_4=&-\frac{1}{\tau_{max}}\mu \alpha i+O(\alpha^3).
\end{align}

In the following, we discuss the stability of the fully discrete version of
\eqref{eq:odenew} with Runge-Kutta  time discretizations
and the time stepsize $\Delta t_n=\tau_{max}$.
If the first-order accurate Euler method is employed, then
the fully discrete scheme becomes
\begin{equation}
\vec{u}^{n+1}=\big(\vec I+\Delta t_n  \vec G \big)\vec{u}^{n}.\label{eq:rk1total}
\end{equation}
It is easy to get that for small
$\alpha$, the inequality $$|1+\Delta t_n \lambda_2|=1+\frac{11}{6}\mu^2\alpha^2+O(\alpha^3)>1,$$
holds, thus the fully discrete scheme \eqref{eq:rk1total} is unstable.
It is similar to the old $P^1$-based \CDG{}.

If the second-order
Runge-Kutta method
\begin{align}\label{eq:RK2}\begin{aligned}
\vec{u}^{(1)}=&\vec{u}^n+\Delta t_n  \vec G \vec{u}^n,\\
\vec{u}^{n+1}=&\frac{1}{2}\vec{u}^n+\frac{1}{2}\big(\vec{u}^{(1)}+\Delta
t_n  \vec G \vec{u}^{(1)}\big),
\end{aligned}\end{align}
is used to discretize the time derivative, then the fullly-discrete
scheme    may be formed as follows
\begin{equation}
\vec{u}^{n+1}=\big(\vec I+\Delta t_n  \vec G+ \frac{1}{2}\Delta t_n^2  \vec G^2\big)\vec{u}^{n}.\label{eq:rk2total}
\end{equation}
If the inequality
$$\varrho(\mu):=\max\limits_{\substack{1\le i\le 4\\ \alpha}}|1+\Delta t_n\lambda_i+\frac{\Delta t_n^2
	\lambda_i^2}{2}|\le 1,$$
holds, then the scheme \eqref{eq:rk2total} is stable.
However, when $\alpha$ is smaller, the previous analysis
tells us that
$$\varrho(\mu)\ge|1+\Delta t_n\lambda_2+\frac{\Delta t_n^2 \lambda_2^2}{2}|=1+\frac{4}{3}\mu^2
\alpha^2+O(\alpha^3)>1.$$
It means that
for any small $\mu$, the method \eqref{eq:rk2total} is unstable,
but the old $P^1$-based \CDG{}    with second-order accurate
Runge-Kutta methods \eqref{eq:RK2} are stable under the certain CFL condition.

If the higher-order Runge-Kutta time discretization or $K>1$,
then it is difficult to analyze analytically  its stability.
 For this reason, the CFL numbers are numerically estimated.
Table \ref{tab:cfl} lists the admissible maximum CFL numbers of
new \CDG{} \eqref{eq:semiuchnew} with $\nu$th order Runge-Kutta method,
 and old \CDG{} as well as  \DG{}.
 It is seen that the CFL numbers of new methods are smaller than the old,  especially for the $P^2$-based
  \CDG{}, but the difference between the new and old $P^3$-based methods is very small.
  Thus we  may expect that the new  $P^3$-based method
  is likely to improve the computational efficiency.

\begin{table}
	\centering{}
	\caption{Numerically estimated maximum CFL numbers, where $K$ denotes the degree of polynomial basis,
and $\nu$ is the order of Runge-Kutta method.}
	\begin{tabular}{|c|c|c|c|c|c|c|c|c|c|}
		\hline
		&\multicolumn{3}{|c|}{non-central DG}    &\multicolumn{3}{|c|}{old CDG
			} & \multicolumn{3}{|c|}{new CDG}\\
		\hline
		$K$& 1&2&3 &1 &2 & 3  &1 &2 & 3 \\
		\hline{}
		$\nu=2$  & 0.333 & - & - & 0.439 & - & -  & - &- &- \\
		\hline{}
		$\nu=3$  & 0.409 & 0.209 & 0.130  &  0.588 & 0.330 & 0.224
		& 0.335 & 0.146 & 0.145 \\
		\hline{}
		$\nu=4$ & 0.464 & 0.235 & 0.145  & 0.791 & 0.472 & 0.316
		& 0.306 &        0.162 &0.149 \\
		\hline
	\end{tabular}\label{tab:cfl}
\end{table}

\begin{Remark}
The CFL number of 	\CDG{} is dependent on the size of $\theta=\Delta t_n/\tau_{max}$.
In general, if $\theta$ is smaller, the CFL number  may become bigger.
For example.  if $\theta=0.3$ and the third-order explicit
	Runge-Kutta time discretization is employed,
then the CFL number of old $P^1,P^2,P^3$-based \CDG{} may be taken as $2.57,1.5$, and 1,
respectively.
\end{Remark}

\begin{Remark}
The maximum CFL number of the old $P^1$-based \CDG{} with second-order explicit Runge-Kutta time discretization
is about 0.439, which is lesser than that in \cite{LiuCDGError}.
Moreover, our numerical experiments show that when
the CFL number  $\mu=0.44$,
the old $P^1$-based \CDG{} with second-order explicit Runge-Kutta time discretization
becomes unstable in solving ~\eqref{eq:oneDlinear} because $\varrho\approx 1.00014$.
\end{Remark}
\begin{Remark}
	It is worth mentioning that the CFL number of new $P^1$-based \CDG{} with fourth-order Runge-Kutta
time discretization is lesser than with the third-order  Runge-Kutta method.
This situation is not too common.
\end{Remark}

In order to demonstrate
the accuracy of   new methods and further compare them to
the old, 
the \CDG{} are used to solve
the initial-boundary value problem of two-dimensional Burgers equation
\begin{equation}
u_t+{\big(\frac{u^2}{2}\big)}_x+{\big(\frac{u^2}{2}\big)}_y=0,
\label{eq:burgers}
\end{equation}
with the initial data  $u(x,y,0)=0.5+\sin\big(\pi(x+y)/2\big)$,
 the computational domain  $[0,4]\times[0,4]$,  and
 the periodic boundary conditions.
Table~\ref{tab:burODnewold} gives
the errors and orders at $t=0.5/\pi$ obtained by using
the new and old \CDG{} with or without limiter in global.
 Up to the output time  $t=0.5/\pi$,  the solution is still smooth.
Those data show that two kinds of  $P^K$-based \CDG{} may achieve the theoretical order $K+1$,
and the global use of WENO limiter may keep the accuracy of \CDG{}.
Table~\ref{tab:cmpbur} presents the CPU times for  two kinds of  $P^K$-based   \CDG{}.
It is seen that for the scalar equation, the advantage of new methods
is not obvious in comparison to the old,
but we may expect that the new methods may
exhibiting great advantage in solving the RHD equations.

\begin{table}[!htbp]
	\centering
	\caption{The $l^1$ errors and orders $t=0.5/\pi$ of the new and old $P^K$-based \CDG{} for the Burgers equation \eqref{eq:burgers}. The fourth order Runge-Kutta time discretization and $N\times N$ cells are used. }
	\begin{tabular}{|c|c|c|c|c|c|c|c|c|c|}
		\hline
		\multirow{3}{20pt}{}
		&\multirow{2}{2pt}{}
		&\multicolumn{4}{|c|}{without limiter}&\multicolumn{4}{|c|}{with limiter in global}\\
		\cline{3-10}
		& & \multicolumn{2}{|c|}{new method}&\multicolumn{2}{|c|}{old method}&
		\multicolumn{2}{|c|}{new method}&\multicolumn{2}{|c|}{old method}\\
		\cline{2-10}
		& N &  $l^1$  error &  order &  $l^1$ error & order &  $l^1$ error & order
&  $l^1$ error & order \\
		\hline
		\multirow{6}{20pt}{$P^{1}$}
		&10& 4.60e-01& --& 4.74e-01& --& 1.64e+00& --& 1.16e+00& --\\
		\cline{2-10}
		&20&1.10e-01& 2.06&1.13e-01& 2.07&4.82e-01& 1.77&3.07e-01& 1.92\\
		\cline{2-10}
		&40&2.80e-02& 1.98&2.85e-02& 1.98&1.18e-01& 2.03&6.93e-02& 2.15\\
		\cline{2-10}
		&80&6.97e-03& 2.00&7.09e-03& 2.01&3.09e-02& 1.94&1.60e-02& 2.11\\
		\cline{2-10}
		&160&1.75e-03& 2.00&1.77e-03& 2.00&7.94e-03& 1.96&4.10e-03& 1.97\\
		\cline{2-10}
		&320&4.36e-04& 2.00&4.43e-04& 2.00&1.99e-03& 2.00&1.04e-03& 1.98\\
		\hline
		\multirow{6}{20pt}{$P^{2}$}
		&10& 8.13e-02& --& 7.98e-02& --& 4.40e-01& --& 3.10e-01& --\\
		\cline{2-10}
		&20&1.18e-02& 2.78&1.20e-02& 2.74&6.02e-02& 2.87&3.95e-02& 2.97\\
		\cline{2-10}
		&40&1.41e-03& 3.07&1.52e-03& 2.98&5.67e-03& 3.41&3.12e-03& 3.67\\
		\cline{2-10}
		&80&1.74e-04& 3.02&1.91e-04& 2.99&3.18e-04& 4.15&2.15e-04& 3.86\\
		\cline{2-10}
		&160&2.16e-05& 3.01&2.40e-05& 2.99&2.74e-05& 3.54&2.28e-05& 3.24\\
		\cline{2-10}
		&320&2.69e-06& 3.00&3.01e-06& 3.00&3.11e-06& 3.14&2.74e-06& 3.06\\
		\hline
		\multirow{6}{20pt}{$P^{3}$}
		&10& 3.25e-02& --& 2.95e-02& --& 3.45e-01& --& 2.44e-01& --\\
		\cline{2-10}
		&20&1.84e-03& 4.15&1.80e-03& 4.04&3.04e-02& 3.51&2.04e-02& 3.58\\
		\cline{2-10}
		&40&1.26e-04& 3.87&1.28e-04& 3.81&9.61e-04& 4.98&5.89e-04& 5.12\\
		\cline{2-10}
		&80&7.75e-06& 4.02&8.40e-06& 3.94&1.43e-05& 6.07&9.74e-06& 5.92\\
		\cline{2-10}
		&160&4.84e-07& 4.00&5.40e-07& 3.96&4.18e-07& 5.10&3.90e-07& 4.64\\
		\cline{2-10}
		&320&3.03e-08& 4.00&3.43e-08& 3.98&2.39e-08& 4.12&2.40e-08& 4.02\\
		\hline
	\end{tabular}
	\label{tab:burODnewold}
\end{table}

\begin{table}[!htbp]
	\centering
	\caption{CPU times (second) for new and old  \CDG{} solving the initial-boundary problem of Burgers equation \eqref{eq:burgers}.
$320\times 320$ cells. }
	\begin{tabular}{|c|c|c|c|c|}
		\hline
		\multirow{2}{20pt}{} & \multicolumn{2}{|c|}{ without limiter
			}&\multicolumn{2}{|c|}{with limiter in global}\\
		\cline{2-5}
		& new   & old  & new   & old   \\
		\hline
		$P^1$& 130.3 & 88.8 & 167.6  & 104.4\\
		\hline
		$P^2$& 878.5 & 585.4 &  1251.6 &  729.9\\
		\hline
		$P^3$& 2590.5 & 2560.7 & 3580.4  & 3004.1\\
		\hline
	\end{tabular}
	\label{tab:cmpbur}
\end{table}

\section{Numerical results}
\label{sec:rhdcdgnum}

 This section uses our $P^K$-based \CDG{} with WENO limiter
 presented in the last section, $K=1,2,3$,
 to solve several initial value problems or initial-boundary-value problems
 of one- and two-dimensional RHD equations
 in order to demonstrate the accuracy and effectiveness of \CDG{}.
 The \CDG{} will be compared to the \DG{}. Moreover,
because the solutions of \CDG{} on two mutually dual meshes are almost identical each other,
only the solution on one mesh $\{C_j\}$ or $\{C_{j,k}\}$ will be shown in the following.

\subsection{1D case}

For the 1D computations, the uniform mesh is used, that is, the spatial step size $h_{j+\frac12}$ is constant.
The CFL numbers $\mu$ of $P^1$-, $P^2$-, $P^3$-based \CDG{}
are taken as $0.4,~0.3,~0.2$, respectively,
respectively, and $\theta=\Delta t_n /\tau_{n} =1$.
Unless otherwise stated, $M=50$ is used in the TVB modified minmod function and the third-order accurate TVD
Runge-Kutta~\eqref{eq:RK3} is employed and the time step size is determined by
\begin{equation}
\label{eq:CDGCFL1D}
\Delta t_n=\theta \tau_{n}=\frac{\theta \mu h_{j+\frac12}}{\max\limits_{i,j}\big
	\{ |\lambda^{(i)}(\vec U^{C,(0)}_{j})|,|\lambda^{(i)}(\vec
	U^{D,(0)}_{j+\frac{1}{2}})| \big\}},
\end{equation}
where the eigenvalues $\lambda^{(i)}(\vec U)$ may be found in \cite{ZhaoTang2013}. 

\begin{Example}[Riemann Problem 1]\label{exRHDRM1DT1}\rm
	The initial data are 
	$$(\rho,v_1,p)(x,0)=\begin{cases} (1,0.9,1),&\text{$ x<0.5$,}
	\\ (1,0,10),&\text{$x>0.5,$}\end{cases}$$
and $\Gamma=4/3$. As the time increases, the initial discontinuity will be decomposed into a slowly left-moving shock wave, a contact discontinuity, and a right-moving shock wave.

Fig.~\ref{fig:RHDRM1DT1rho} presents the densities at $t=0.4$ calculated by using the \CDG{} and \DG{}.
As can be seen from those plots, the numerical solutions  of \CDG{} and \DG{}
are in good agreement with the exact solutions, but there exist obvious oscillations in the densities behind the left-moving shock wave obtained by using the $P^2$- and $P^3$-based \DG{}, while no obvious oscillation is observed in the densities obtained by \CDG{}.
Such phenomenon 
is also observed in the velocities and pressures, see
  Figs.~\ref{fig:RHDRM1DT1vel}~ and~\ref{fig:RHDRM1DT1pre}.
  The ``troubled'' cells identified by the \DG{}  is more than  the \CDG{}, see Fig.~\ref{fig:RHDRM1DT1cell}.

\end{Example}
\ifx\outnofig\undefine
\begin{figure}[!htbp]
	\centering{}
	\begin{tabular}{cc}
		\includegraphics[width=0.35\textwidth]{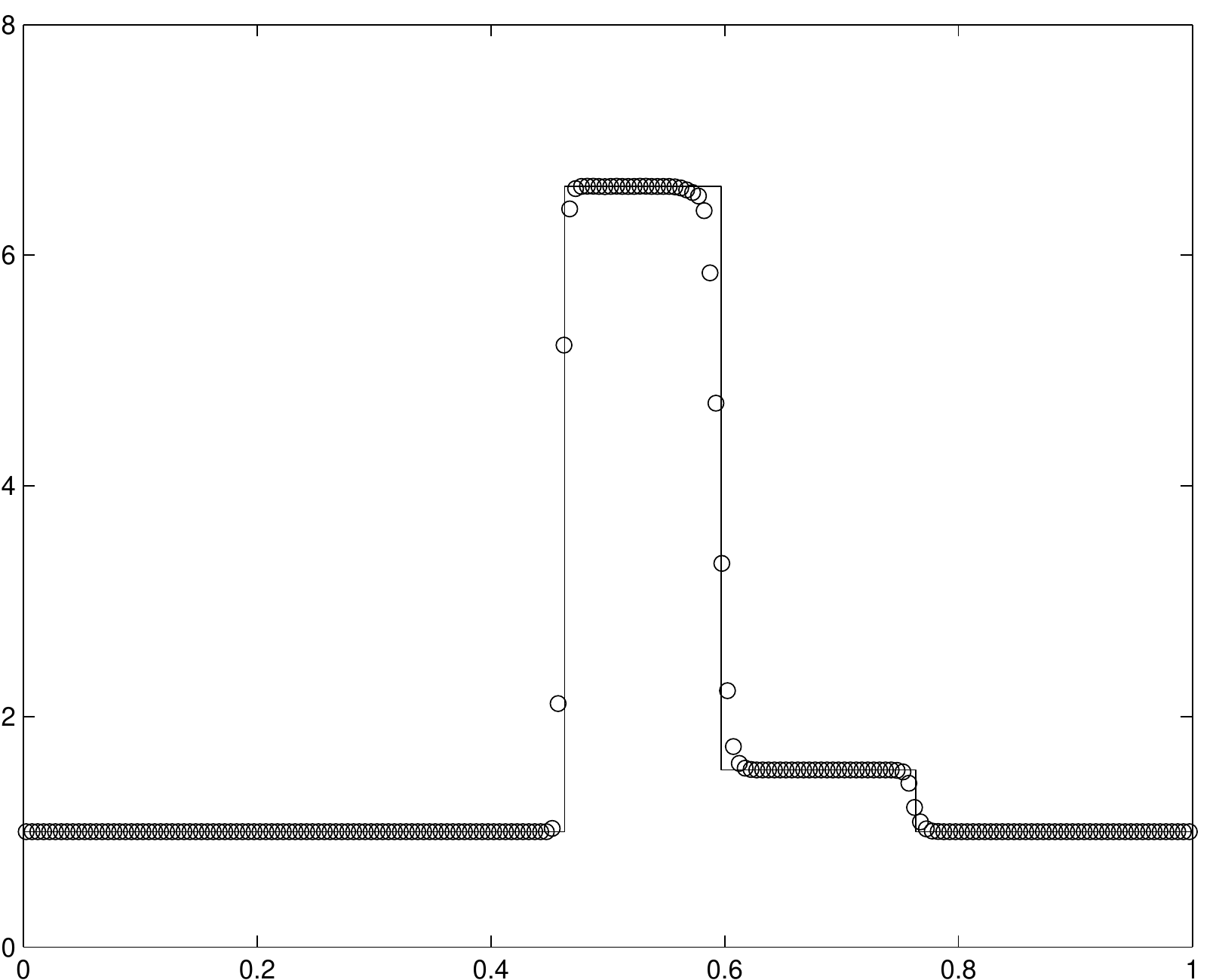}&
		\includegraphics[width=0.35\textwidth]{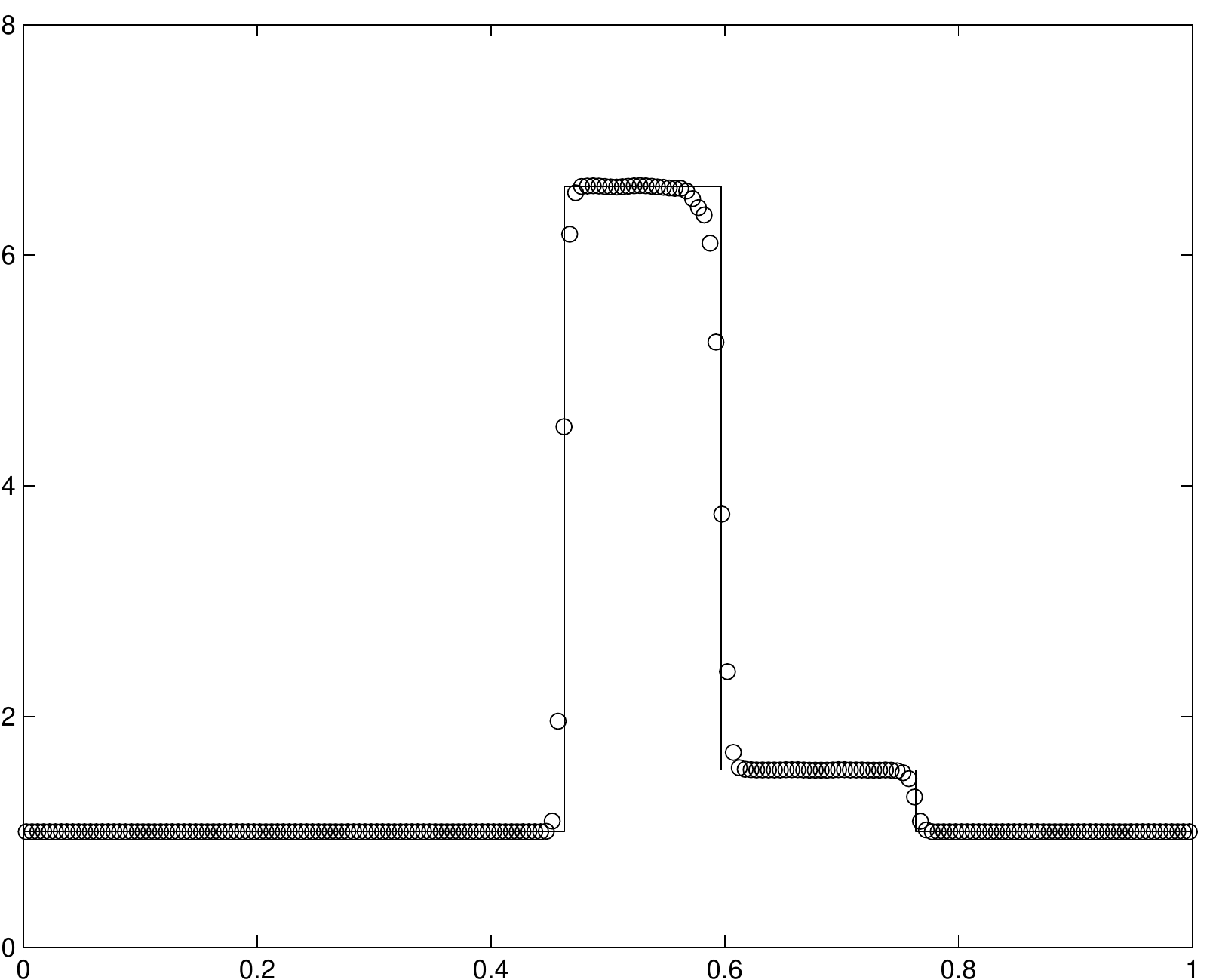}\\
		\includegraphics[width=0.35\textwidth]{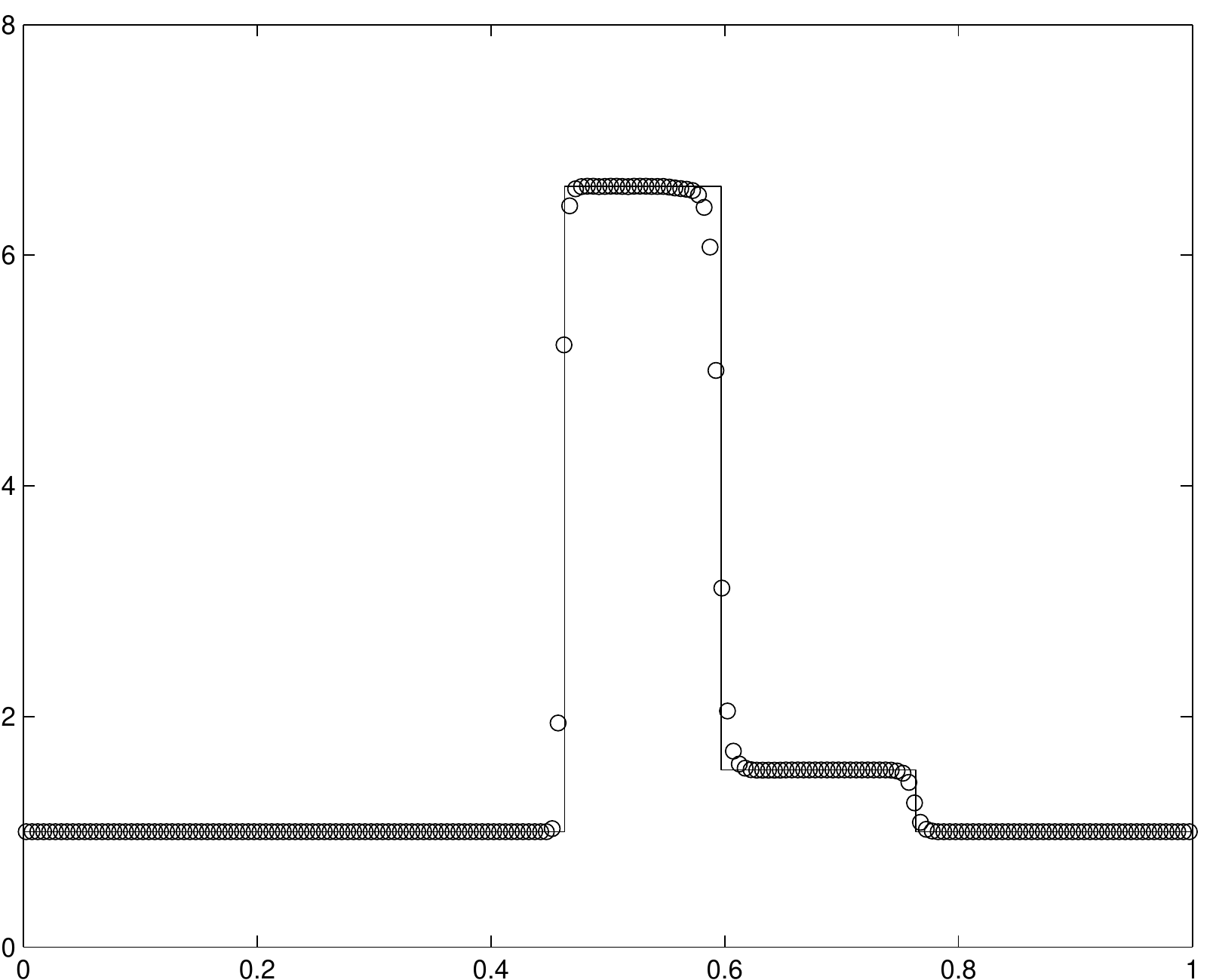}&
		\includegraphics[width=0.35\textwidth]{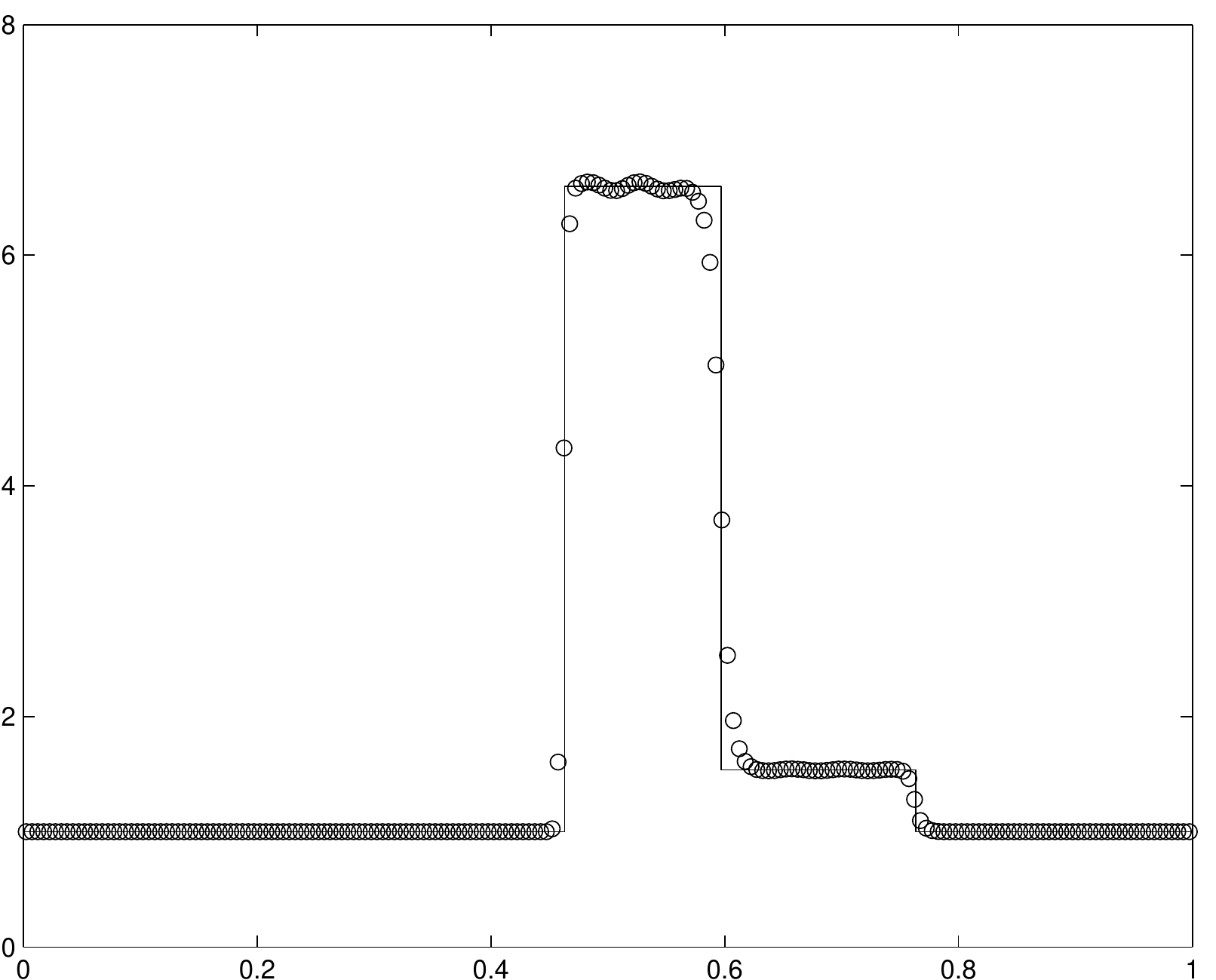}\\
		\includegraphics[width=0.35\textwidth]{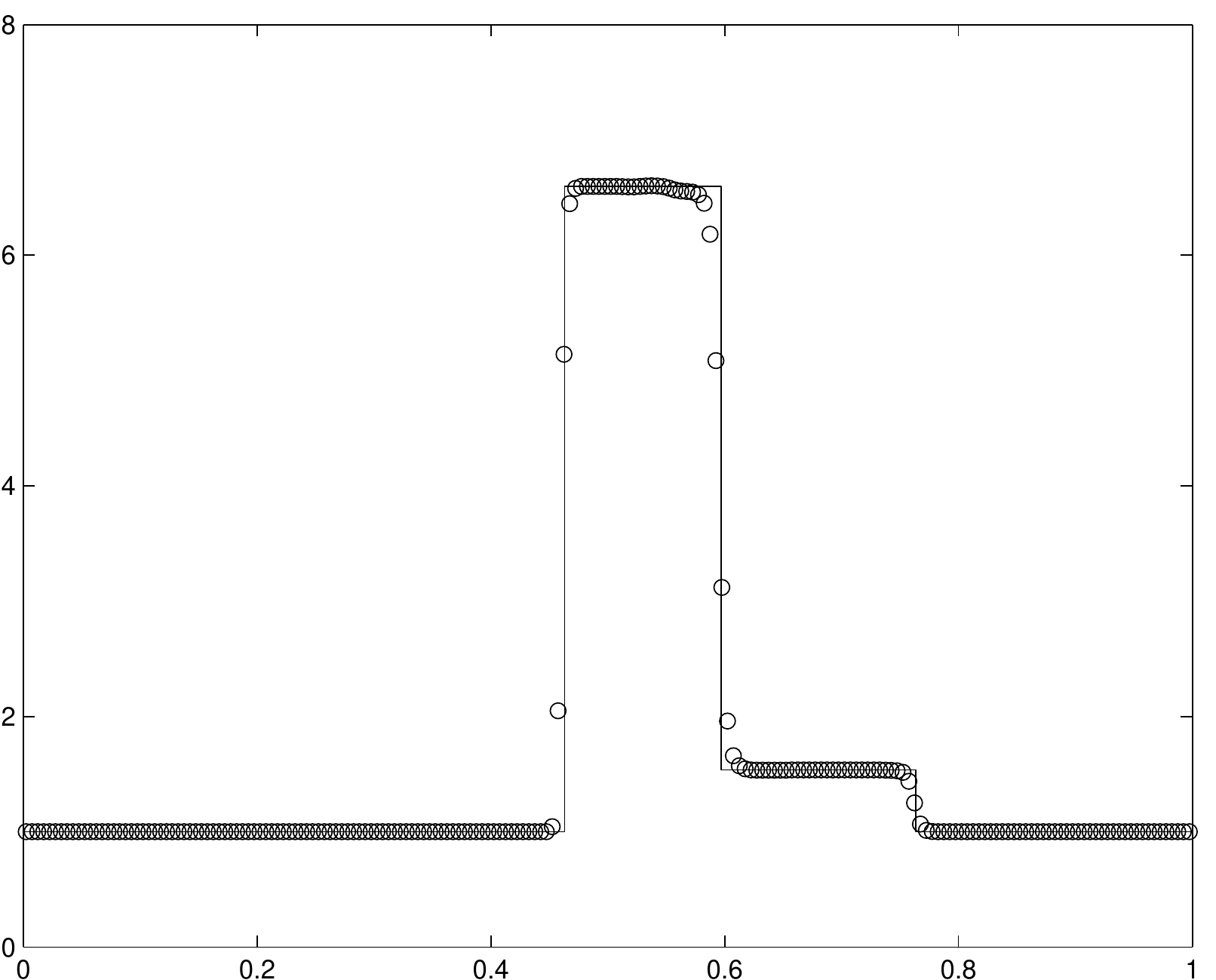}&
		\includegraphics[width=0.35\textwidth]{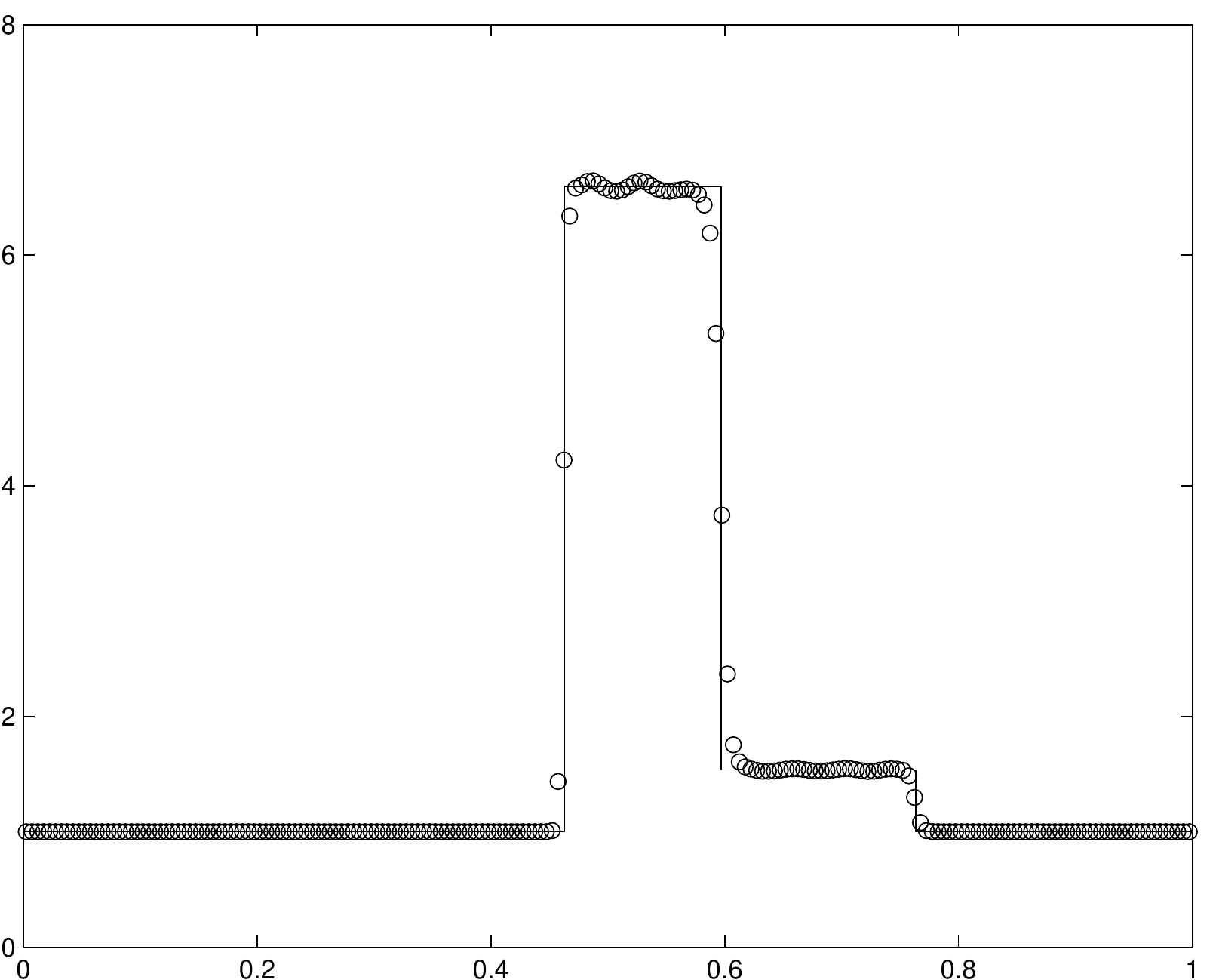}\\
	\end{tabular}
	\caption{Example~\ref{exRHDRM1DT1}: the densities $\rho$ at $t=0.4$. The symbol
		¡°$\circ$¡± denotes numerical solution with $200$ cells, while the solid line is exact solution.
Left: $P^K$-based \CDG{}; right: $P^K$-based \DG{}.  From top to bottom: $K=1,~2,~3$. }
	\label{fig:RHDRM1DT1rho}
\end{figure}

\begin{figure}[!htbp]
	\centering{}
	\begin{tabular}{cc}
		\includegraphics[width=0.35\textwidth]{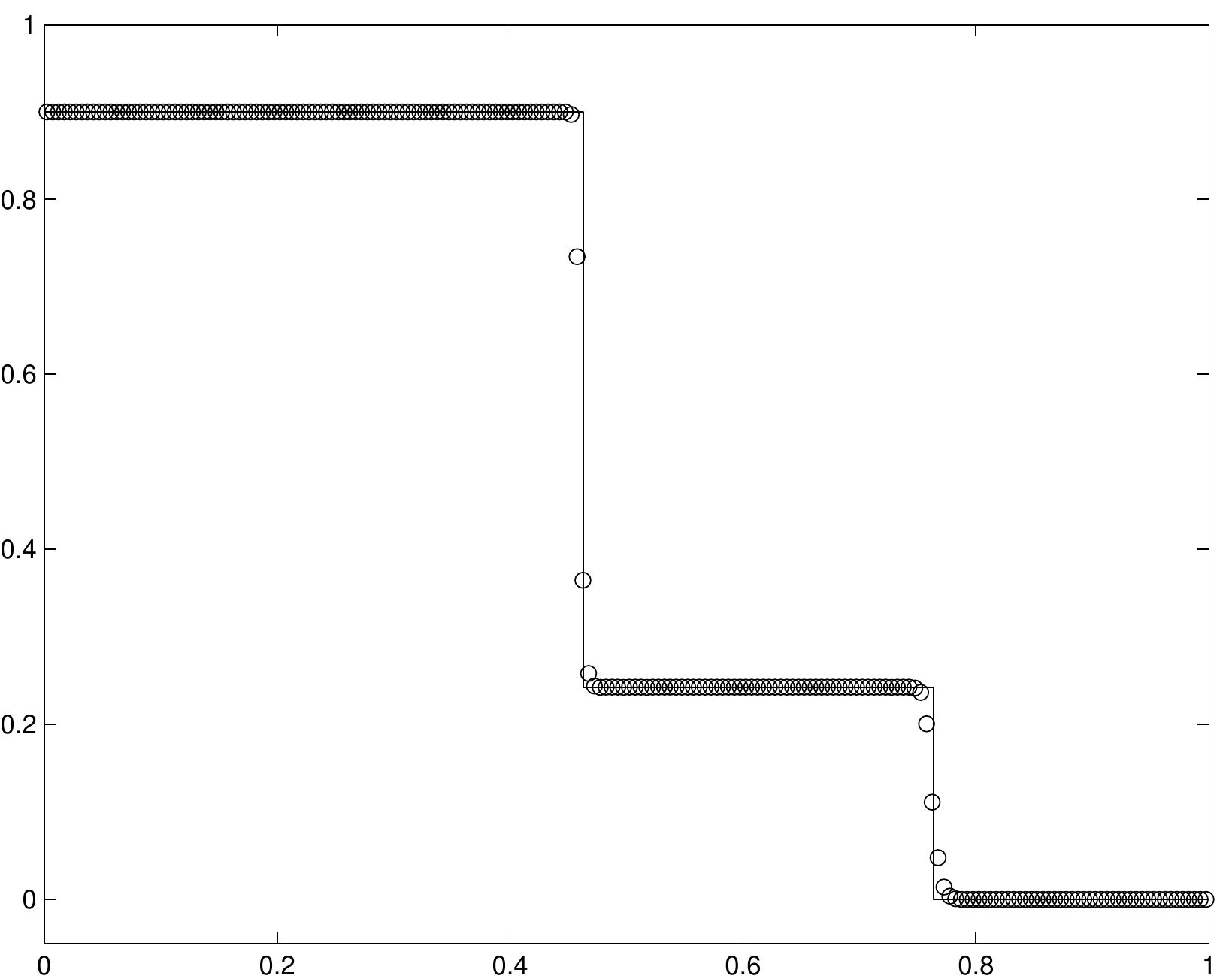}&
		\includegraphics[width=0.35\textwidth]{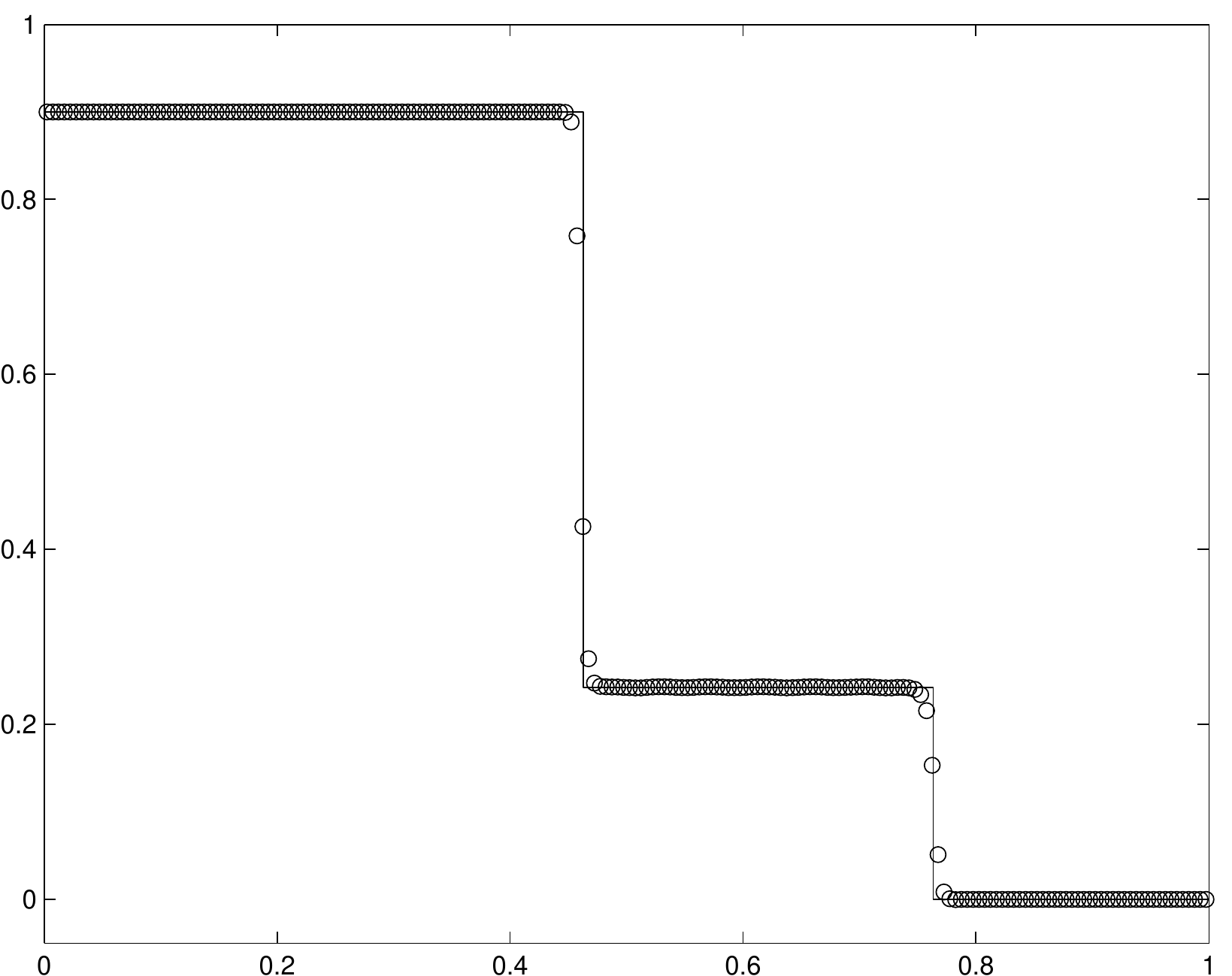}\\
		\includegraphics[width=0.35\textwidth]{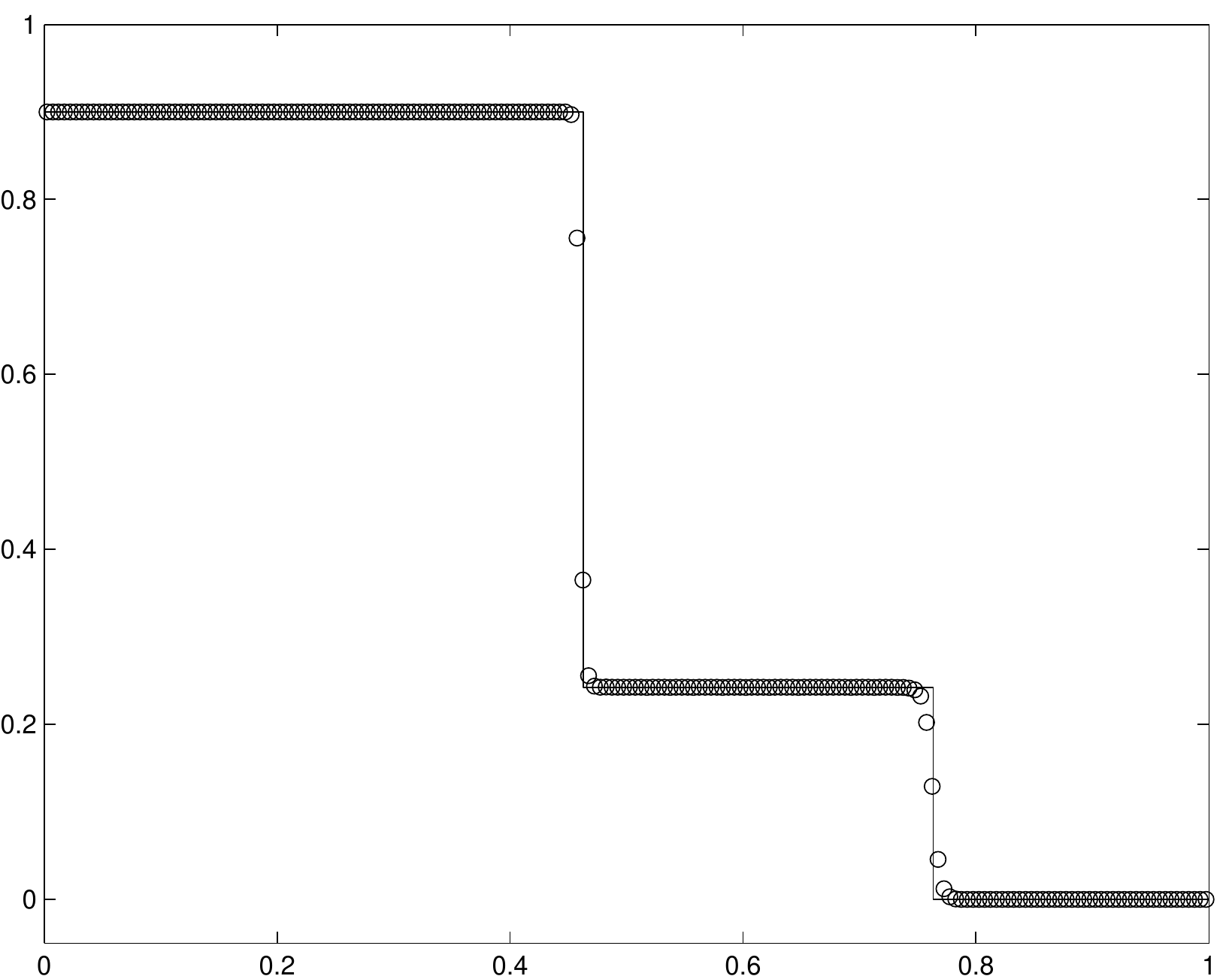}&
		\includegraphics[width=0.35\textwidth]{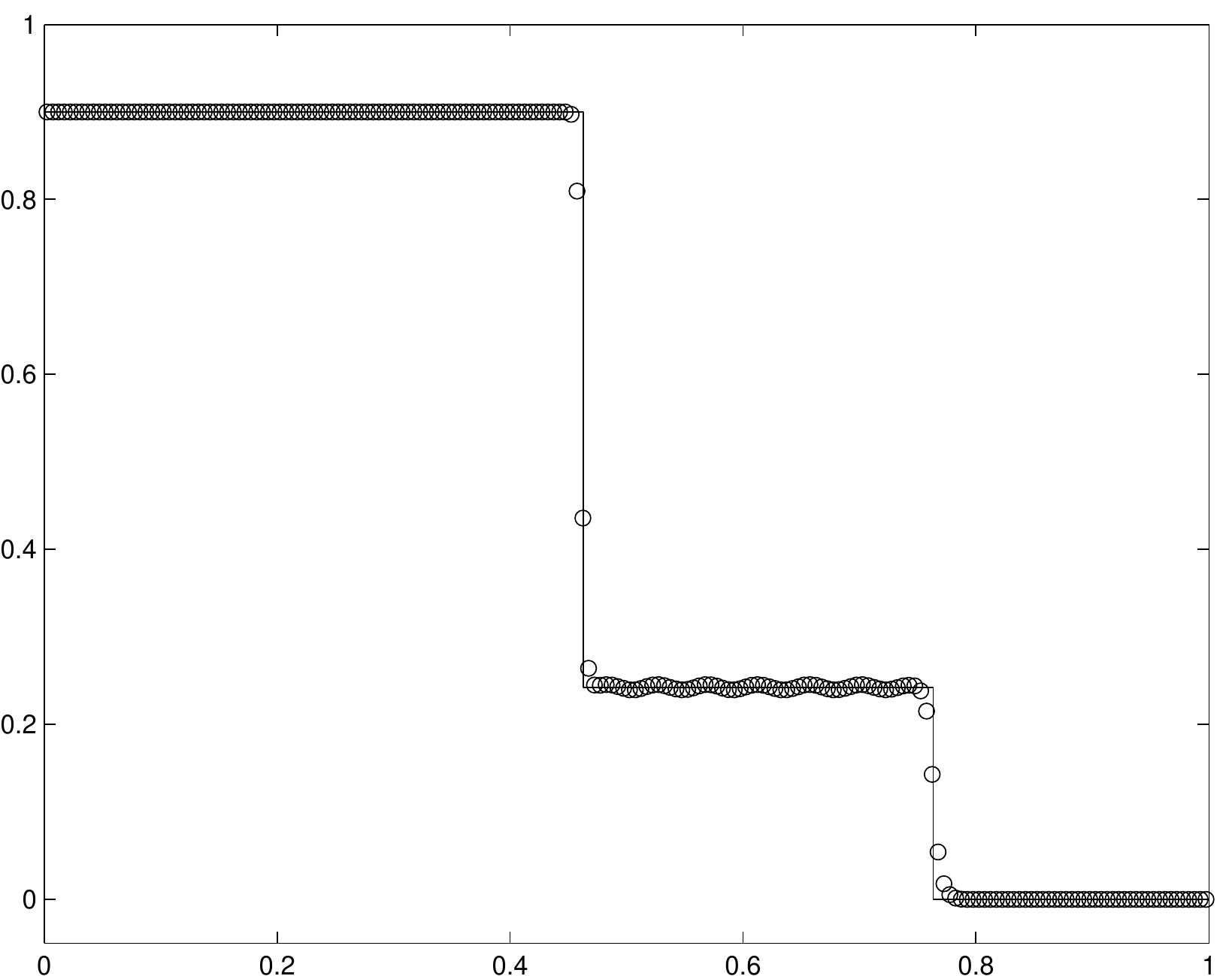}\\
		\includegraphics[width=0.35\textwidth]{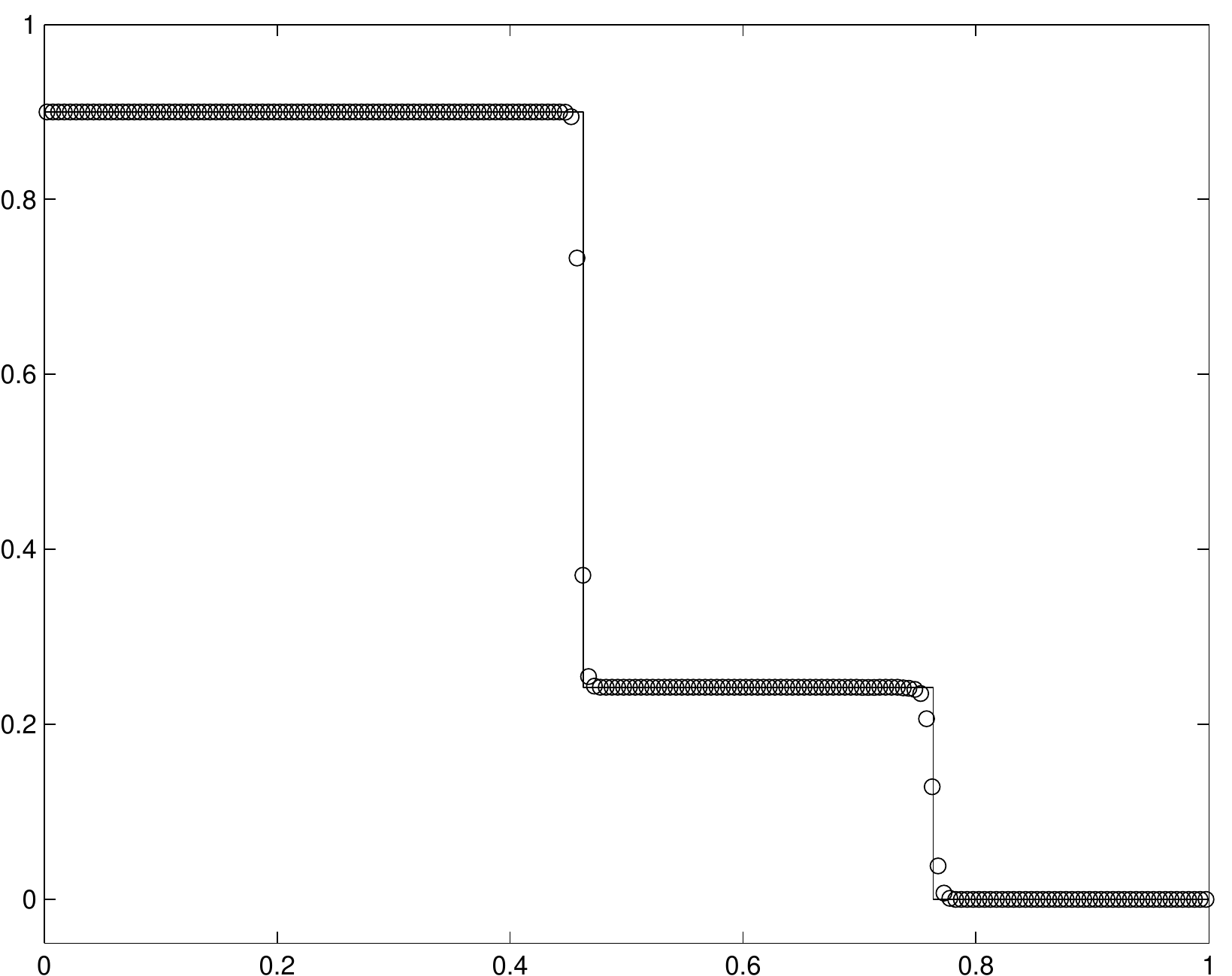}&
		\includegraphics[width=0.35\textwidth]{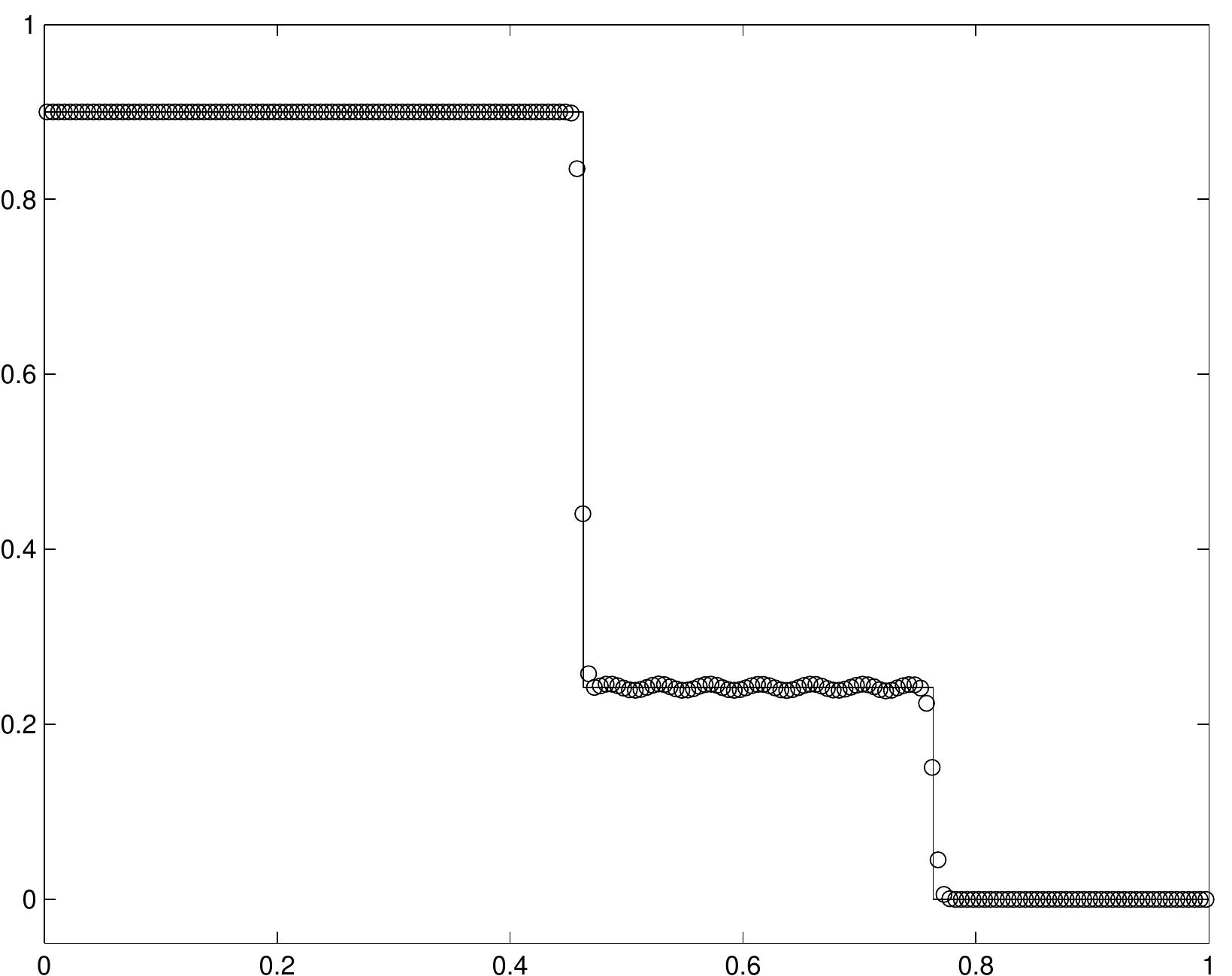}\\
	\end{tabular}
	\caption{Same as Fig. ~\ref{fig:RHDRM1DT1rho} except for the  velocity $v_1$. }
	\label{fig:RHDRM1DT1vel}
\end{figure}

\begin{figure}[!htbp]
	\centering{}
	\begin{tabular}{cc}
		\includegraphics[width=0.35\textwidth]{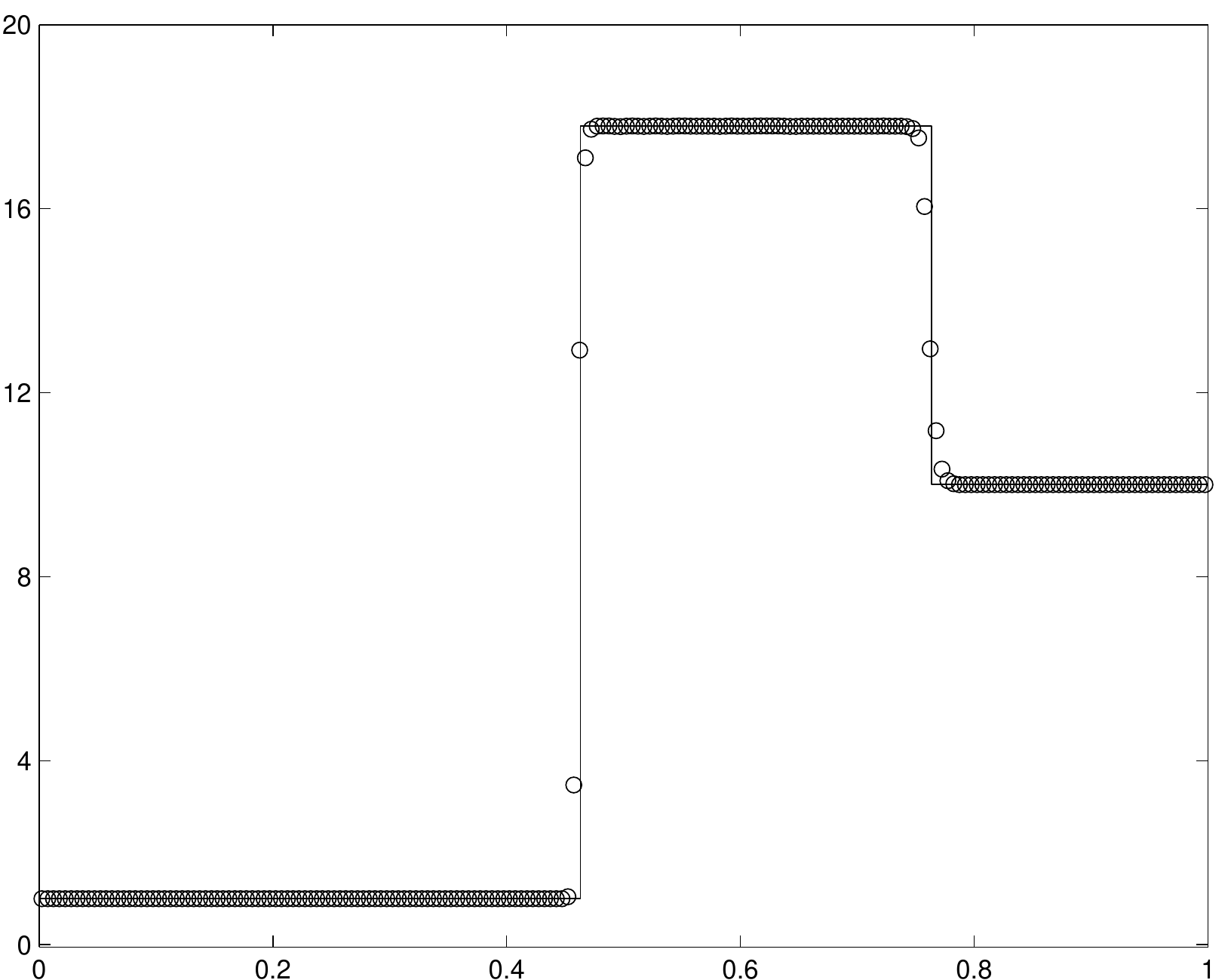}&
		\includegraphics[width=0.35\textwidth]{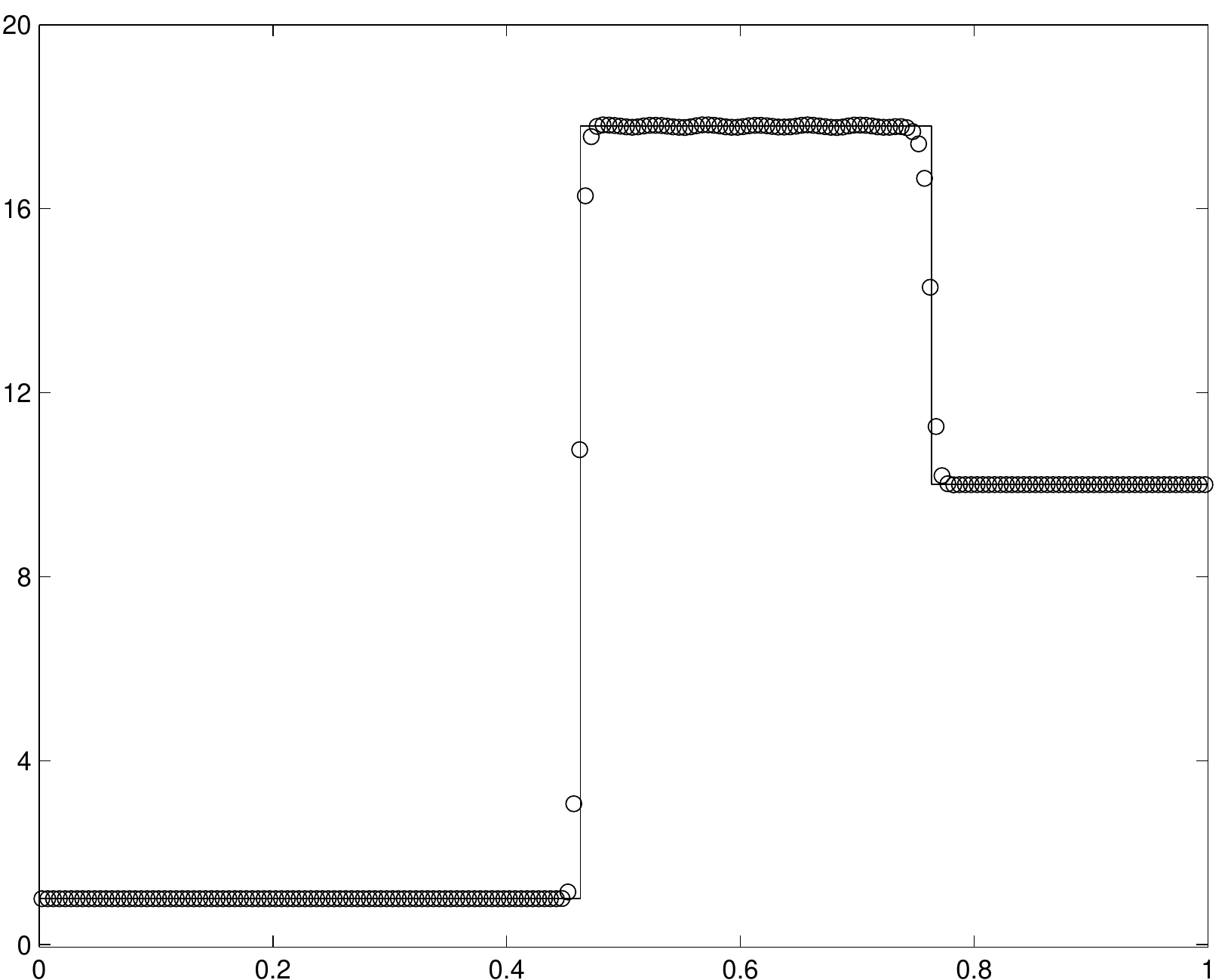}\\
		\includegraphics[width=0.35\textwidth]{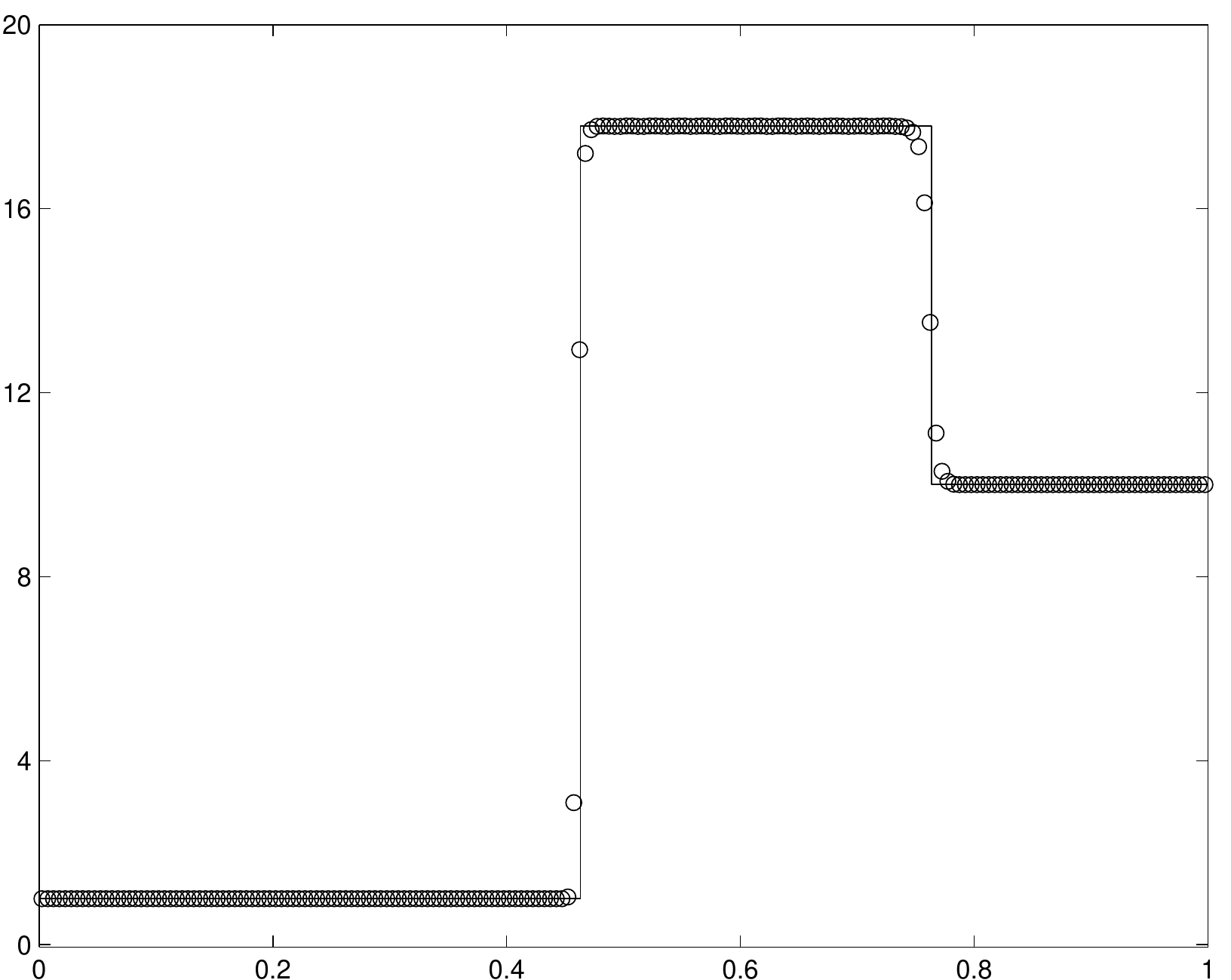}&
		\includegraphics[width=0.35\textwidth]{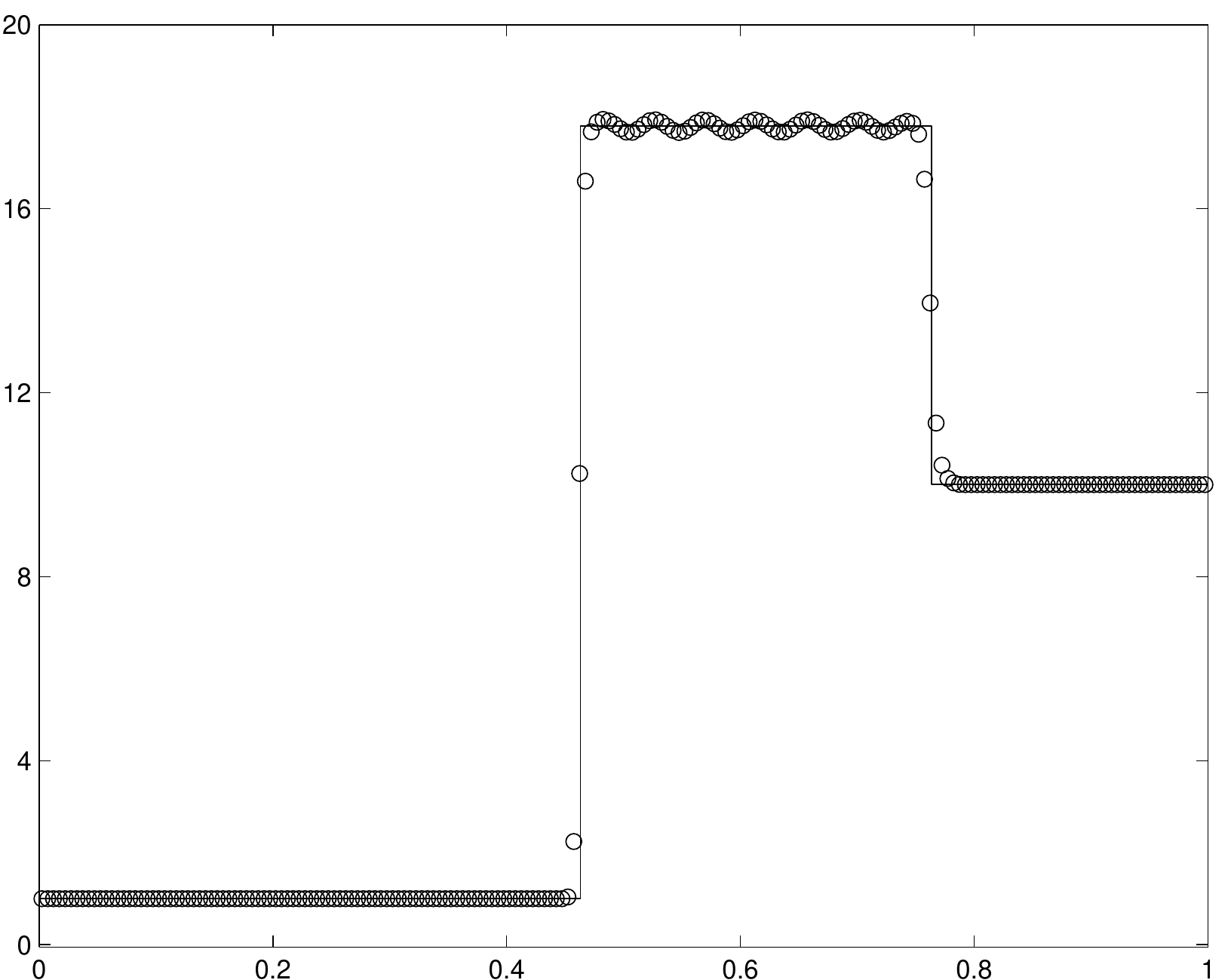}\\
		\includegraphics[width=0.35\textwidth]{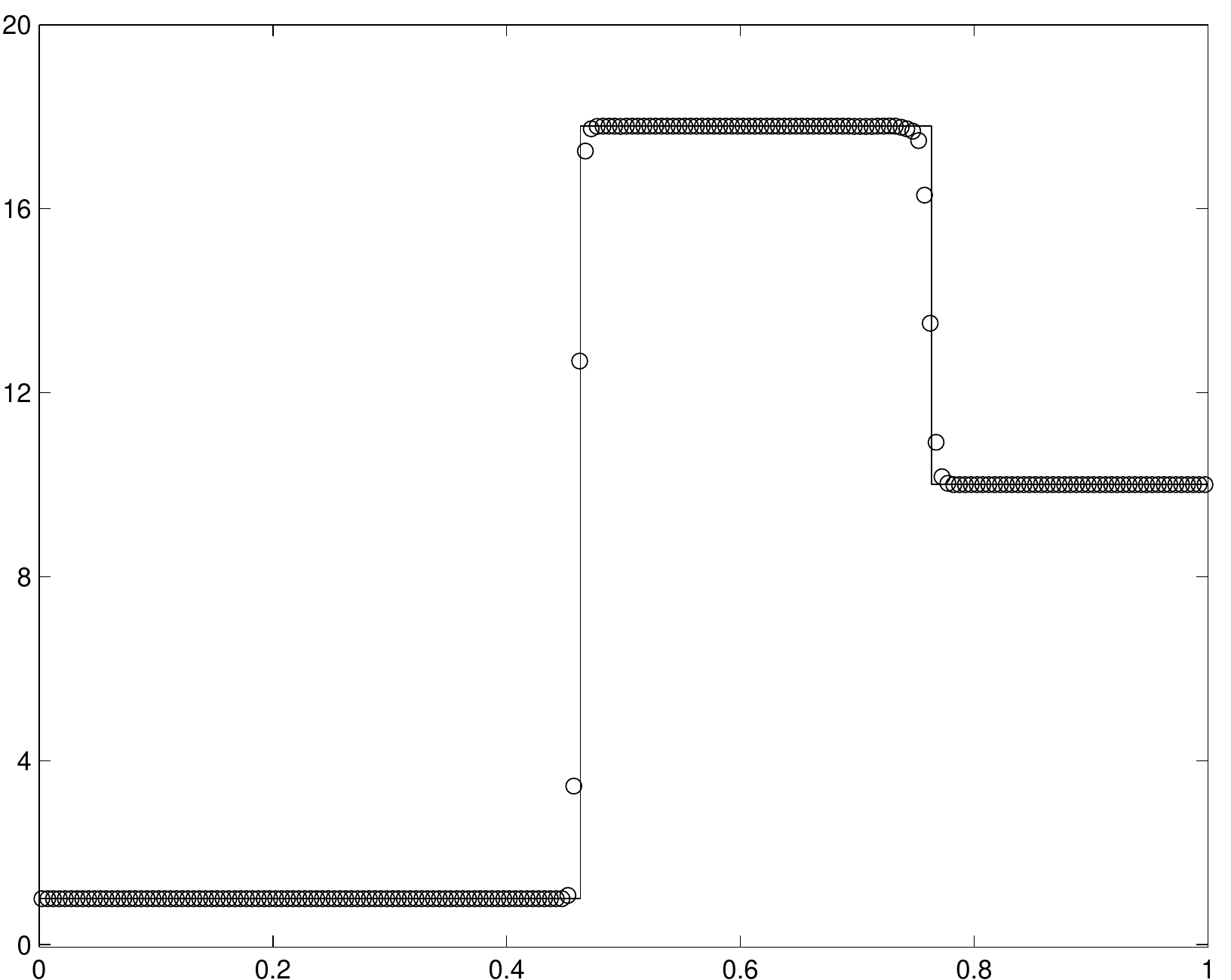}&
		\includegraphics[width=0.35\textwidth]{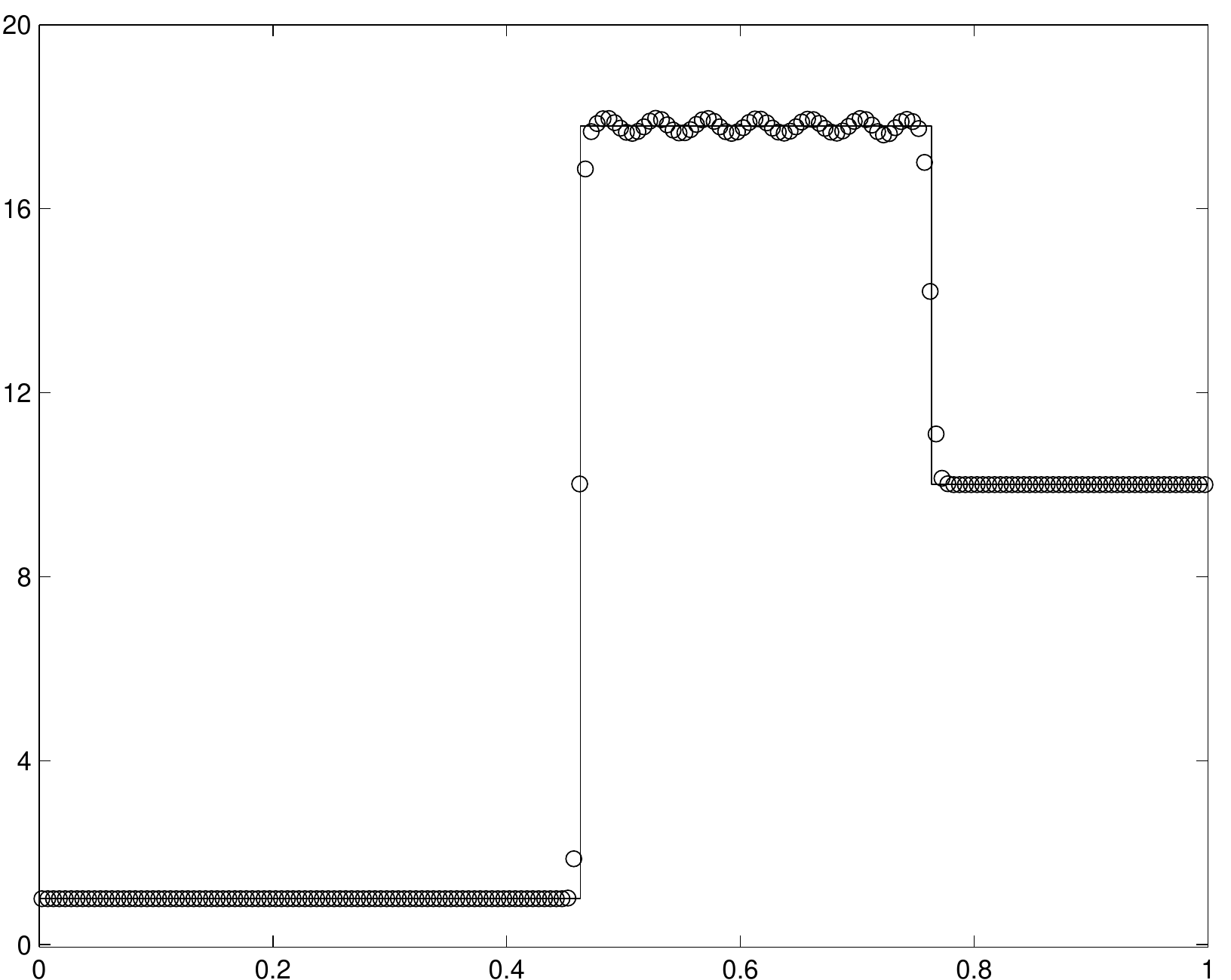}\\
	\end{tabular}
	\caption{Same as Fig. ~\ref{fig:RHDRM1DT1rho} except for the  pressure $p$. }
	\label{fig:RHDRM1DT1pre}
\end{figure}

\begin{figure}[!htbp]
	\centering{}
	\begin{tabular}{cc}
		\includegraphics[width=0.35\textwidth]{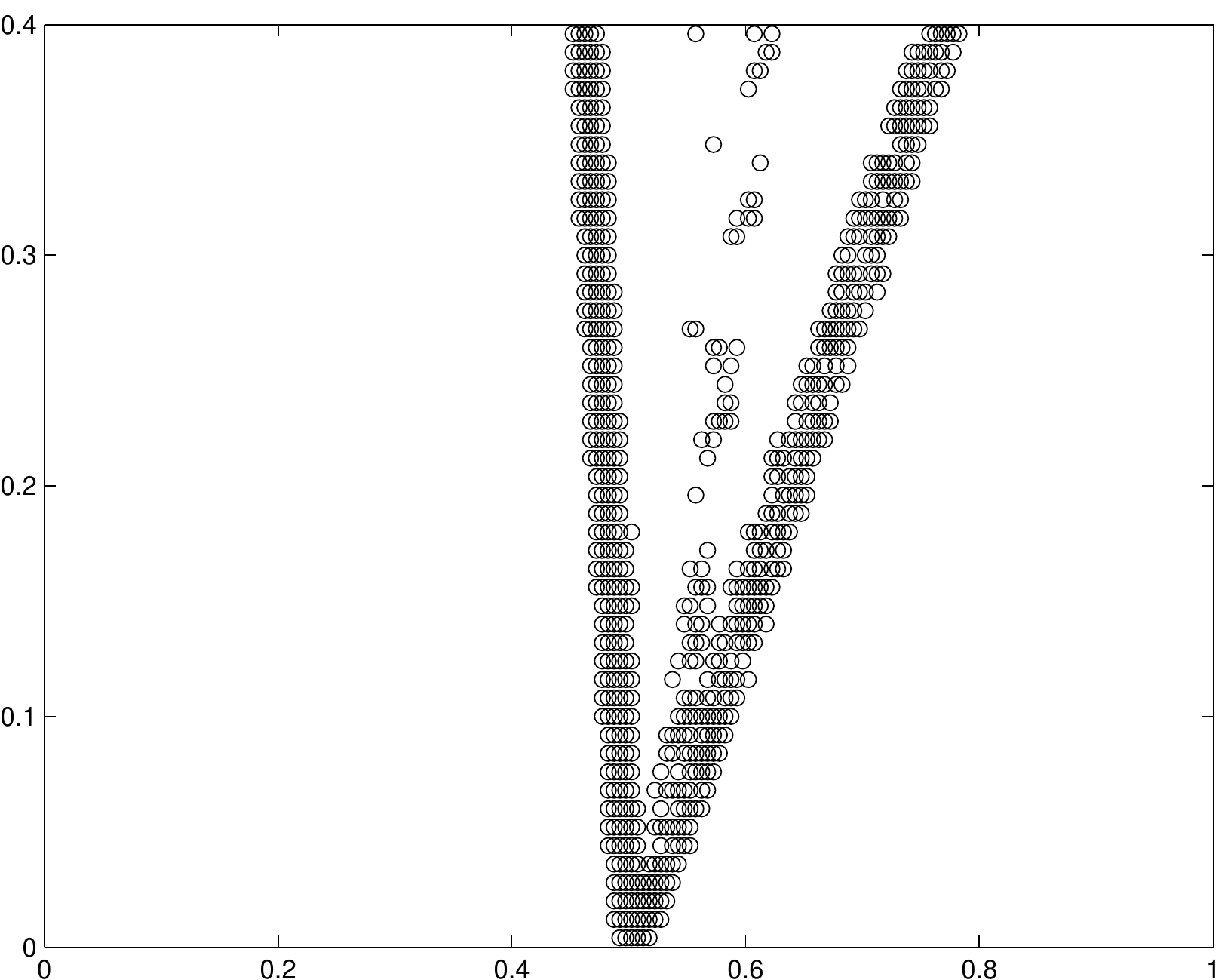}&
		\includegraphics[width=0.35\textwidth]{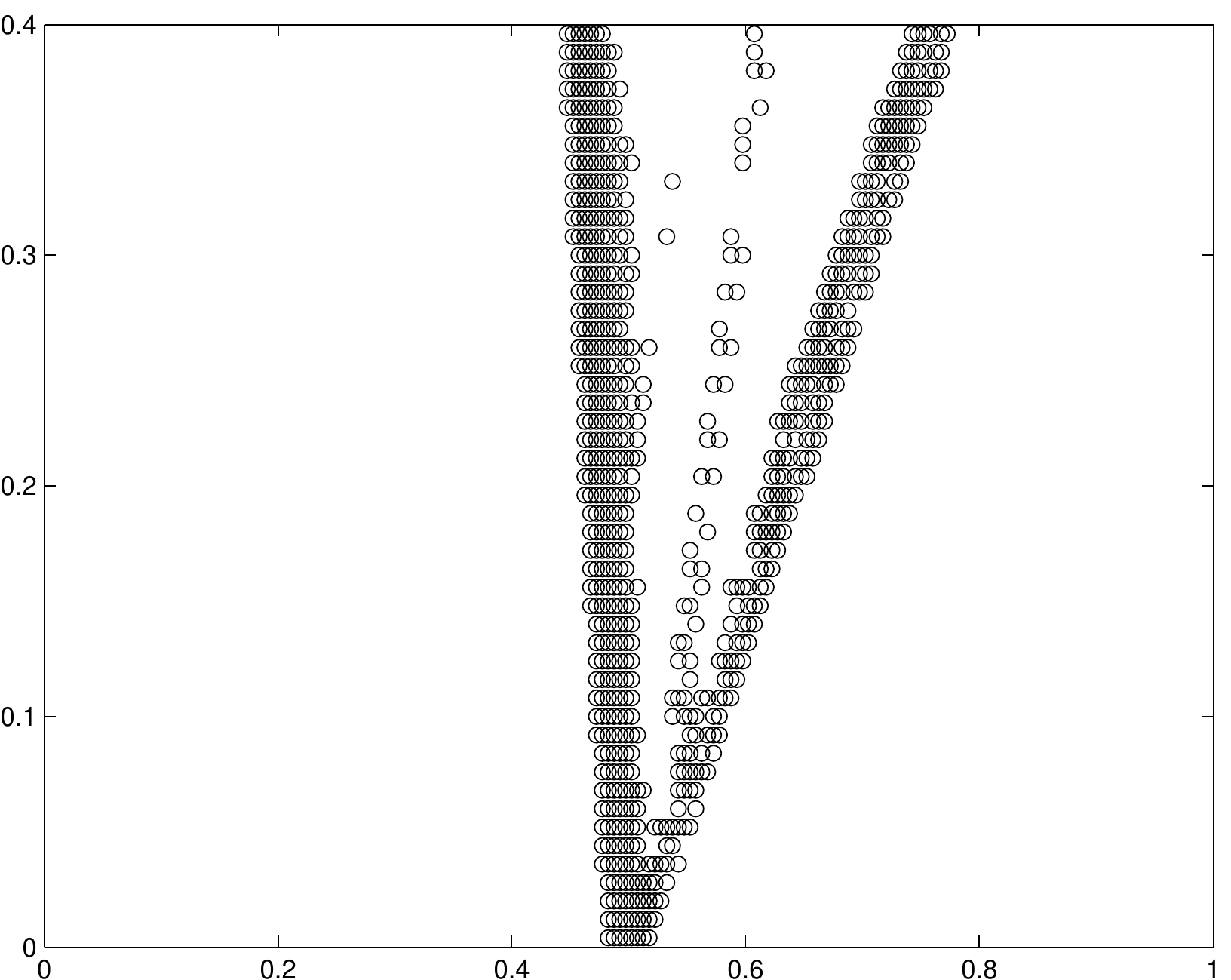}\\
		\includegraphics[width=0.35\textwidth]{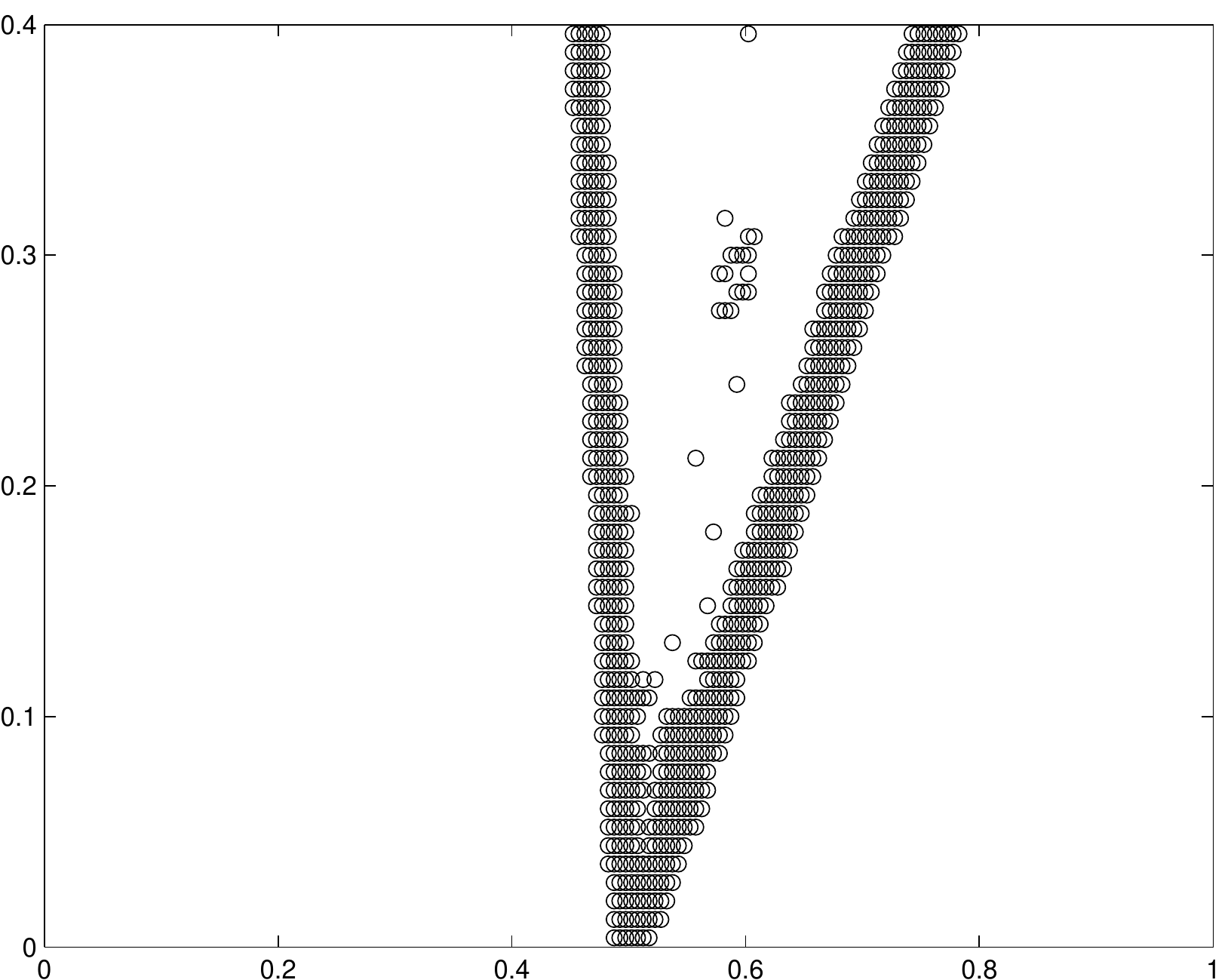}&
		\includegraphics[width=0.35\textwidth]{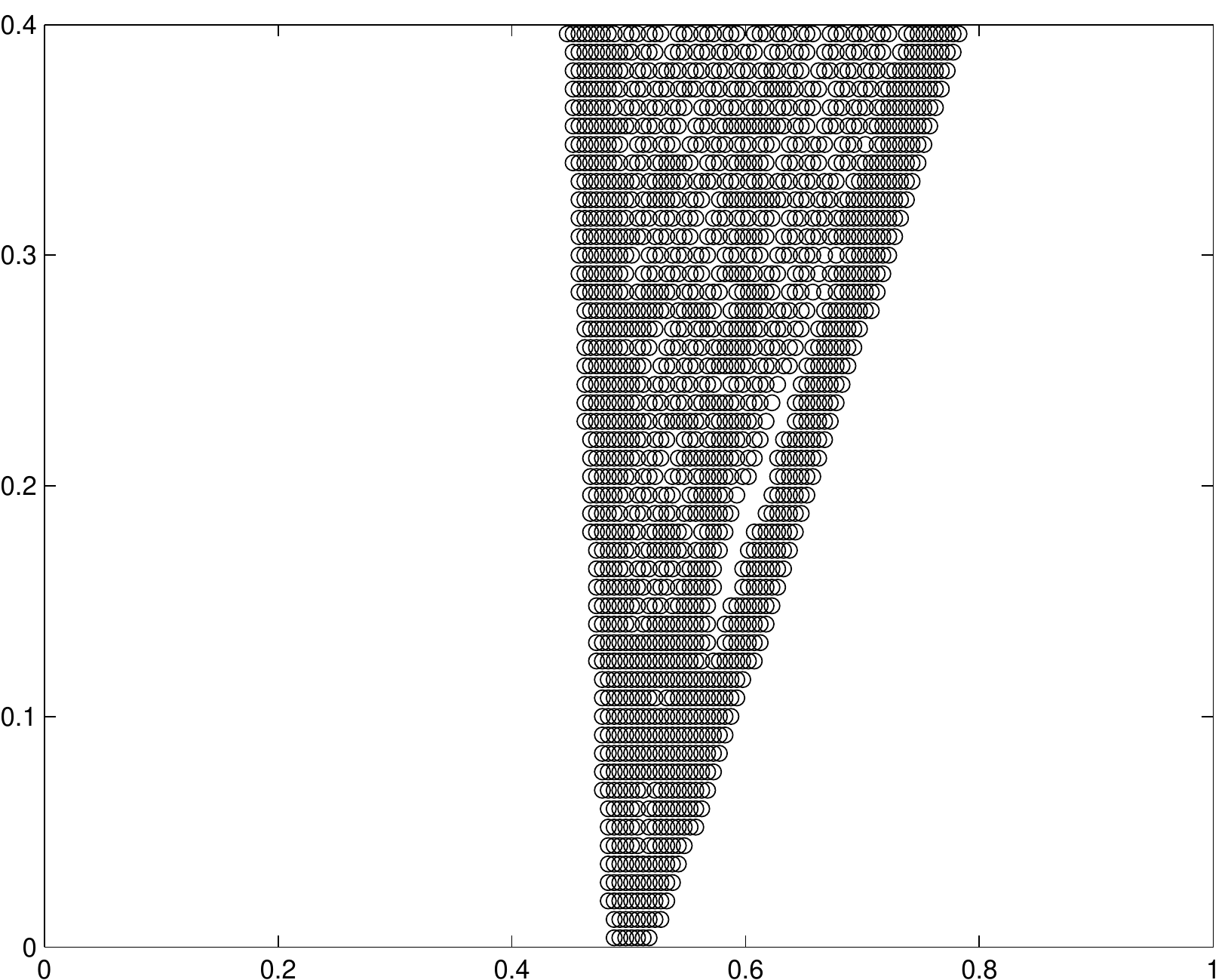}\\
		\includegraphics[width=0.35\textwidth]{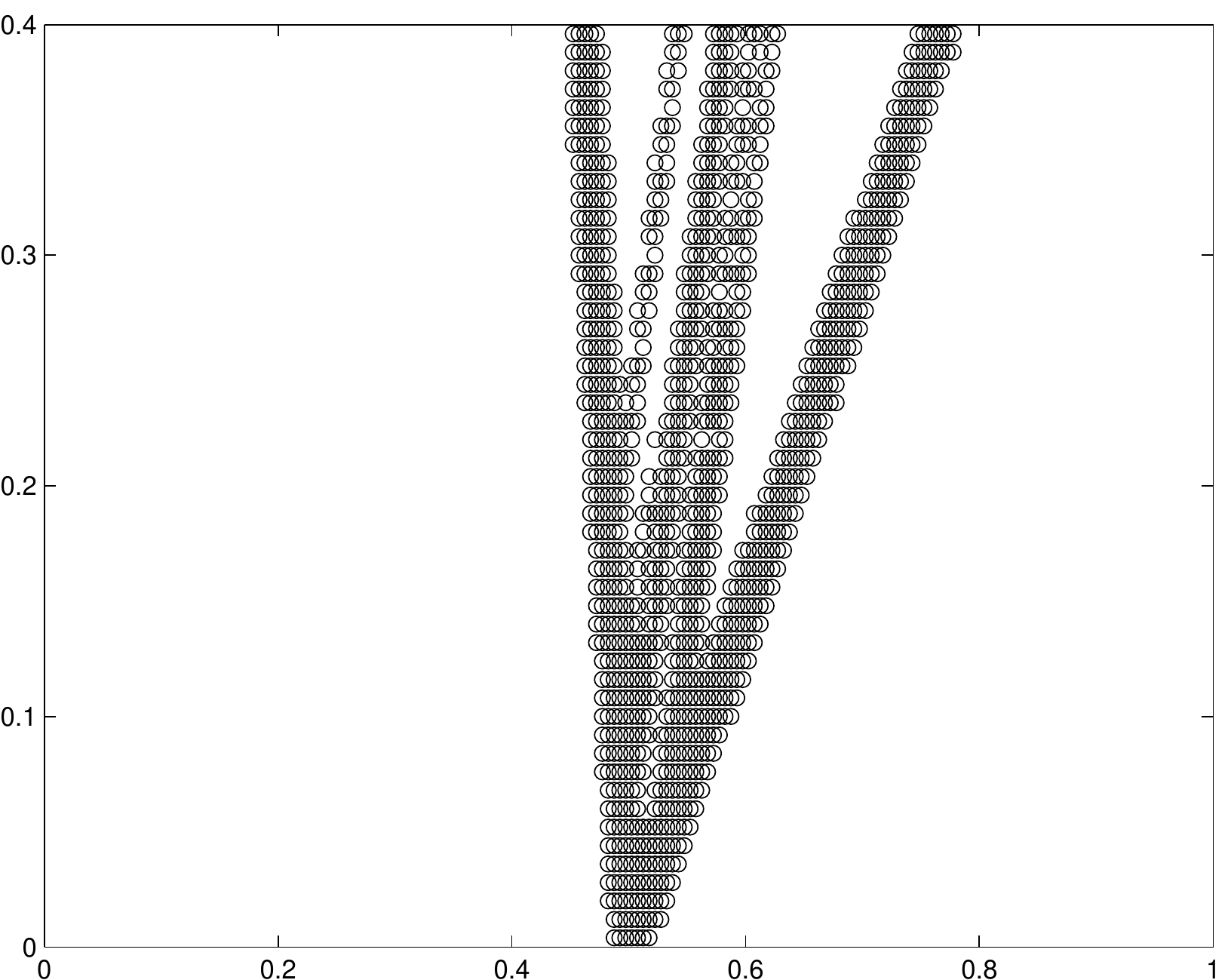}&
		\includegraphics[width=0.35\textwidth]{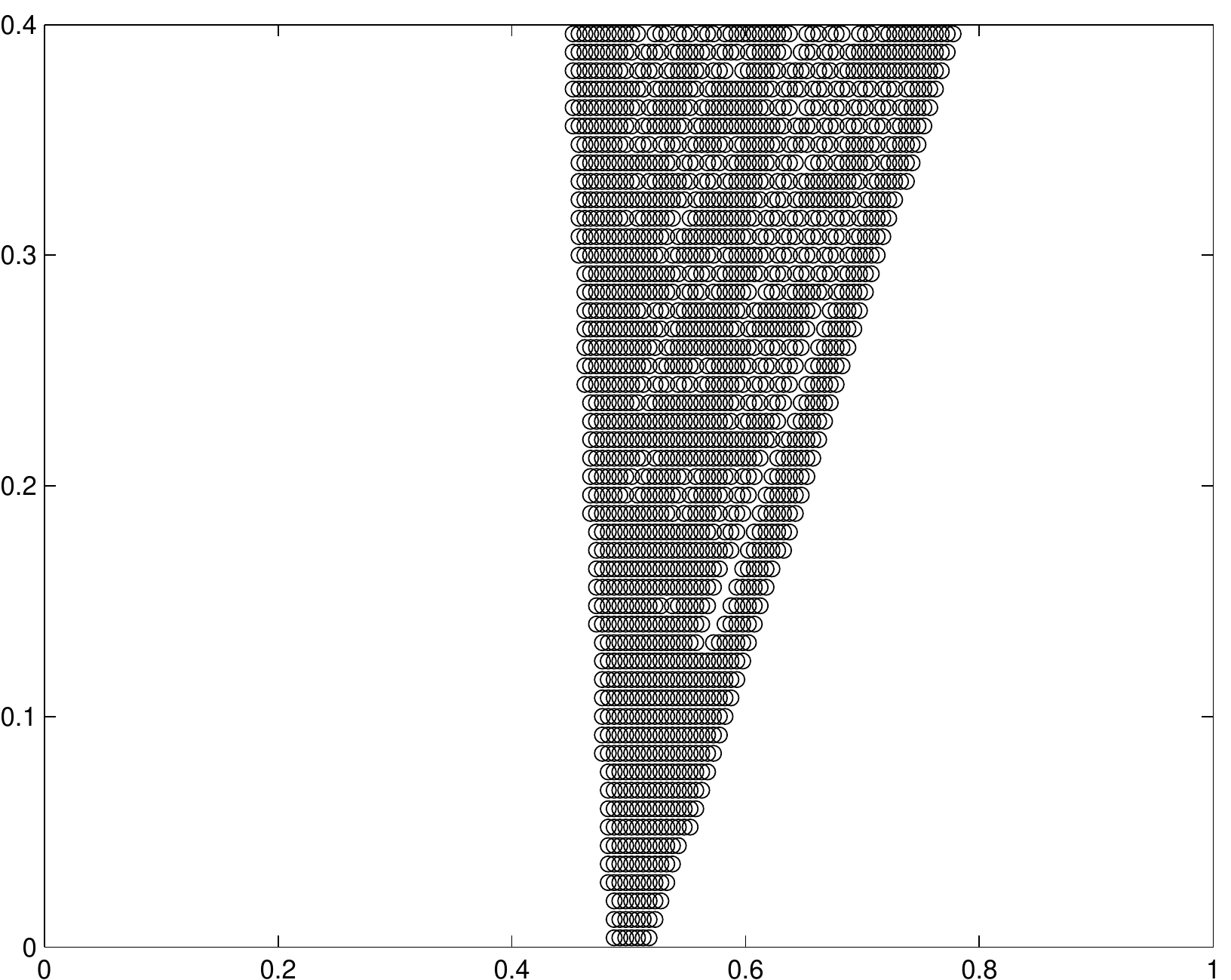}\\
		
	\end{tabular}
	\caption{
Same as Fig.~\ref{fig:RHDRM1DT1rho} except for the ``troubled'' cells in the $(x,t)$ plane. }
	\label{fig:RHDRM1DT1cell}
\end{figure}
\fi


\begin{Example}[Riemann problem 2]\label{exRHDRM1DT4}\rm
	The initial data of second Riemann problem is
	$$(\rho,v_1,p)(x,0)=\begin{cases} (10,0.0,1000), & x<0.5,
	\\ (1,0.0,0.01), &x>0.5,\end{cases}$$
and $\Gamma=5/3$. The solution of this problem will contain  a left-moving rarefaction wave, a contact discontinuity, and a right-moving shock wave as $t>0$. The speed of the contact is almost identical to
 the shock wave so that it is much more challenging for the numerical methods than the first Riemann problem. 

The \CDG{} with new and old calculations of the flux integral over the cell
are considered here.
The maximum CFL numbers are taken as those in Table~\ref{tab:cfl}.
Figs.~\ref{fig:RHDRM1DT4rhovel} and \ref{fig:RHDRM1DT4preene}
display the solutions at $t=0.4$ obtained by new and old \CDG{} with $200$ cells.
The density obtained by the old \CDG{} is slightly better than the new.
The CPU times for them with $800$ cells
are estimated in Table~\ref{tab:RHDRM1DT4cmpRC}.
The data show that the new $P^3$-based CDG method is   faster than the old,
there is no obvious difference between two $P^1$-based CDG method,
but the new $P^2$-based  CDG method is   slower than the old due to a
relatively harsh stability condition for the new $P^2$-based CDG method.
\end{Example}

\ifx\outnofig\undefine
\begin{figure}[!htbp]
	\centering{}
	\begin{tabular}{cc}
		\includegraphics[width=0.35\textwidth]{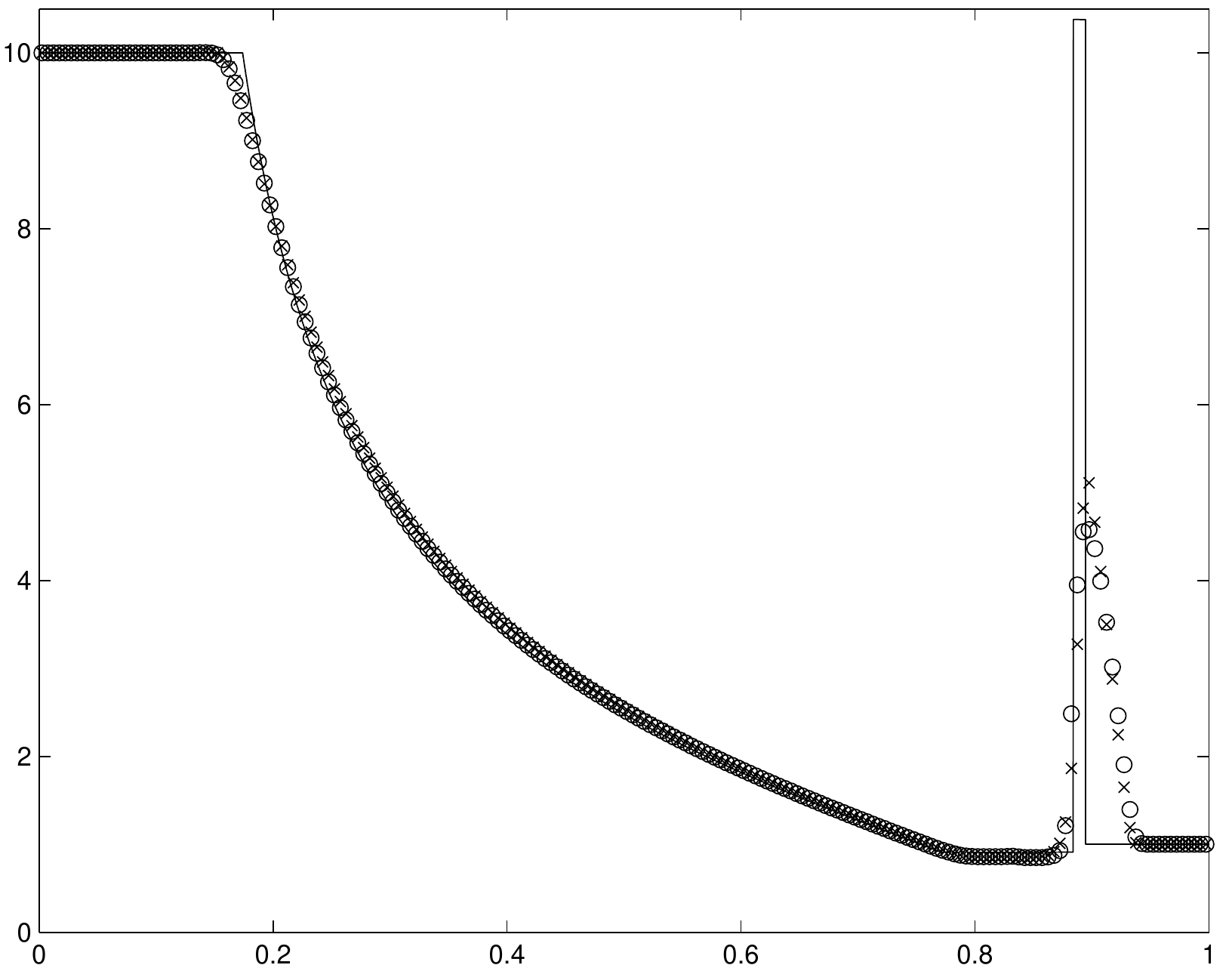}&
		\includegraphics[width=0.35\textwidth]{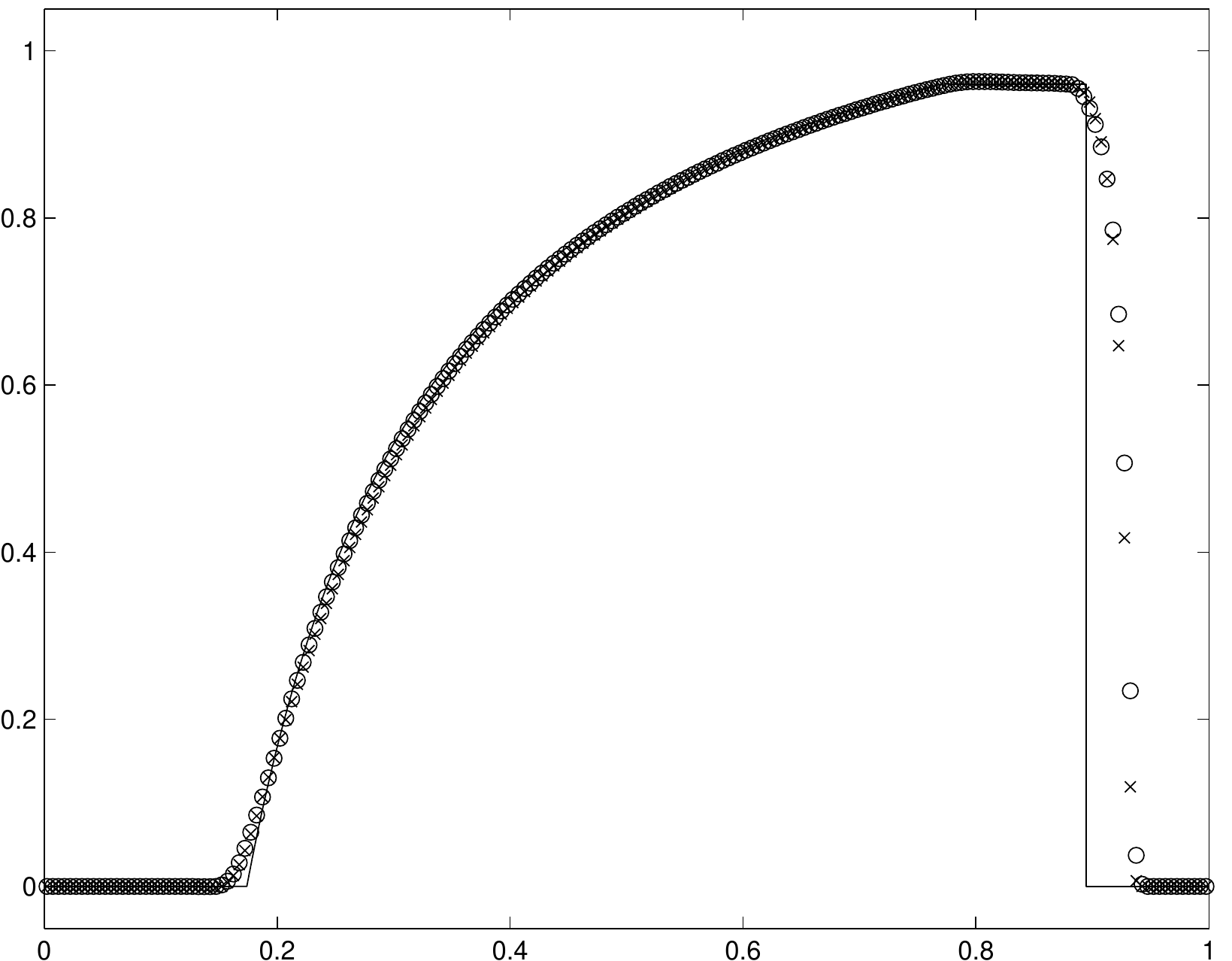}\\
		\includegraphics[width=0.35\textwidth]{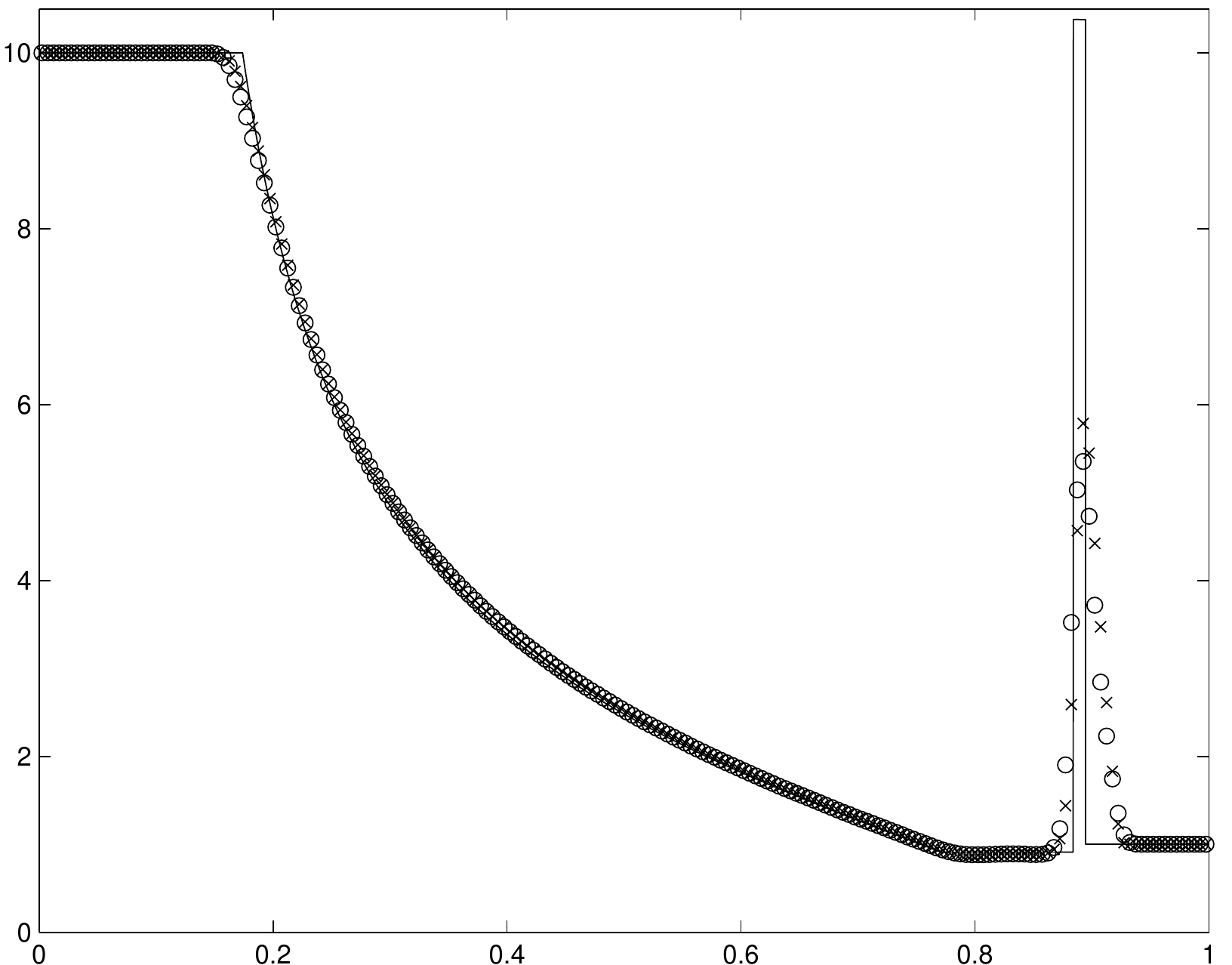}&
		\includegraphics[width=0.35\textwidth]{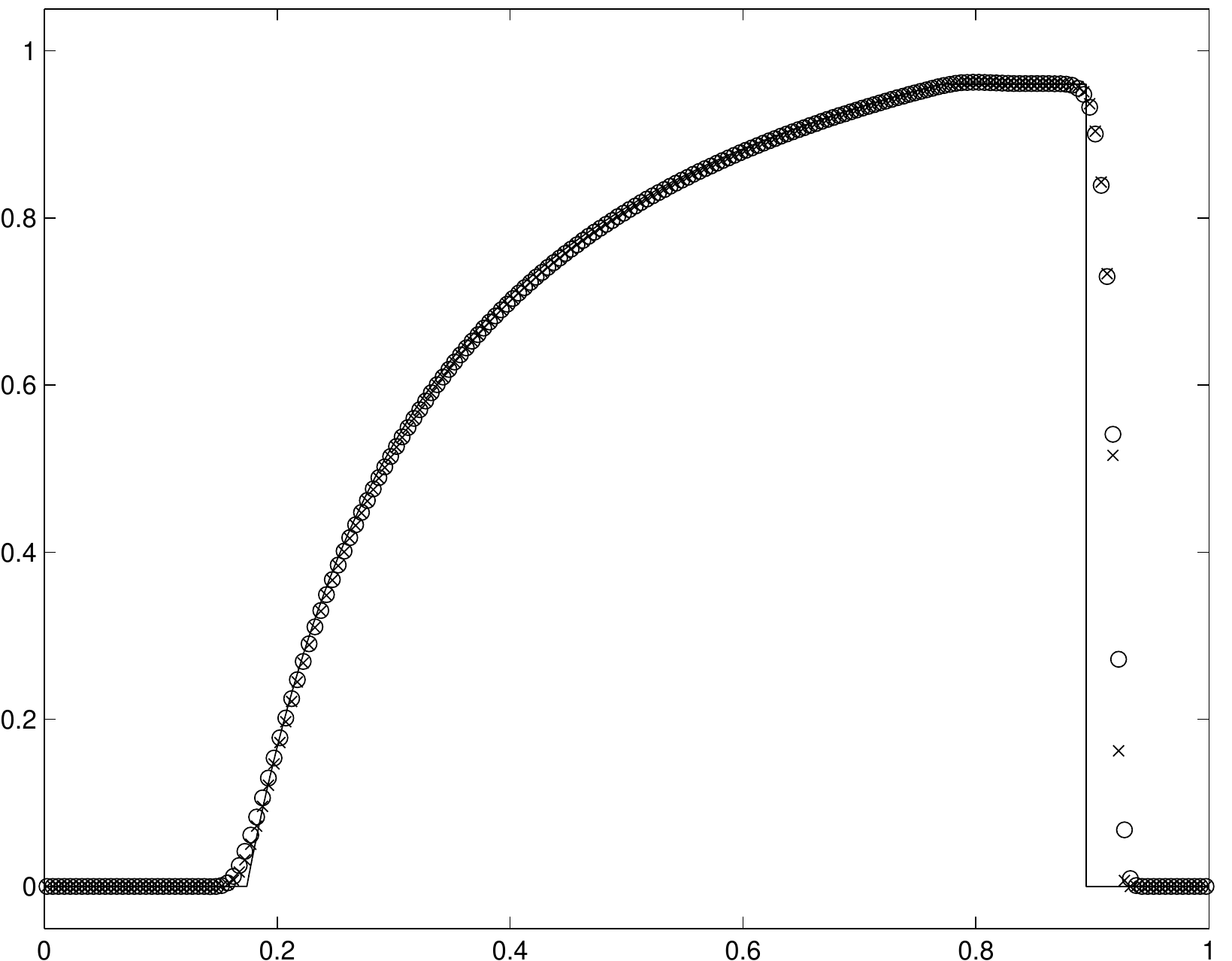}\\
		\includegraphics[width=0.35\textwidth]{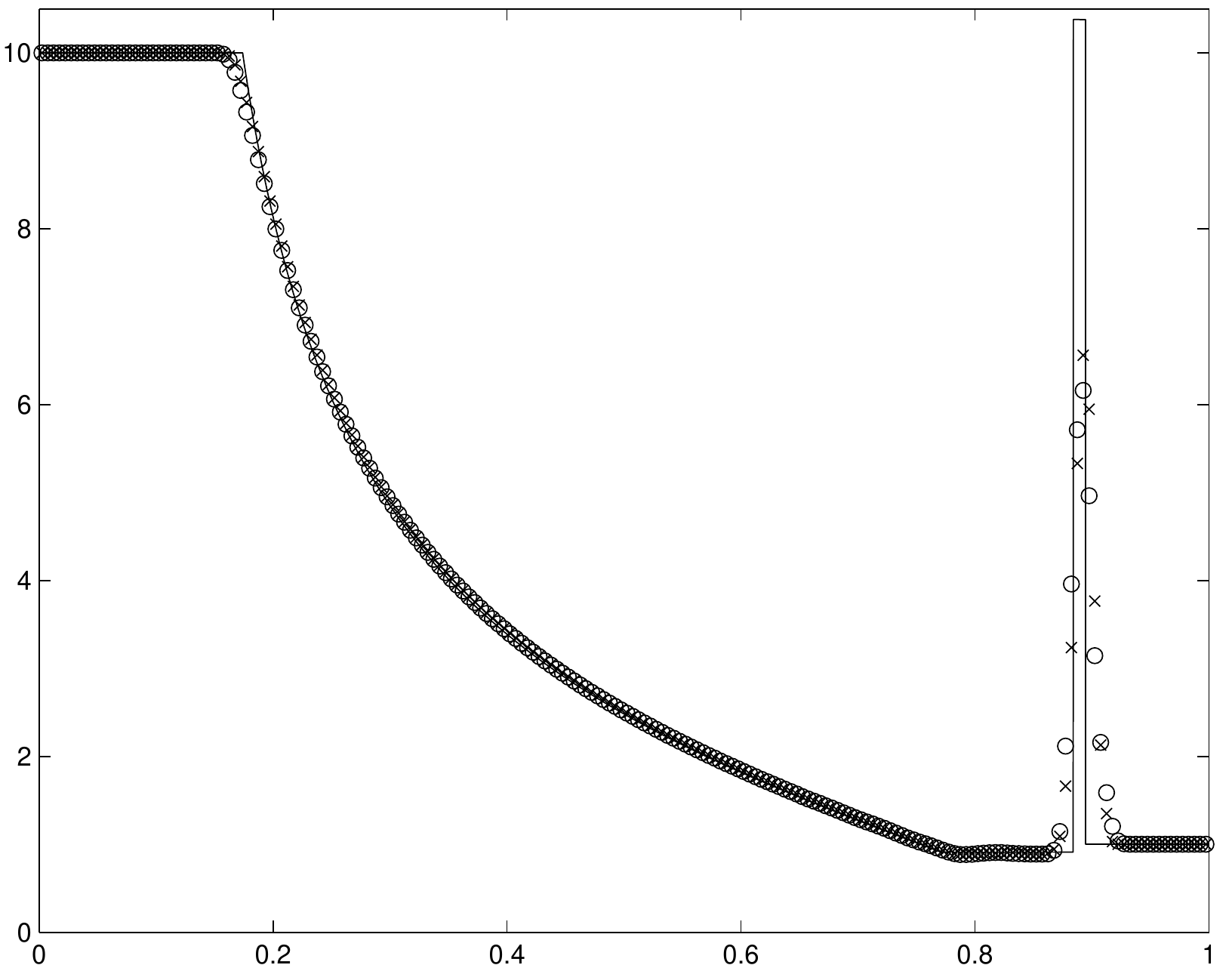}&
		\includegraphics[width=0.35\textwidth]{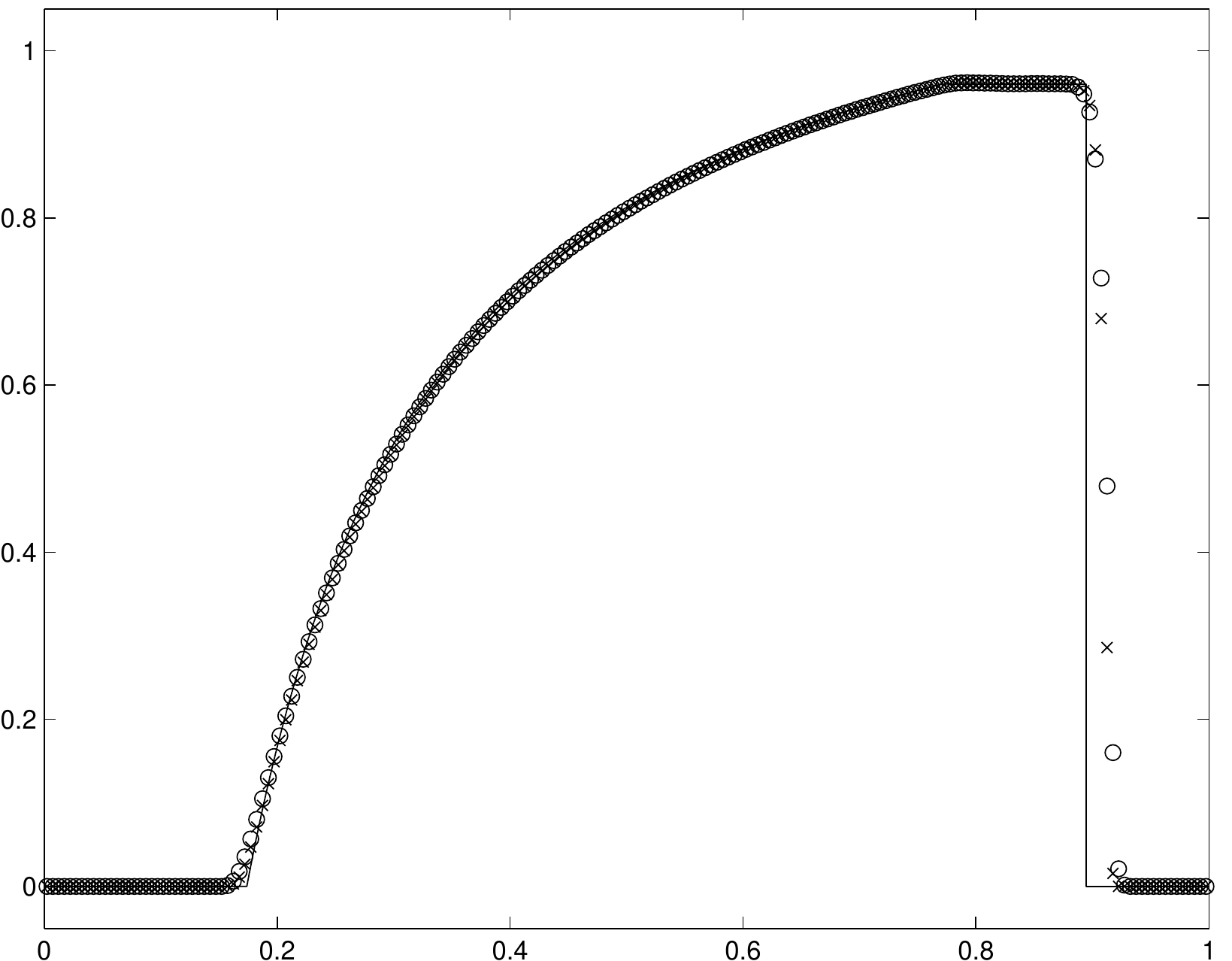}\\
	\end{tabular}
	\caption{Example \ref{exRHDRM1DT4}: The densities  $\rho$ (left) and velocities
		$v_1$ (right) at $t=0.4$. The solid line denotes the exact solution, while the symbol
``$\circ$''
and ``$\times$'' are the solutions obtained by using the new and old \CDG{} with $200$ uniform cells.
From top to bottom:
$K=1,~2,~3$.  }
	\label{fig:RHDRM1DT4rhovel}
\end{figure}

\begin{figure}[!htbp]
	\centering{}
	\begin{tabular}{cc}
		\includegraphics[width=0.35\textwidth]{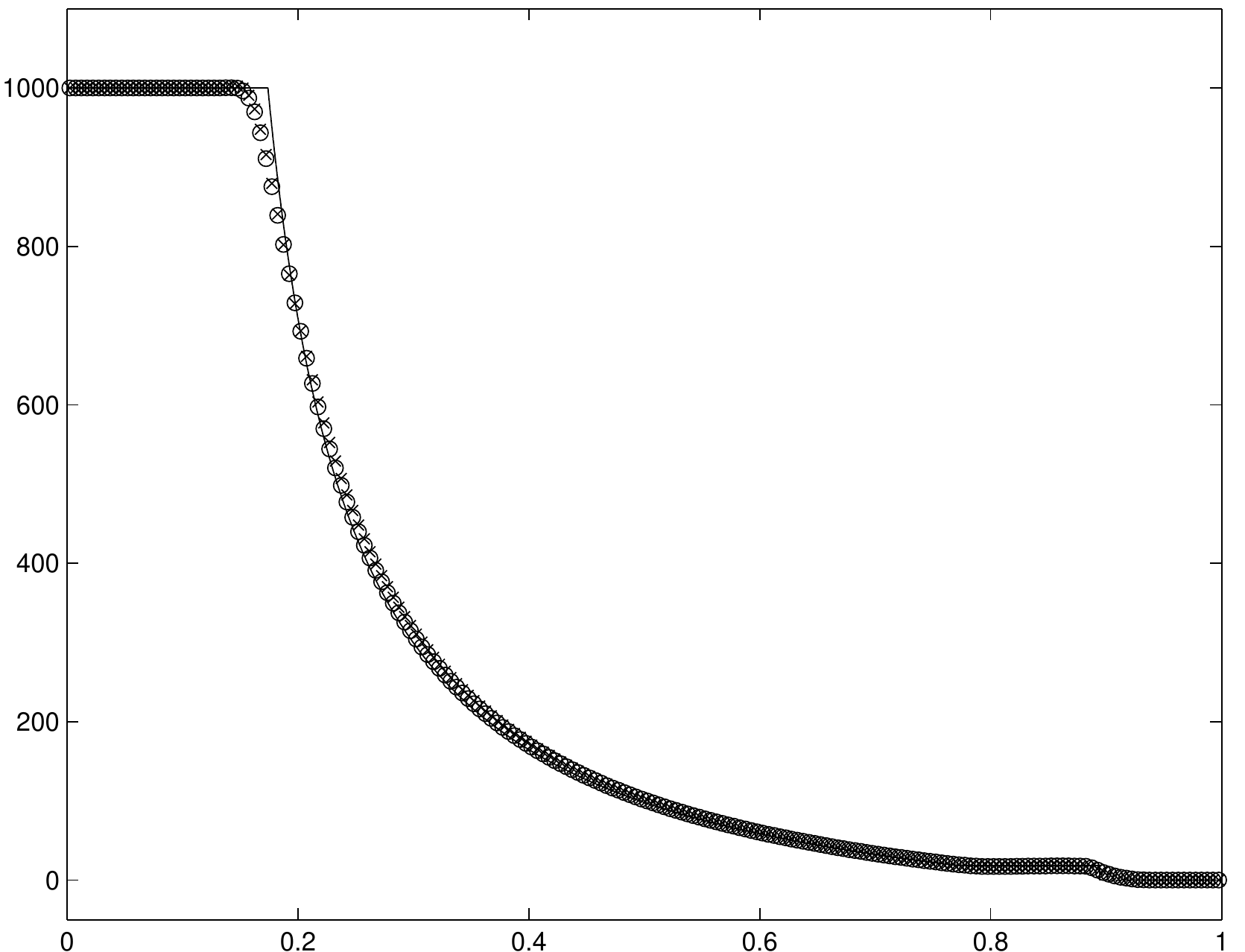}&
		\includegraphics[width=0.35\textwidth]{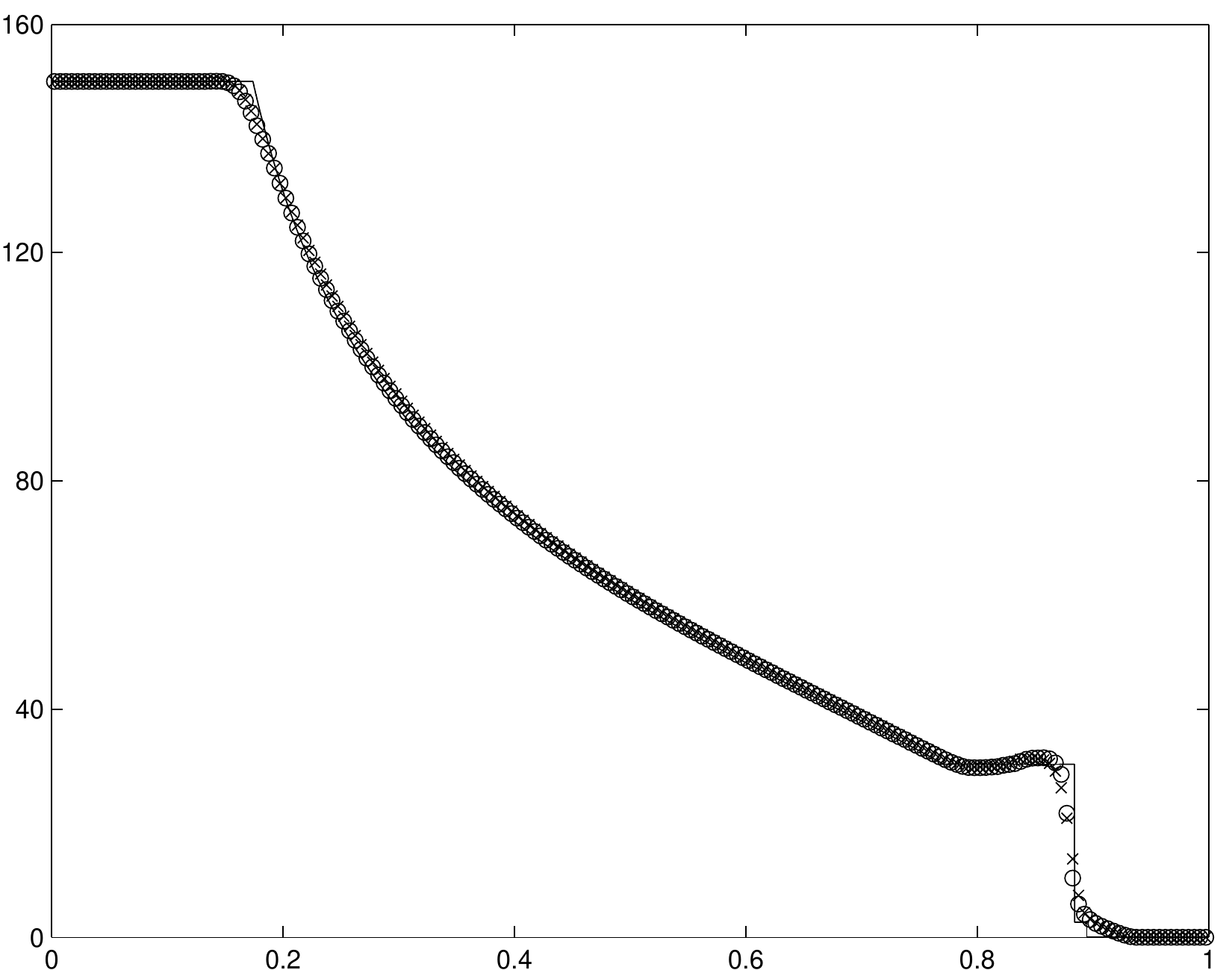}\\
		\includegraphics[width=0.35\textwidth]{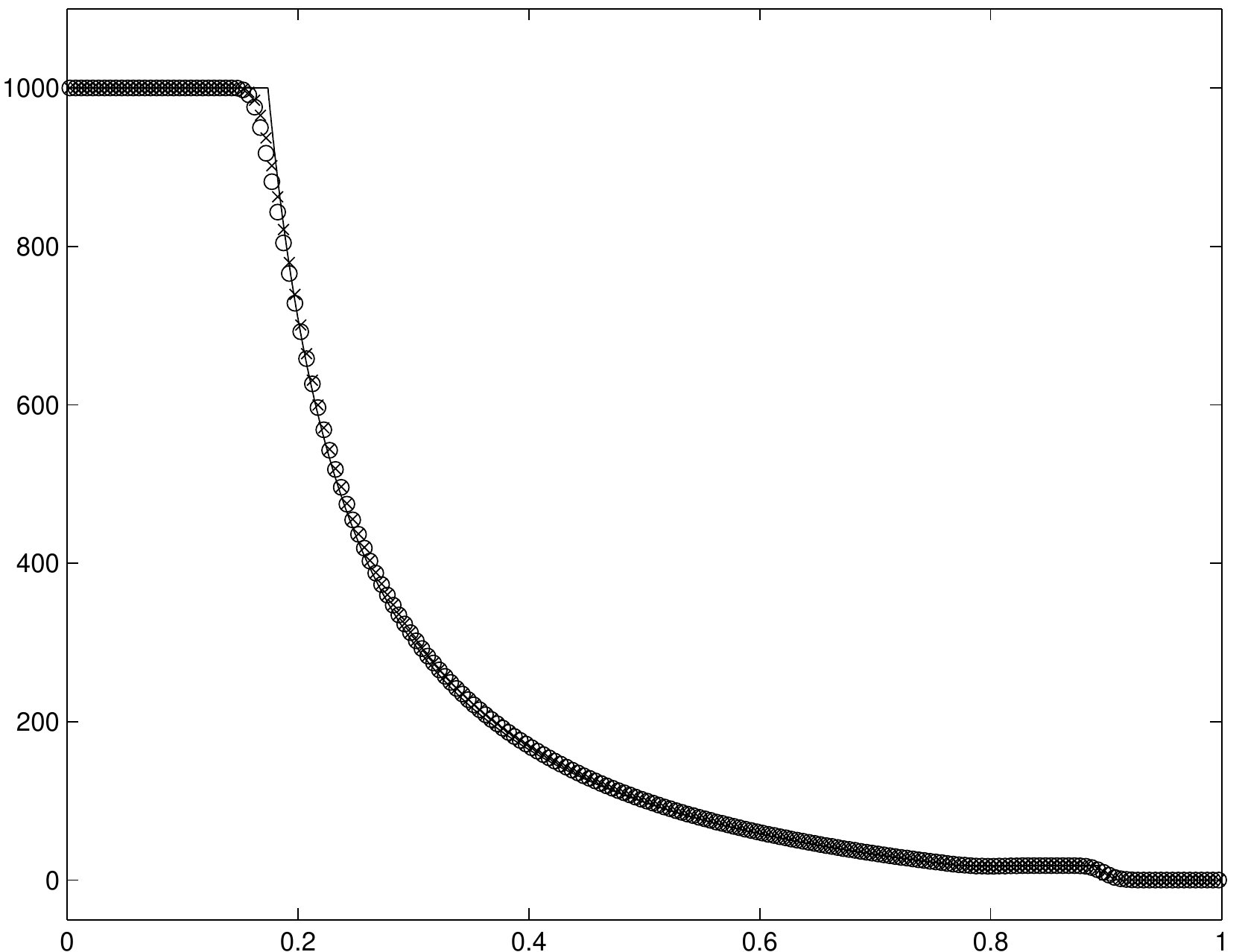}&
		\includegraphics[width=0.35\textwidth]{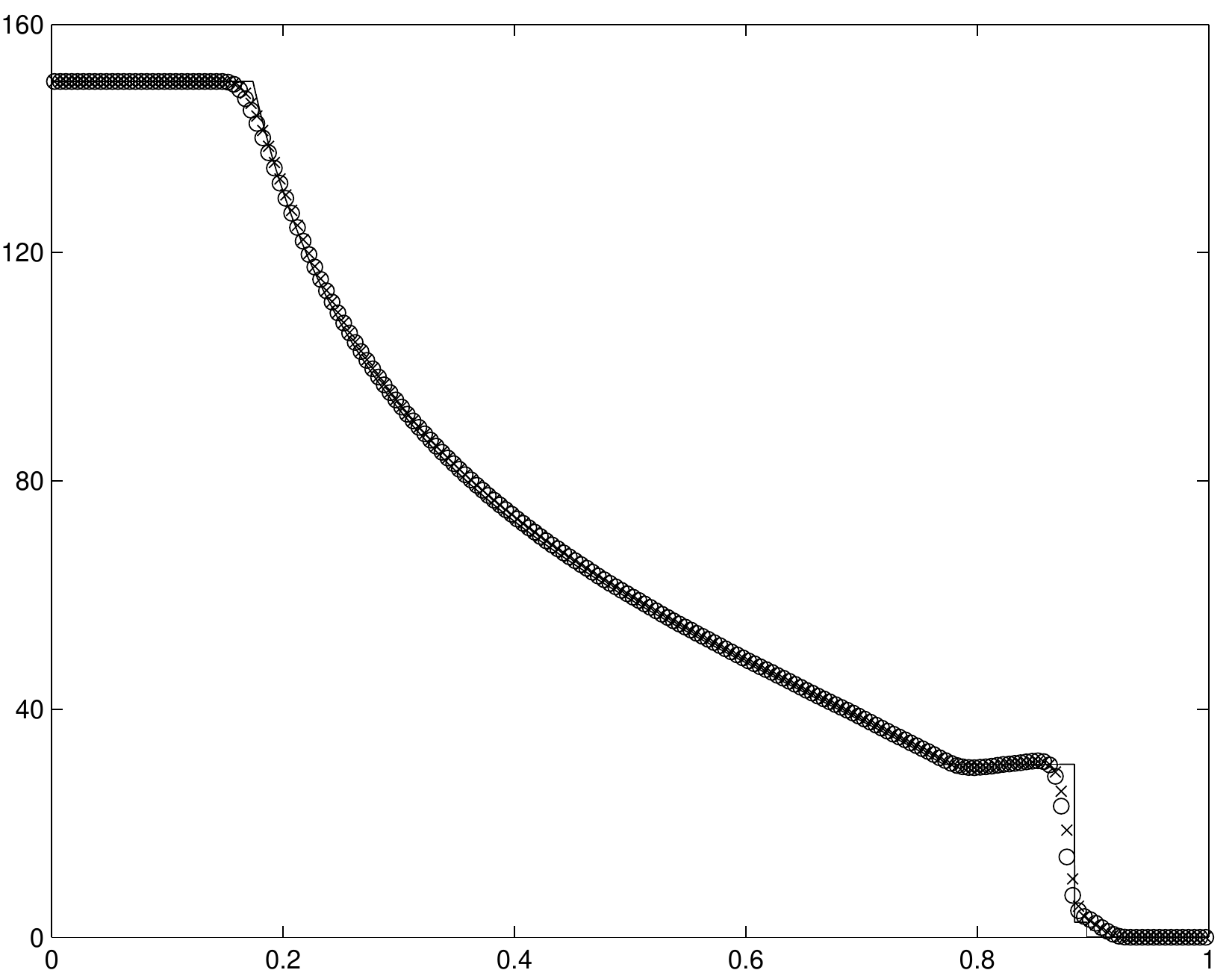}\\
		\includegraphics[width=0.35\textwidth]{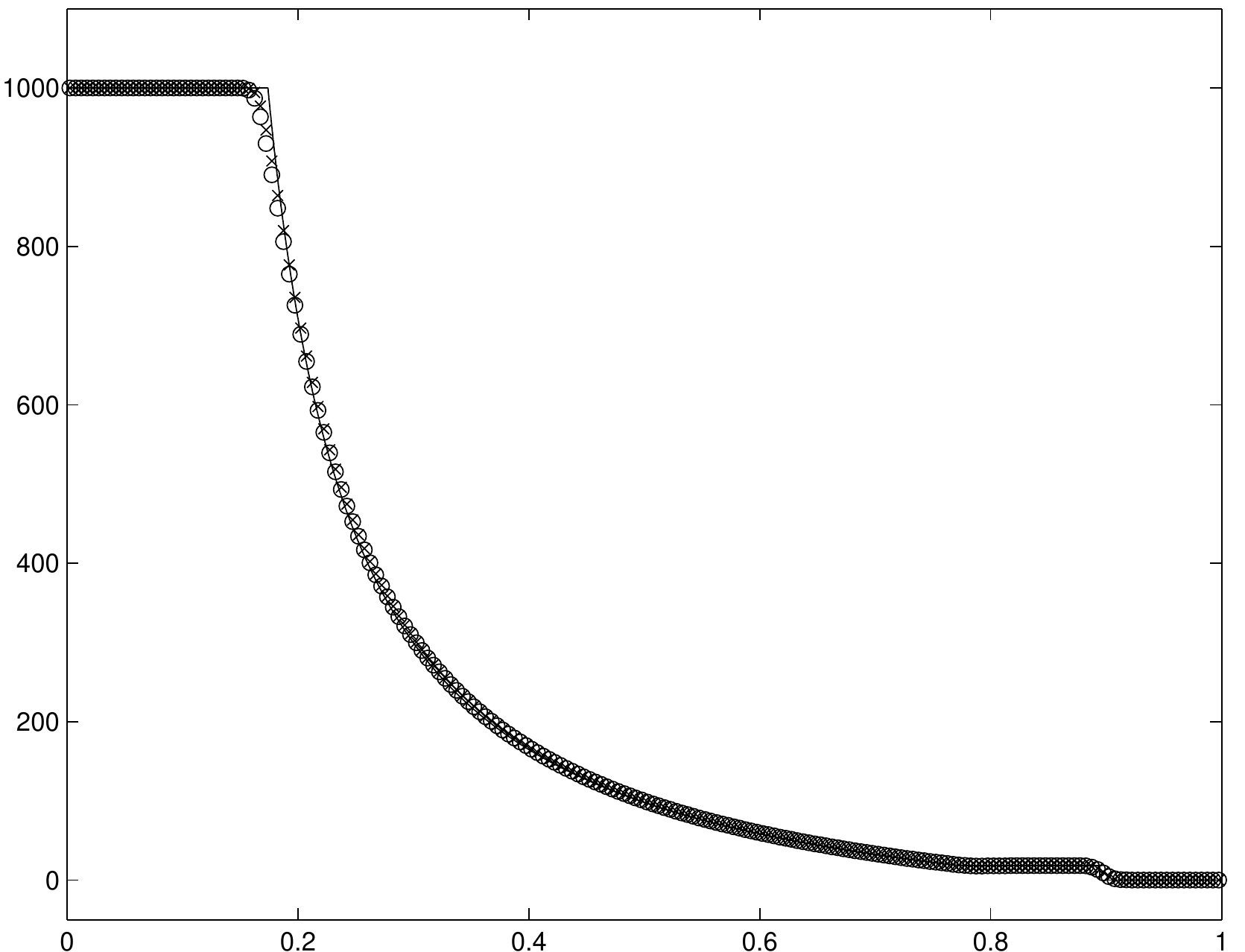}&
		\includegraphics[width=0.35\textwidth]{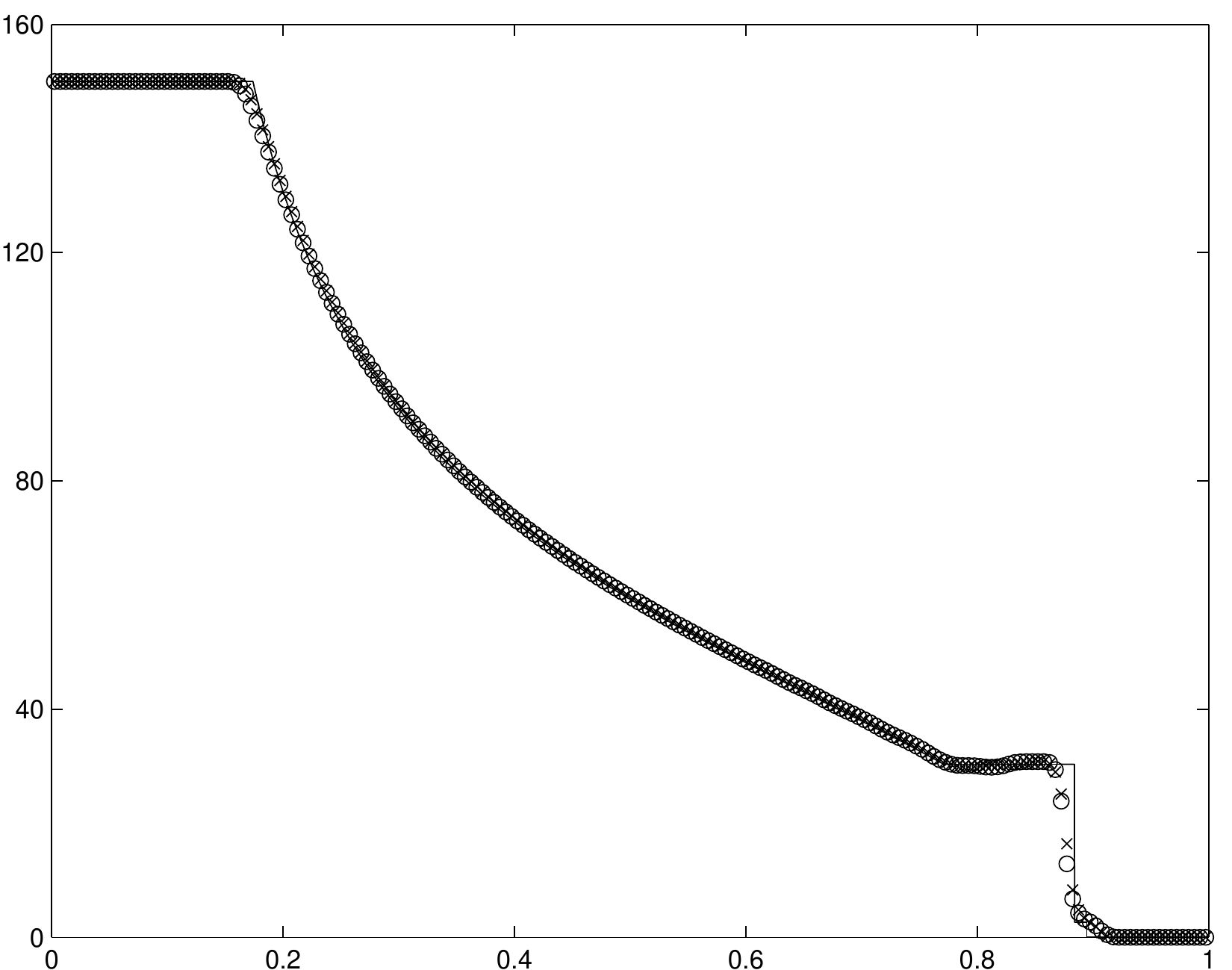}\\
	\end{tabular}
	\caption{Same as Fig.~\ref{fig:RHDRM1DT4rhovel} except for the pressure $p$ (left) and
   specific internal energy $e$ (right). }
	\label{fig:RHDRM1DT4preene}
\end{figure}

{}

\fi

\begin{table}[!htbp]
	\centering
	\caption{Example \ref{exRHDRM1DT4}: CPU times (second) of new and old \CDG{} with 800 cells. }
	\begin{tabular}{|c|c|c|}
		\hline
		& new  method   & old method \\
		\hline
		$P^1$& 15.7 & 15.3\\
		\hline
		$P^2$& 55.3 & 38.4\\
		\hline
		$P^3$& 74.2 & 81.7\\
		\hline
	\end{tabular}
	\label{tab:RHDRM1DT4cmpRC}
\end{table}

\subsection{2D case}

This section solves some 2D RHD problems by using \CDG{} on the uniform rectangular  meshes.
Those problems are the two-dimensional smooth problem, Riemann problems, and
 shock-bubble interaction problem.
The spatial stepsizes in the $x$ and $y$ directions are denoted
by $h^x$ and $h^y$ respectively.
Unless otherwise stated,
only the third-order accurate TVD Runge-Kutta time discretization \eqref{eq:RK3} is employed and
and the time step size is taken as
\begin{equation}
\label{eq:CDGCFL2D}
\Delta
t_n=\theta\tau_{n}=\frac{\theta\mu}{\max\limits_{j,k}\Big(\frac{\max_{i}|\lambda_1^{(i)}(\vec
		U^{(0)})|}{h^x}+\frac{\max_{i}|\lambda_2^{(i)}(\vec
		U^{(0)})|}{h^y}\Big) } ,
\end{equation}
where $\mu$ denotes the CFL number, and ``$\max\limits_{j,k}$'' denotes
 the maximum value over the cells $C_{j,k}$
and $D_{j+\frac{1}{2},k+\frac{1}{2}}$, while
the values of $\mu$ and $\theta$ will be given in
the coming examples.


\begin{Example}[Smooth problem]\label{exRHDSmooth2DCdg}\rm
This smooth problem has been used in \cite{ZhaoTang2013} to
test the accuracy of numerical methods.
The initial data for the primitive variables $\vec V$
is set as
$$
\vec V(x,y,0)=\big(1+0.2\sin(2\pi(x\cos\alpha +y \sin\alpha)), 0.2,0,1\big)^T,
$$
where $\alpha=30^0$ denotes the angle of the sine wave propagation
direction relative to the $x$-axis.
The computational domain $\Omega=[0,2/\sqrt{3}]\times [0,2]$  is specified with the
periodic boundary conditions, and divided into $N\times 2N$
uniform cells.

Table~\ref{tab:smoothRHDCdg2D}~ lists the $l^1$ errors of density and orders at
$t=1$ obtained by using the \CDG{},
where the fourth-order accurate Runge-Kutta time discretization is employed, $\theta=1$,
and the CFL number $\mu$ is taken as $0.3,~0.25$, and $0.2$ for
$P^1$-, $P^2$-, and $P^3$-based methods, respectively.
Those results show that the theoretical order $K+1$ of the $P^K$-based \CDG{}
may be achieved.
\end{Example}

\begin{table}[!htbp]
	\centering
	\caption{Example~\ref{exRHDSmooth2DCdg}:
$l^1$ errors of the density and orders $t = 1$ of the  $P^K$-based Runge-Kutta CDG methods with
$N\times 2N$ cells.   }
	\begin{tabular}{|c|c|c|c|c|c|}
		\hline
		\multirow{2}{20pt}{}
		&\multirow{2}{2pt}{$N$}
		&\multicolumn{2}{|c|}{without limiter}&\multicolumn{2}{|c|}{with limiter in global}\\
		\cline{3-6}
		&&$l^1$ error & order &$l^1$ error&  order\\
		\hline
		\multirow{6}{20pt}{$P^{1}$}
		&10& 9.09e-03& --& 1.76e-01& --\\
		\cline{2-6}
		&20&1.28e-03& 2.83&5.28e-02& 1.73\\
		\cline{2-6}
		&40&3.02e-04& 2.08&2.40e-02& 1.14\\
		\cline{2-6}
		&80&7.56e-05& 2.00&6.00e-03& 2.00\\
		\cline{2-6}
		&160&1.89e-05& 2.00&1.46e-03& 2.04\\
		\cline{2-6}
		&320&4.72e-06& 2.00&3.43e-04& 2.09\\
		\hline
		\multirow{6}{20pt}{$P^{2}$}
		&10& 3.43e-04& --& 2.40e-02& --\\
		\cline{2-6}
		&20&4.24e-05& 3.02&1.33e-03& 4.17\\
		\cline{2-6}
		&40&5.28e-06& 3.01&5.98e-05& 4.48\\
		\cline{2-6}
		&80&6.59e-07& 3.00&4.17e-06& 3.84\\
		\cline{2-6}
		&160&8.23e-08& 3.00&4.20e-07& 3.31\\
		\cline{2-6}
		&320&1.03e-08& 3.00&4.94e-08& 3.09\\
		\hline
		\multirow{6}{20pt}{$P^{3}$}
		&10& 2.53e-05& --& 2.78e-03& --\\
		\cline{2-6}
		&20&1.55e-06& 4.03&7.26e-05& 5.26\\
		\cline{2-6}
		&40&9.61e-08& 4.01&9.60e-07& 6.24\\
		\cline{2-6}
		&80&5.99e-09& 4.00&1.66e-08& 5.85\\
		\cline{2-6}
		&160&3.75e-10& 4.00&5.07e-10& 5.04\\
		\cline{2-6}
		&320&2.34e-11& 4.00&3.31e-11& 3.94\\
		\hline
	\end{tabular}\label{tab:smoothRHDCdg2D}
\end{table}

\begin{Example}[Riemann problem 1]\label{exRHDRM2DT1}\rm{}
The initial data of the first 2D Riemann problem are
	$$(\rho,v_1,v_2,p)(x,y,0)=\begin{cases}(0.035145216124503,0,0,0.162931056509027),& x>0,y>0,\\
	(0.1,0.7,0,1),&    x<0,y>0,\\
	(0.5,0,0,1),&      x<0,y<0,\\
	(0.1,0,0.7,1),&    x>0,y<0,\end{cases}$$
where  the left and bottom discontinuities are two contact discontinuities
and the top and right are two shock waves with the
speed of $0.934563275373844$.

In our computations, $\theta$ is taken as 1 or 0.5,
and 
the value of $\mu\theta$ is fixed as $0.3,~0.25$, and $0.2$
for the $P^1$-, $P^2$-, $P^3$-based \CDG{}, respectively.
The results at $t=0.8$ obtained by the $P^K$-based \CDG{}
are presented in
Figs.~\ref{fig:RHDRMT1P1cdg}, \ref{fig:RHDRMT1P2cdg}, and
~\ref{fig:RHDRMT1P3cdg}.
Fig.~\ref{fig:cmpallRHD2DT1} gives the density at $t=0.8$ along the line
$y=x$. Table~\ref{tab:cellperRM2DT1} shows the percentage of ``troubled'' cells.
It is seen that the resolution of $P^K$-based \DG{} is better than $P^K$-based \CDG{},
when $\mu\theta$ is fixed,
and the \CDG{} with small $\theta$
 improve the resolution of the discontinuity
better than the case of big $\theta$, especially for the $P^1$-based CDG method. 

\end{Example}

\ifx\outnofig\undefine
\begin{figure}[!htbp]
	\centering{}
	\begin{tabular}{cc}
		\includegraphics[width=0.35\textwidth]{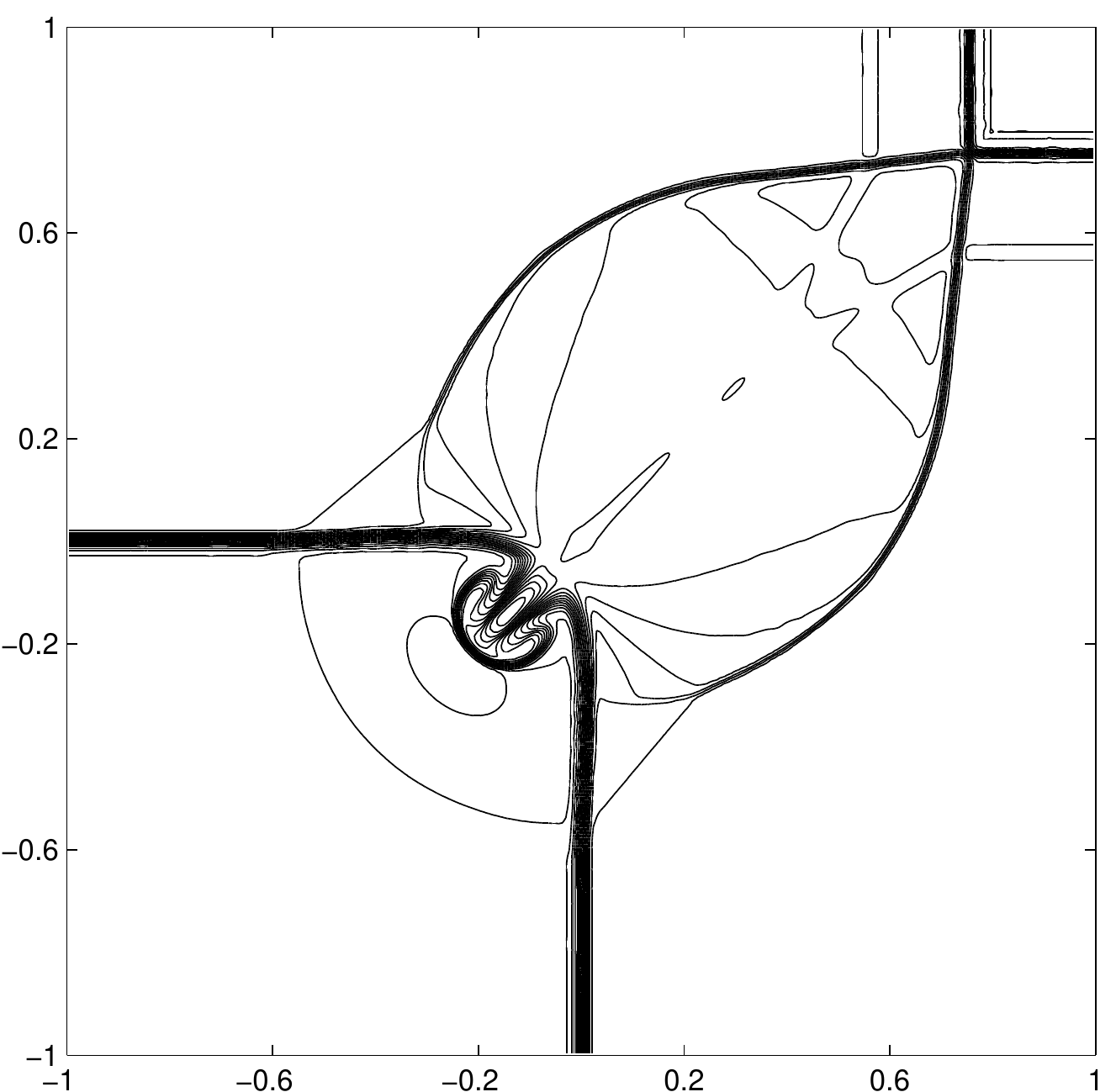}&
		\includegraphics[width=0.35\textwidth]{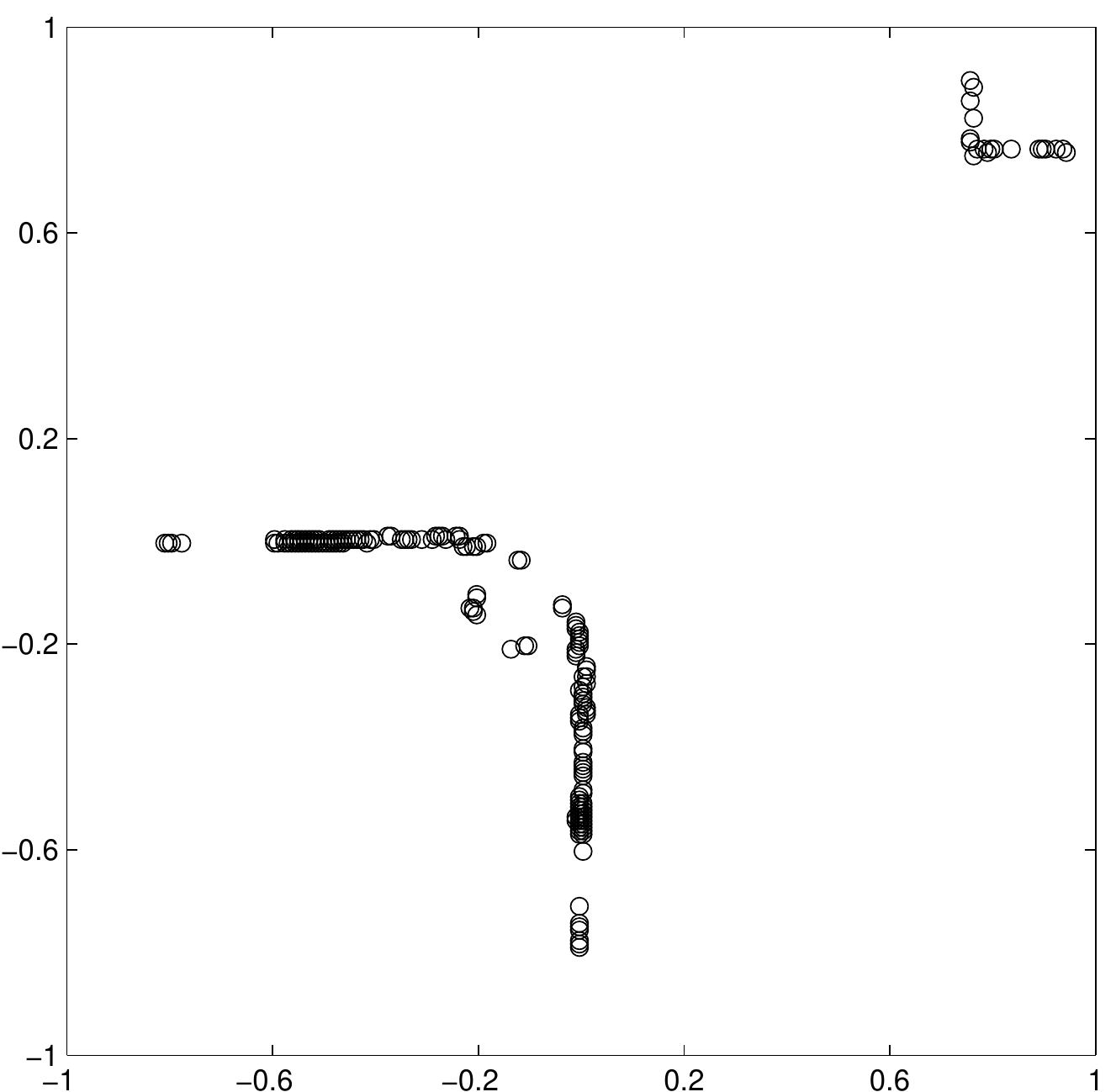}\\
		\includegraphics[width=0.35\textwidth]{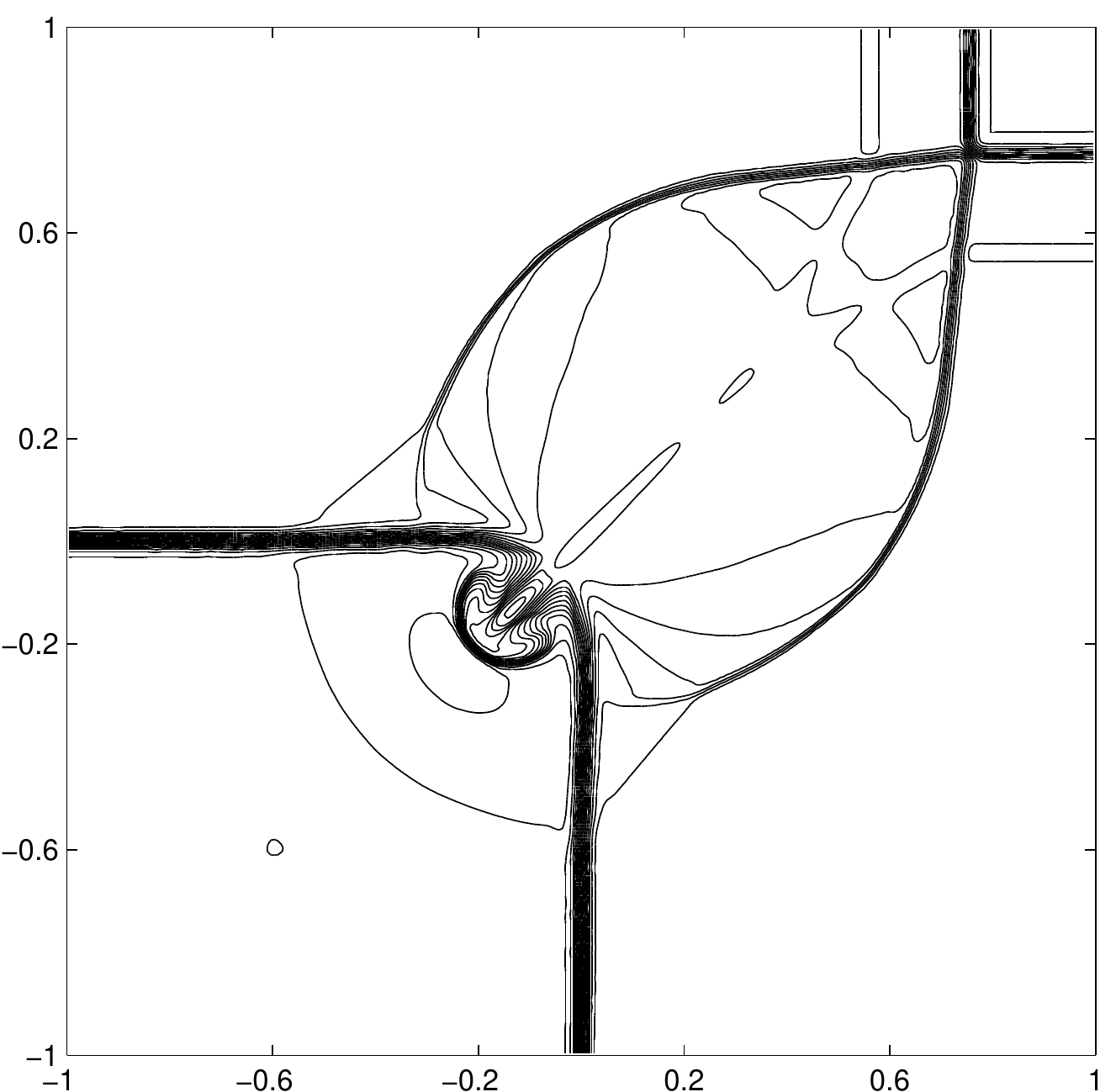}&
		\includegraphics[width=0.35\textwidth]{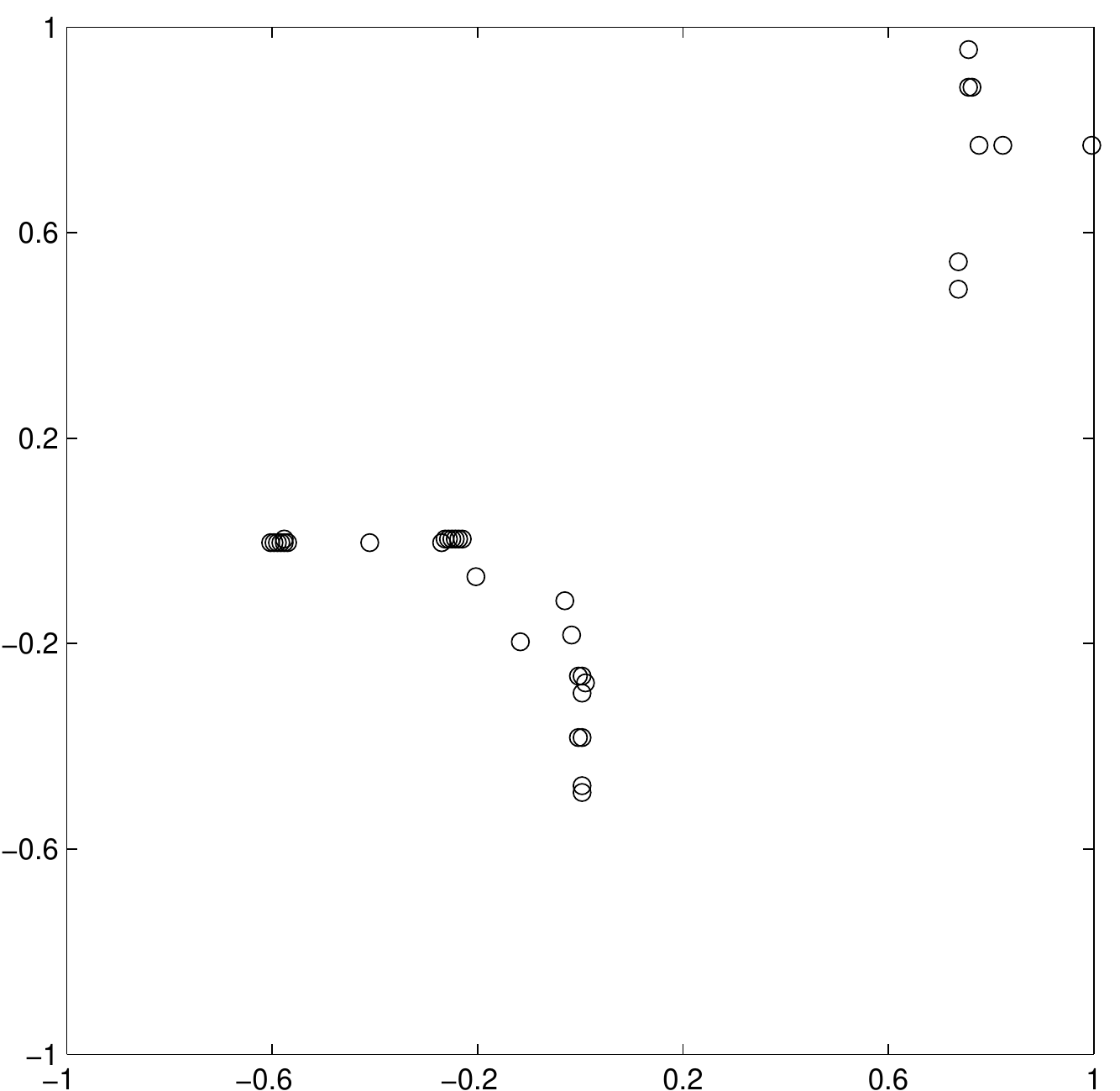}\\
		\includegraphics[width=0.35\textwidth]{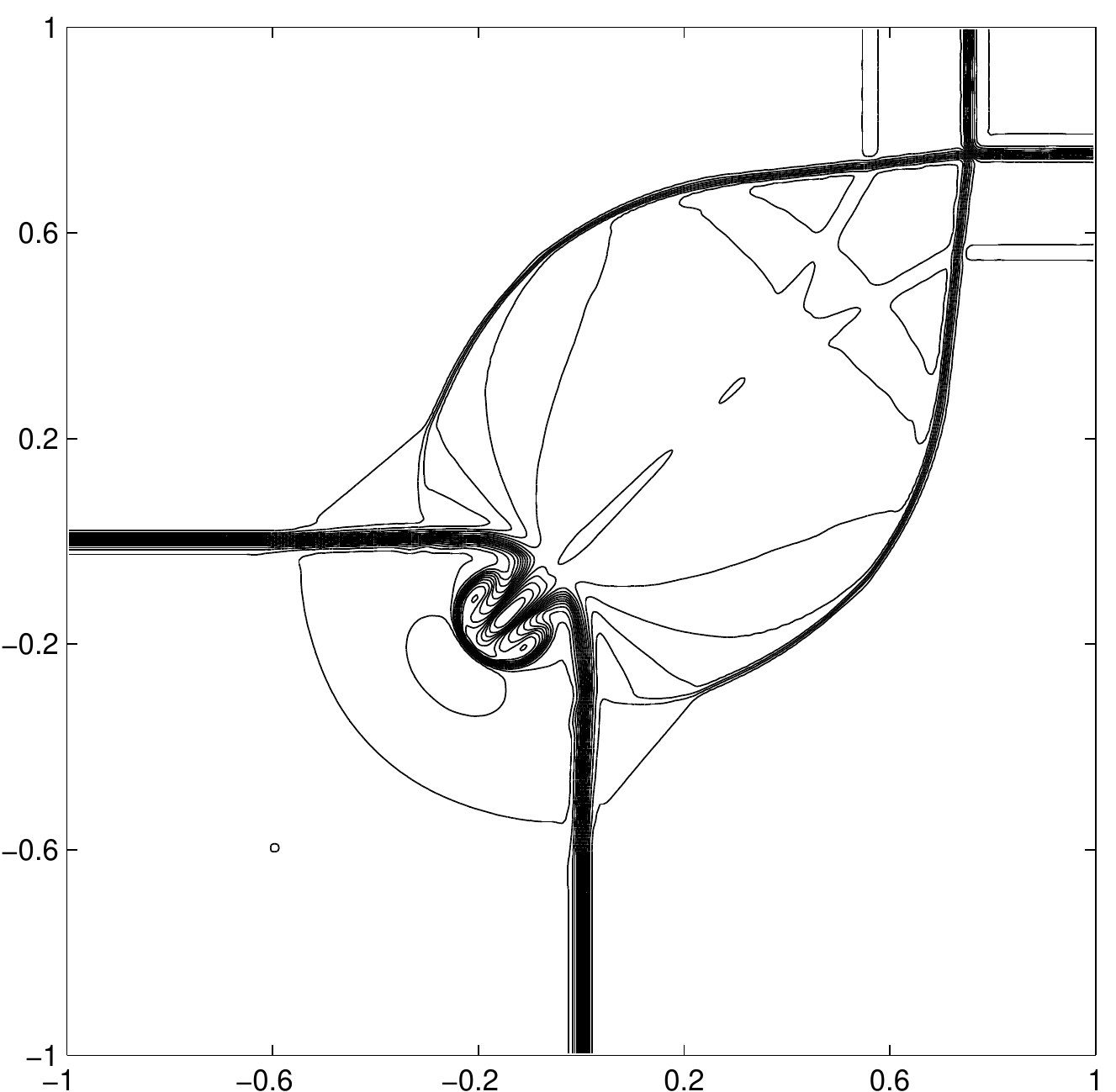}&
		\includegraphics[width=0.35\textwidth]{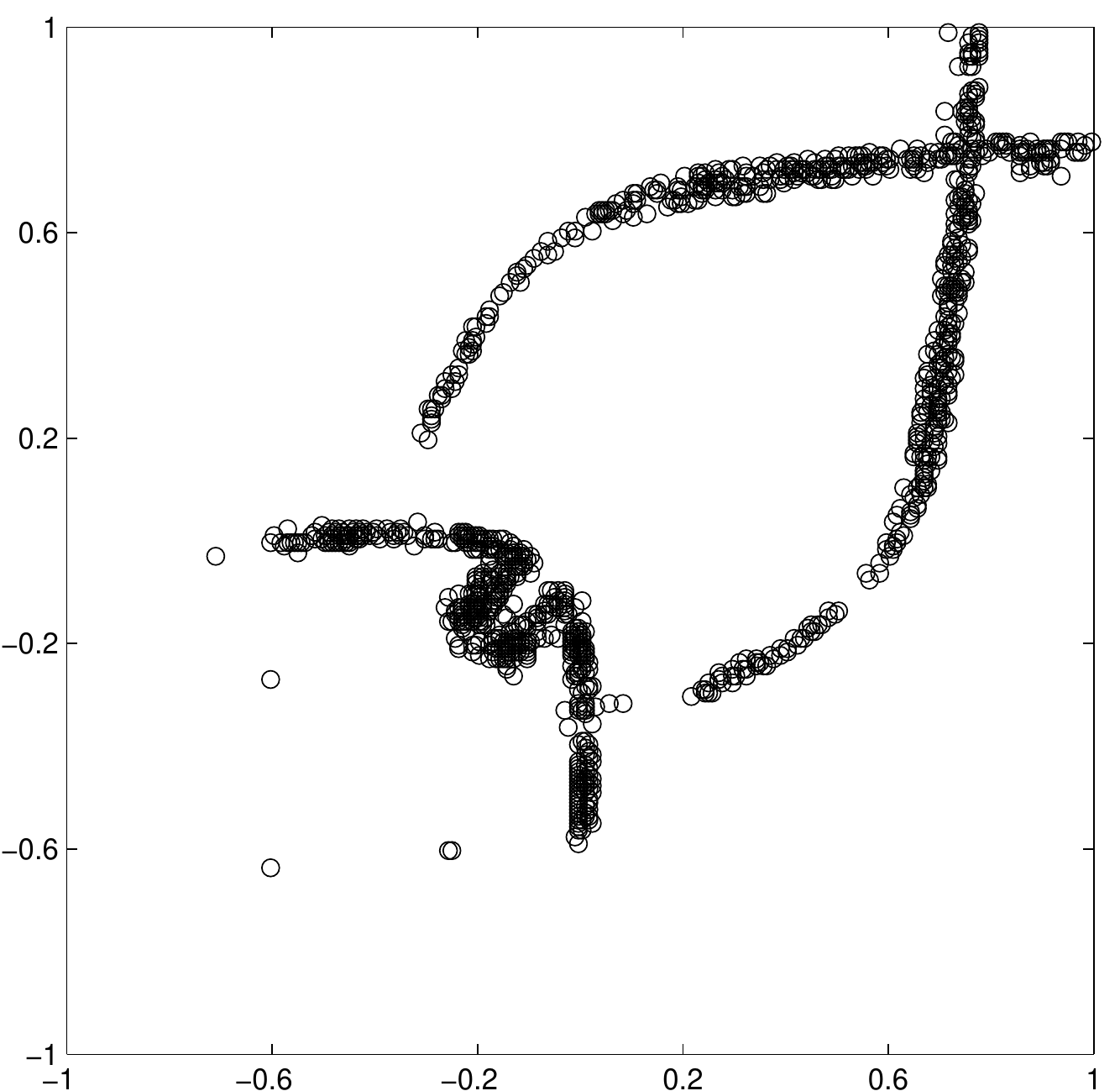}\\
	\end{tabular}
	\caption{Example \ref{exRHDRM2DT1}:
The contour plots of density logarithm  $\log\rho$ at
 (30 equally spaced contour lines from
$-1.46$ to $-0.18$) and the  ``troubled'' cells $t=0.8$ obtained with with $300\times 300$ cells.
Top:  $P^1$-based DG, middle: $P^1$-based CDG with $\theta=1$,
bottom: $P^1$-based CDG with $\theta=0.5$.
}
	\label{fig:RHDRMT1P1cdg}
\end{figure}

\begin{figure}[!htbp]
	\centering{}
	\begin{tabular}{cc}
		\includegraphics[width=0.35\textwidth]{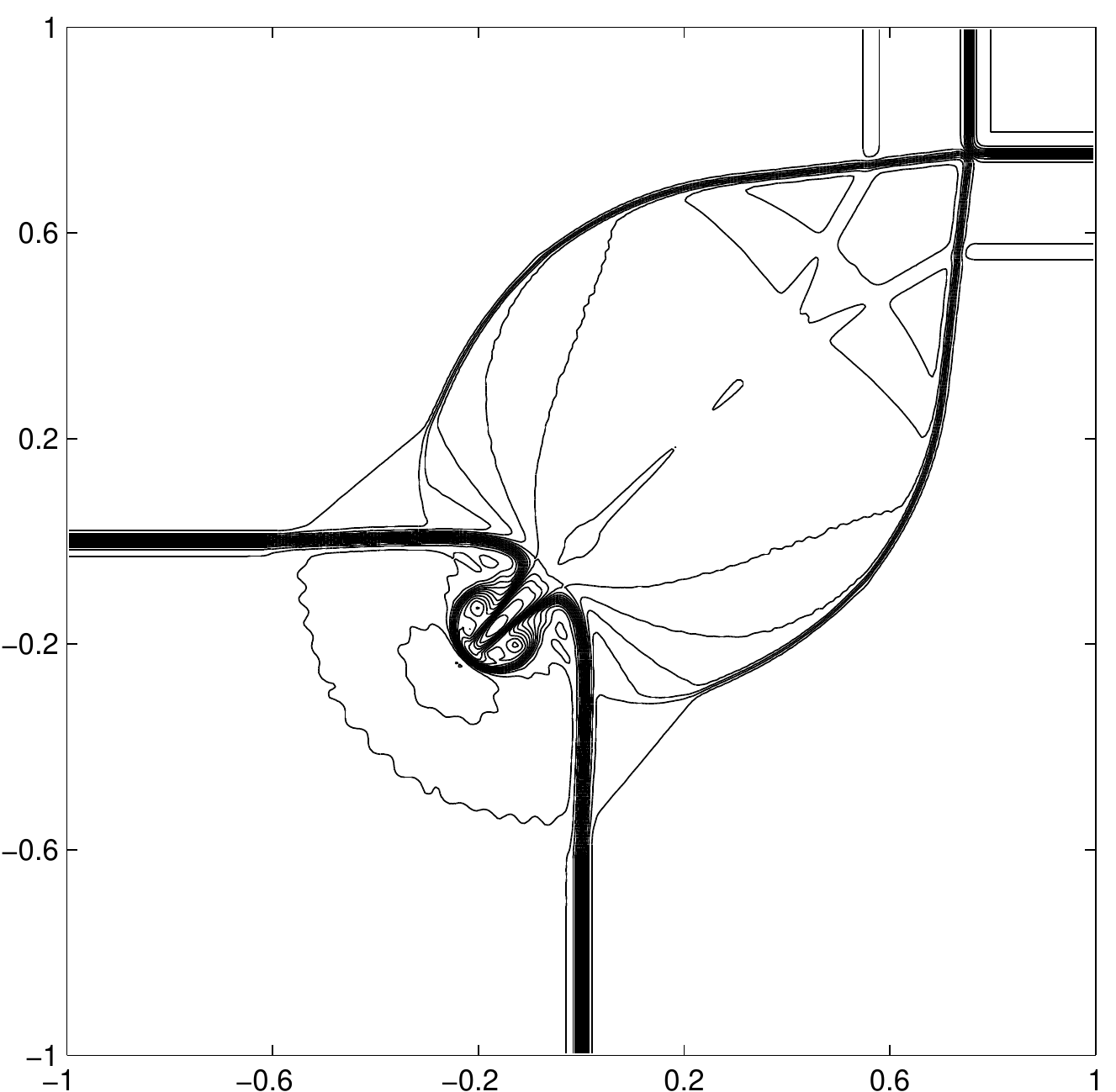}&
		\includegraphics[width=0.35\textwidth]{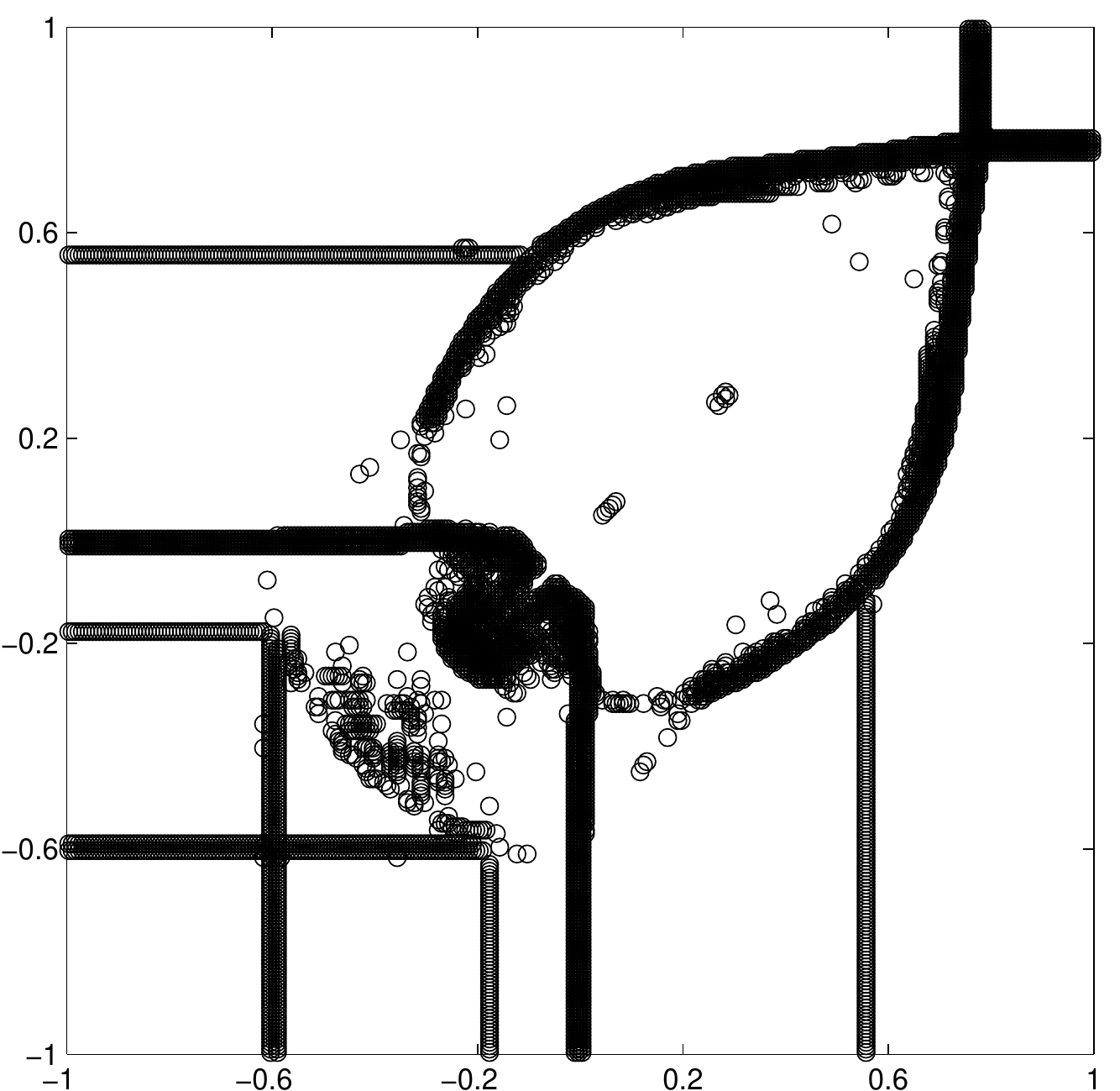}\\
		\includegraphics[width=0.35\textwidth]{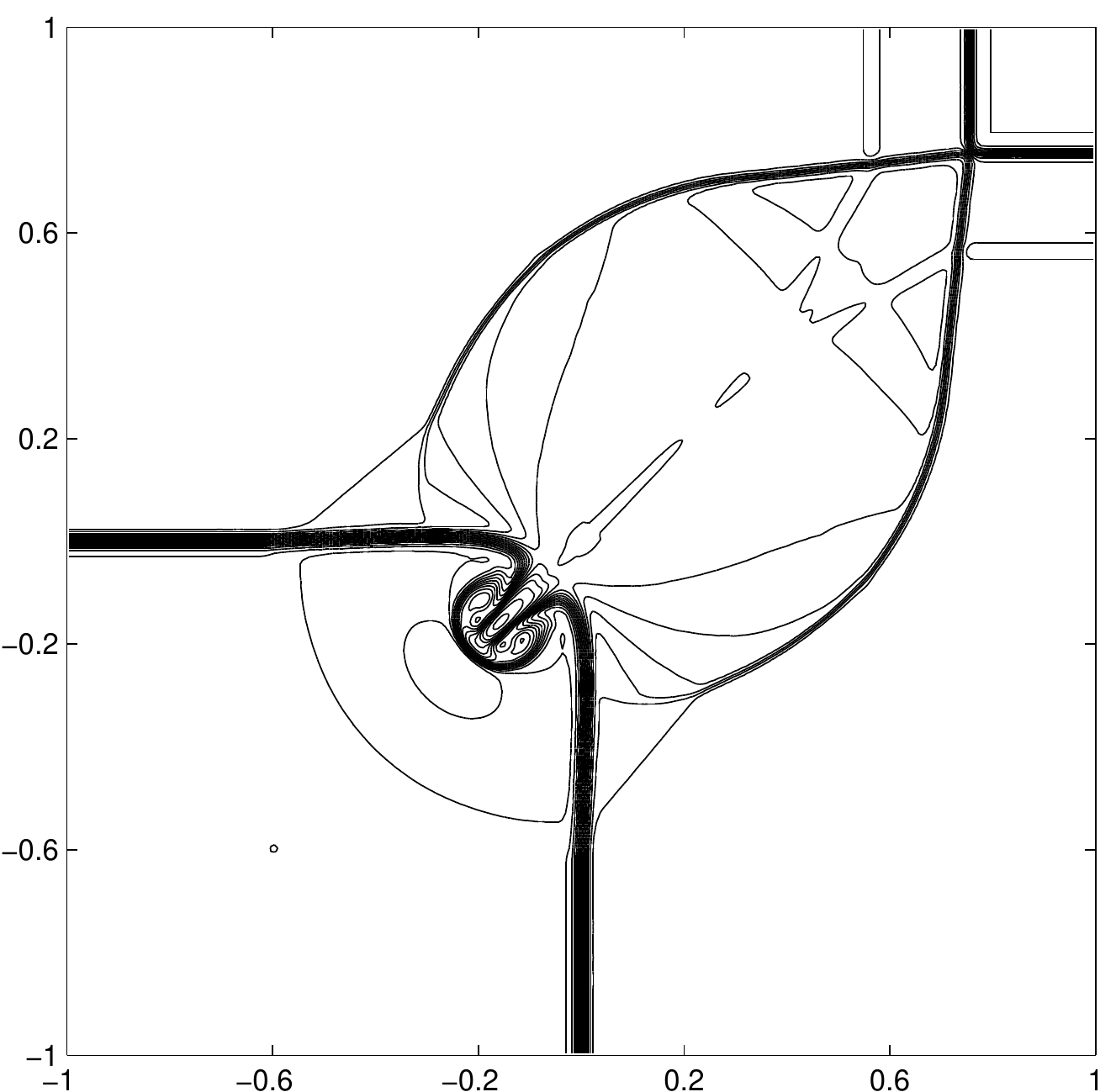}&
		\includegraphics[width=0.35\textwidth]{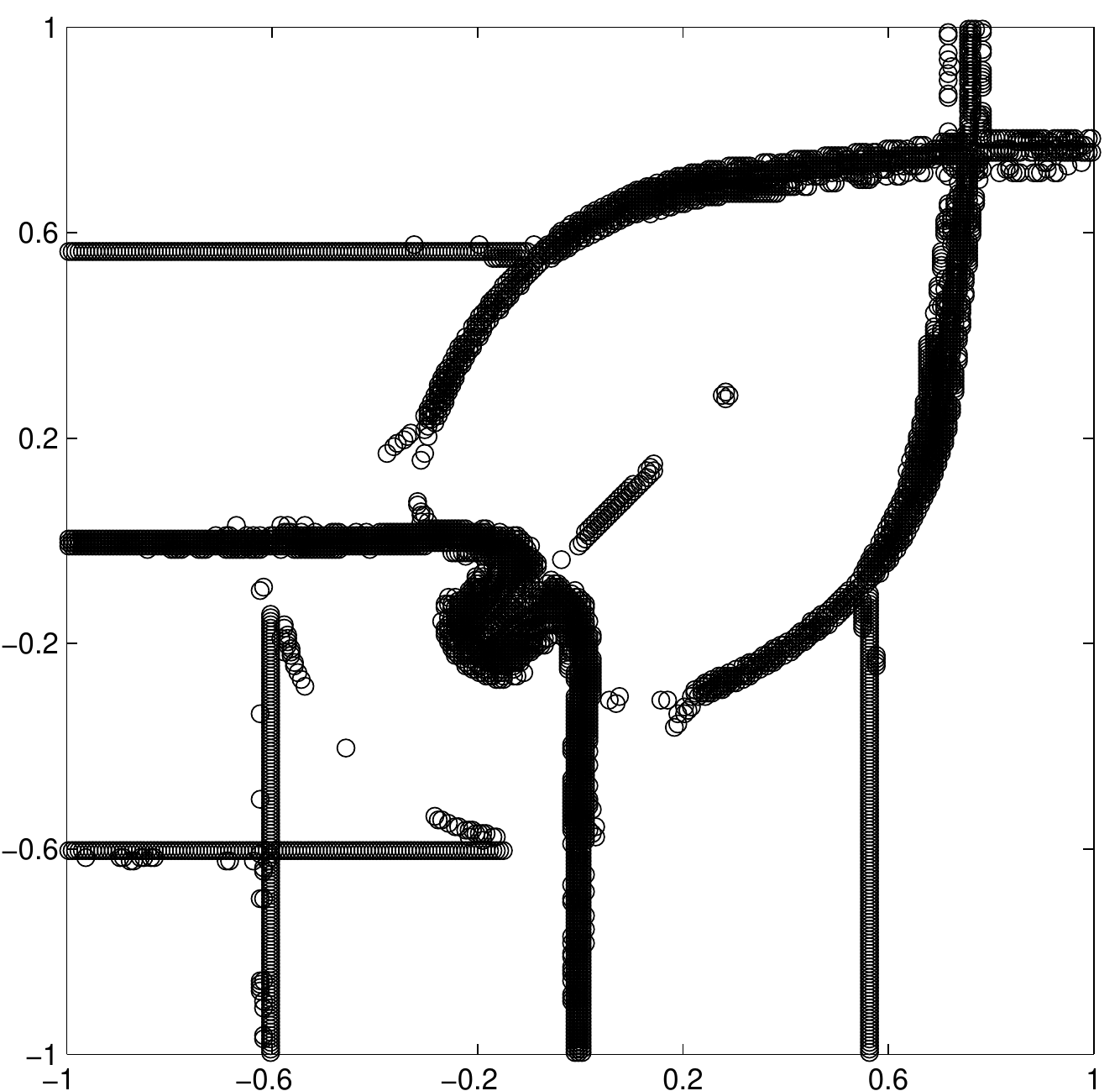}\\
		\includegraphics[width=0.35\textwidth]{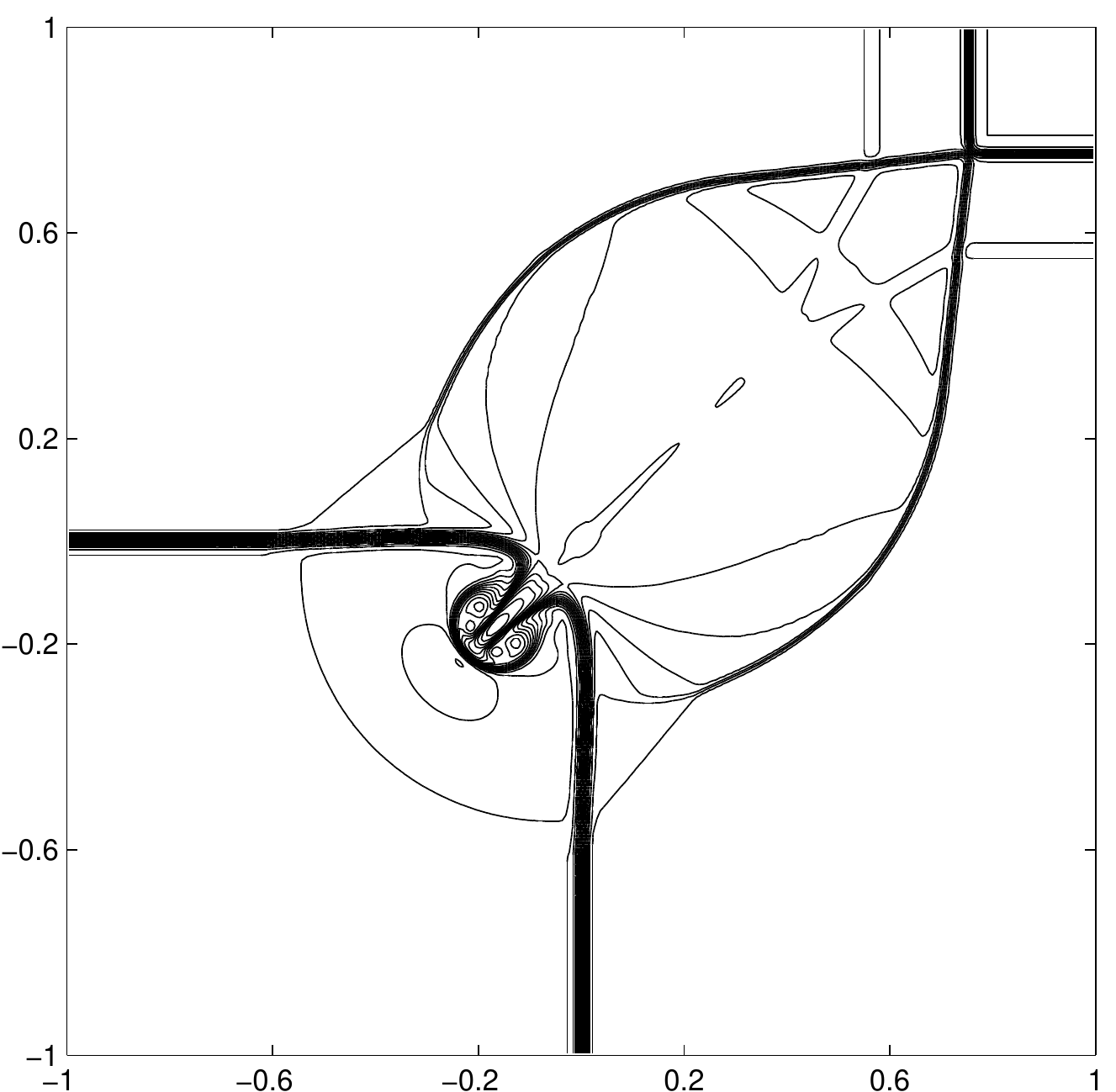}&
		\includegraphics[width=0.35\textwidth]{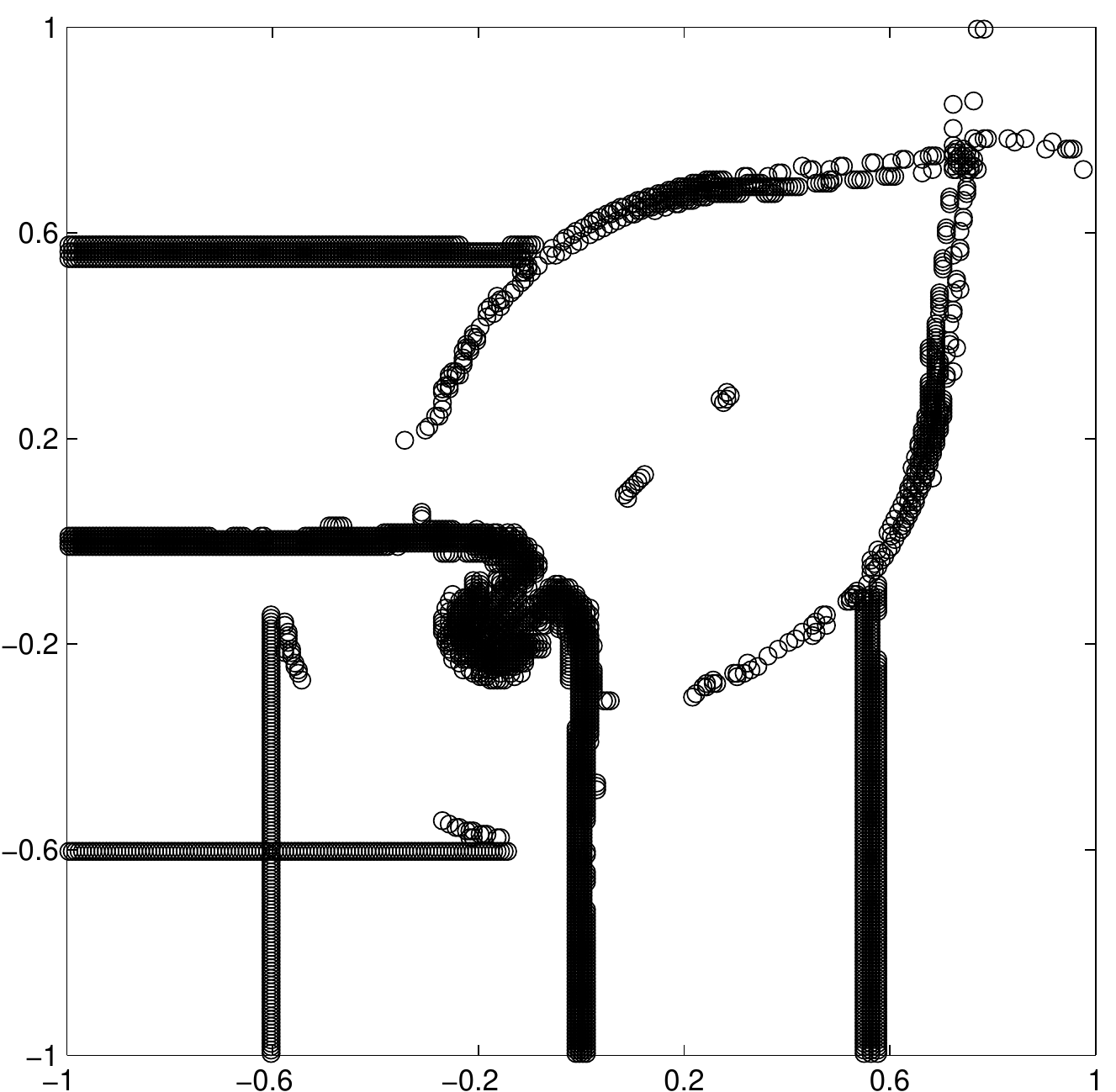}\\
	\end{tabular}
	\caption{Same as Fig.~\ref{fig:RHDRMT1P1cdg} except for the $P^2$-based \CDG{}. }
	\label{fig:RHDRMT1P2cdg}
\end{figure}

\begin{figure}[!htbp]
	\centering{}
	\begin{tabular}{cc}
		\includegraphics[width=0.3\textwidth]{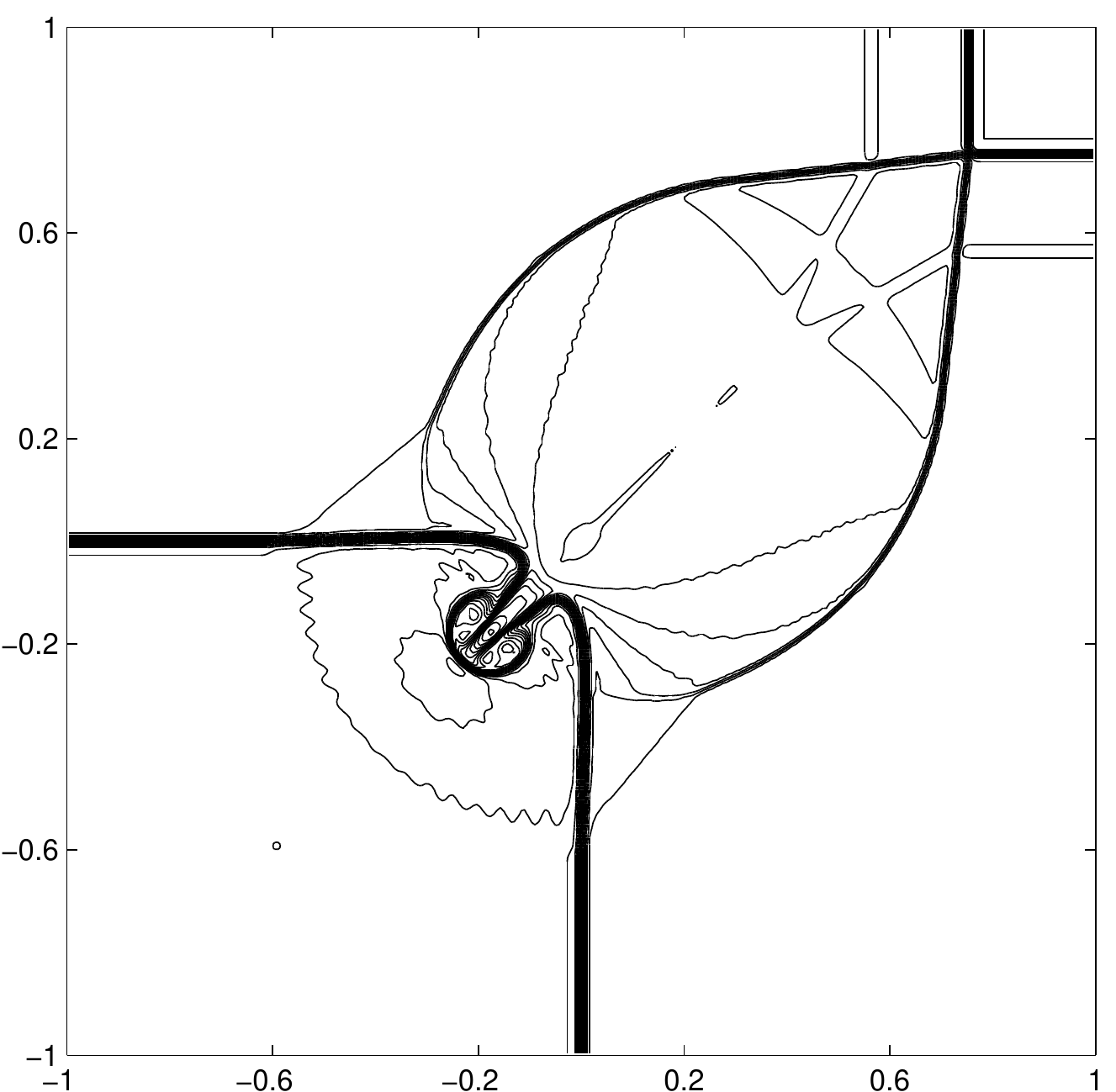}&
		\includegraphics[width=0.3\textwidth]{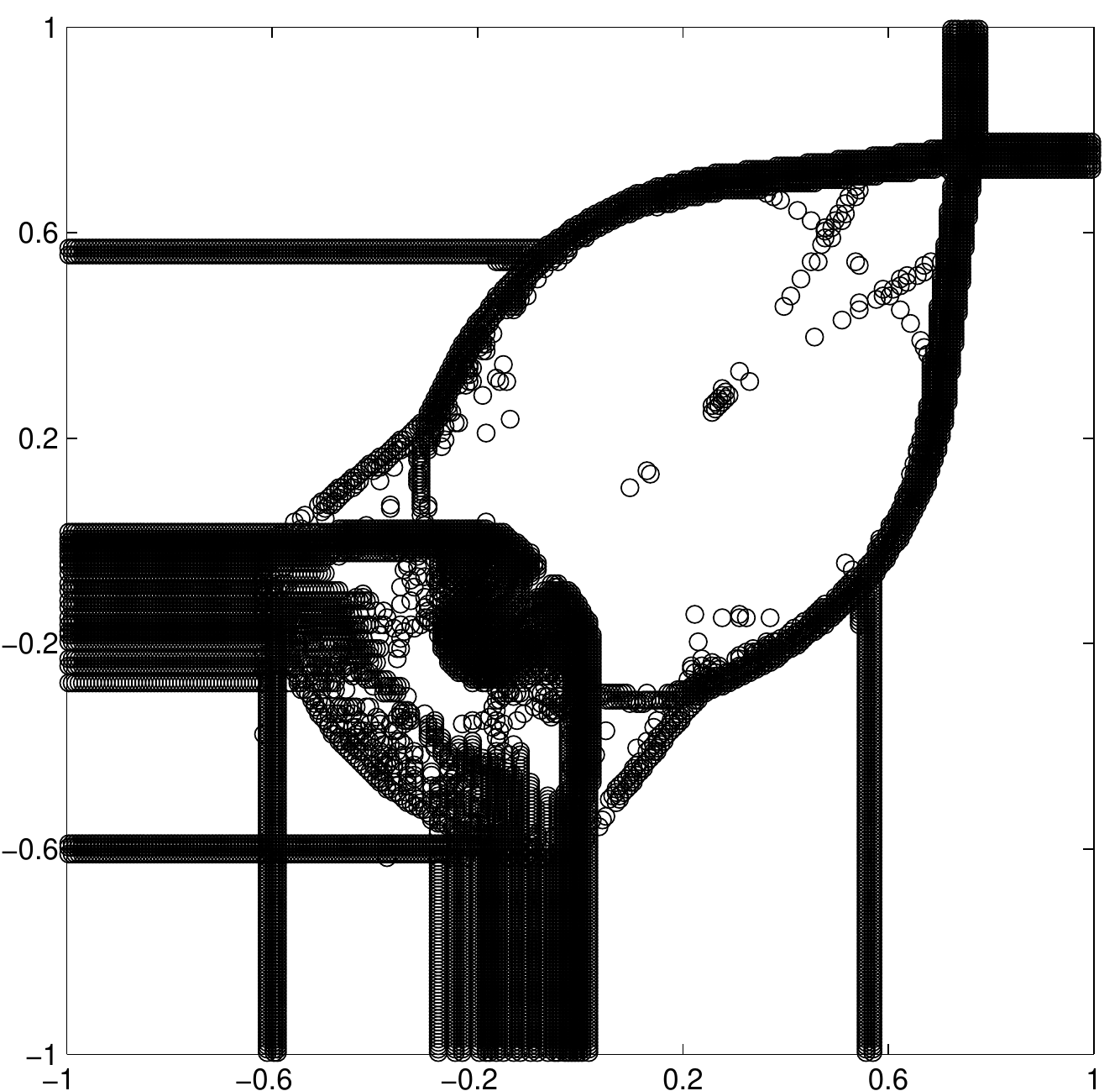}\\
		\includegraphics[width=0.3\textwidth]{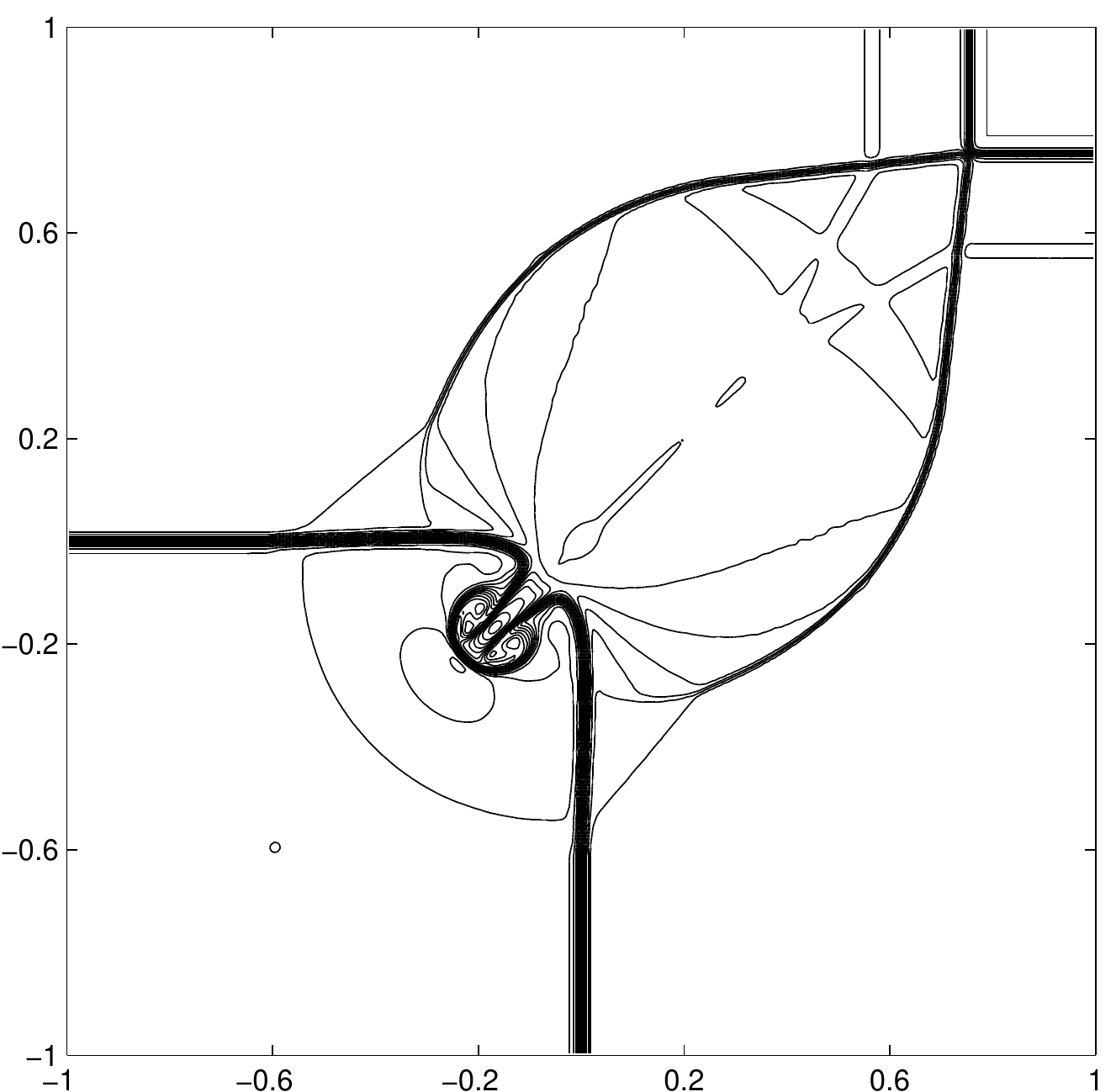}&
		\includegraphics[width=0.3\textwidth]{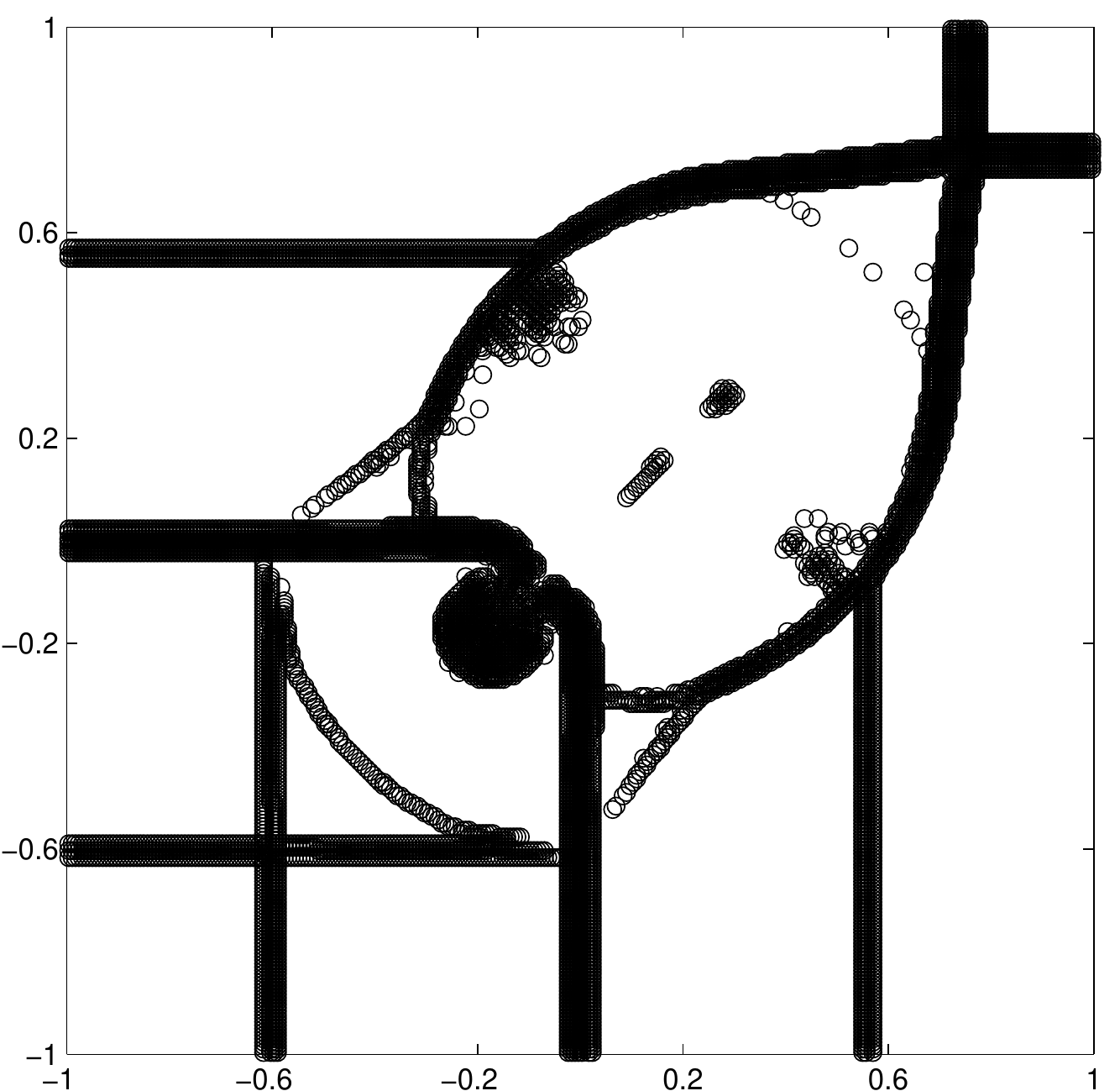}\\
		\includegraphics[width=0.3\textwidth]{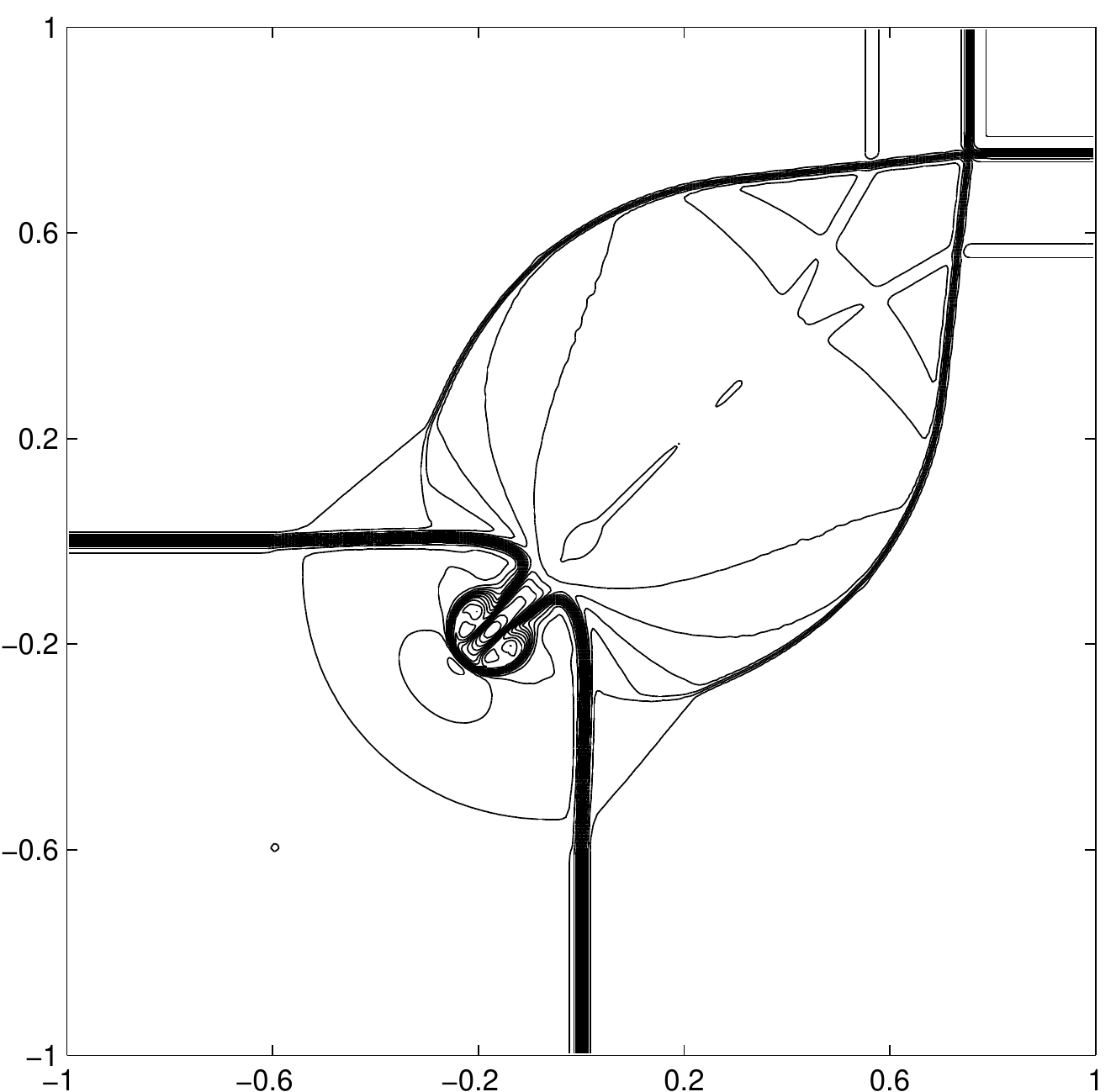}&
		\includegraphics[width=0.3\textwidth]{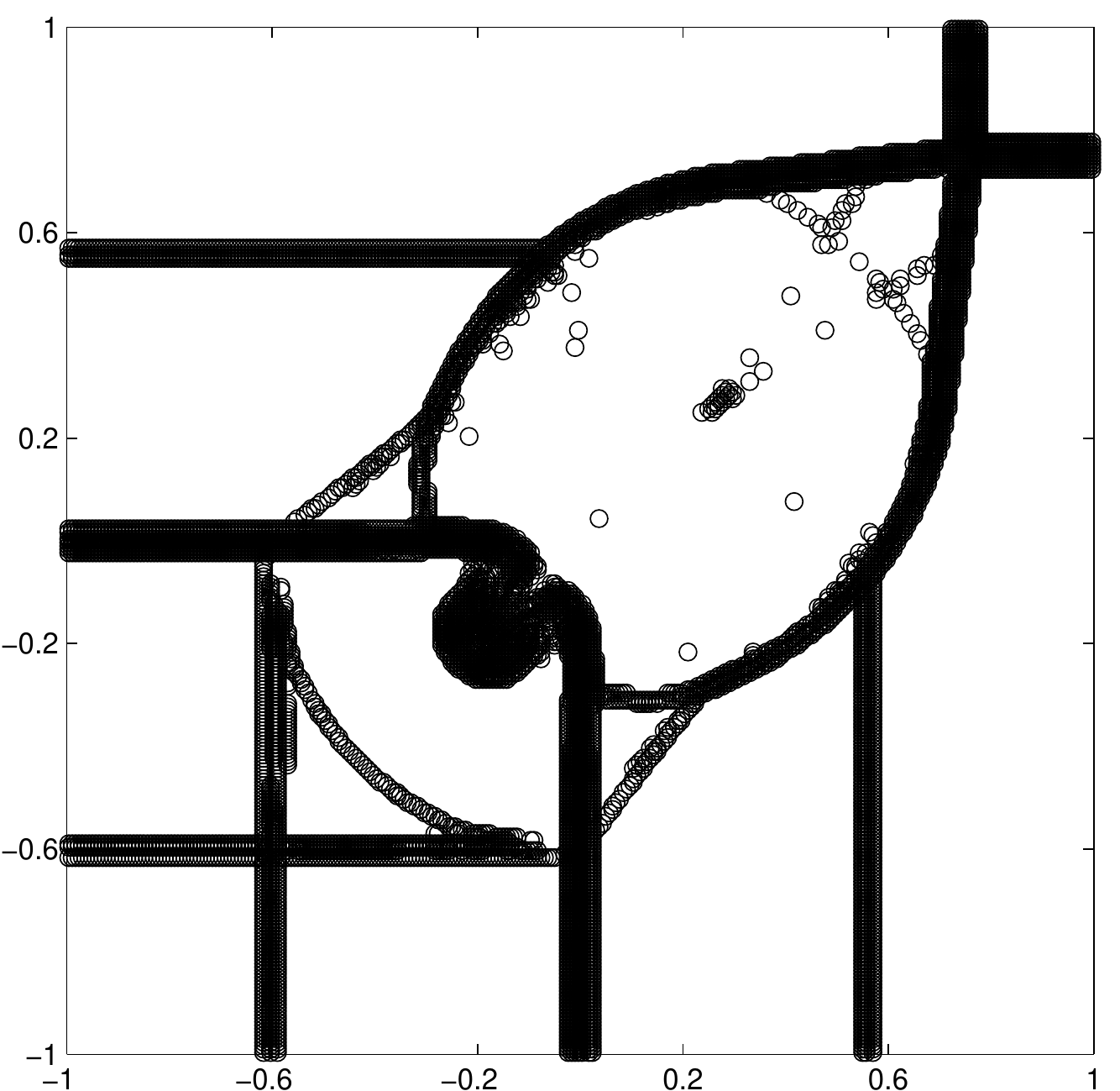}\\
	\end{tabular}
	\caption{Same as Fig.~\ref{fig:RHDRMT1P1cdg} except for the $P^3$-based \CDG{}. }
	\label{fig:RHDRMT1P3cdg}
\end{figure}

\begin{figure}[!htbp]
	\centering{}
	\begin{tabular}{ccc}
		\includegraphics[width=0.33\textwidth]{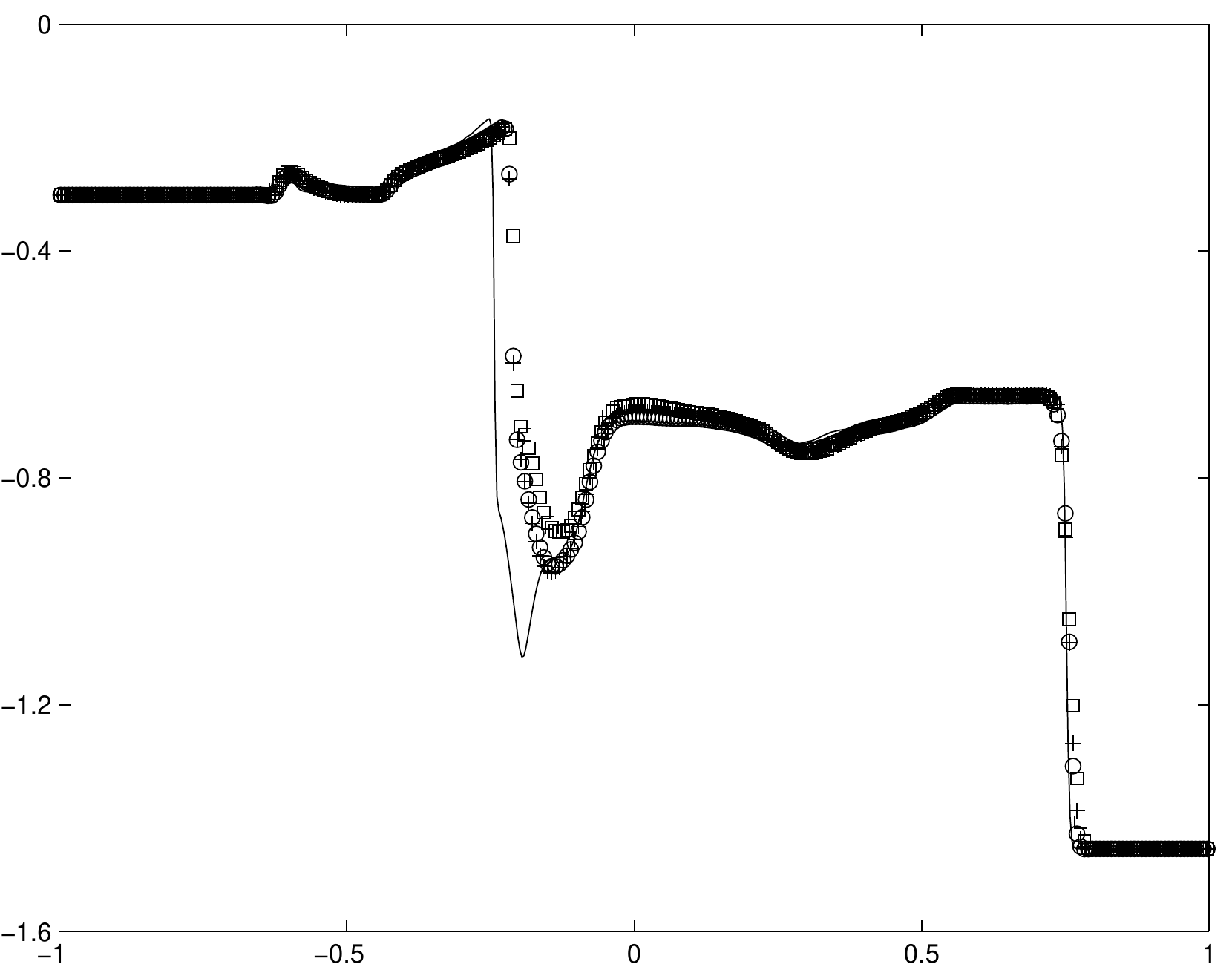}&
		\includegraphics[width=0.33\textwidth]{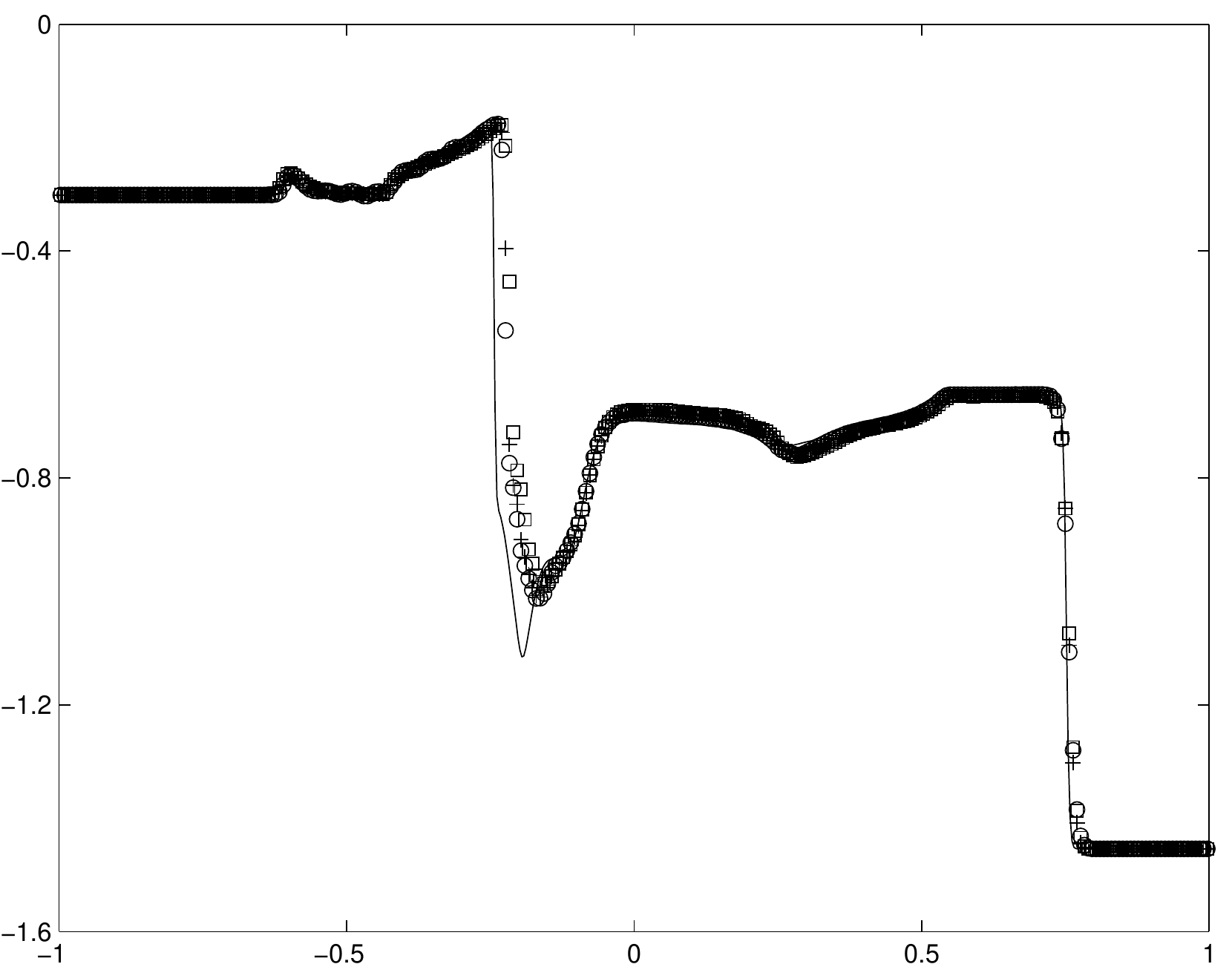}&
		\includegraphics[width=0.33\textwidth]{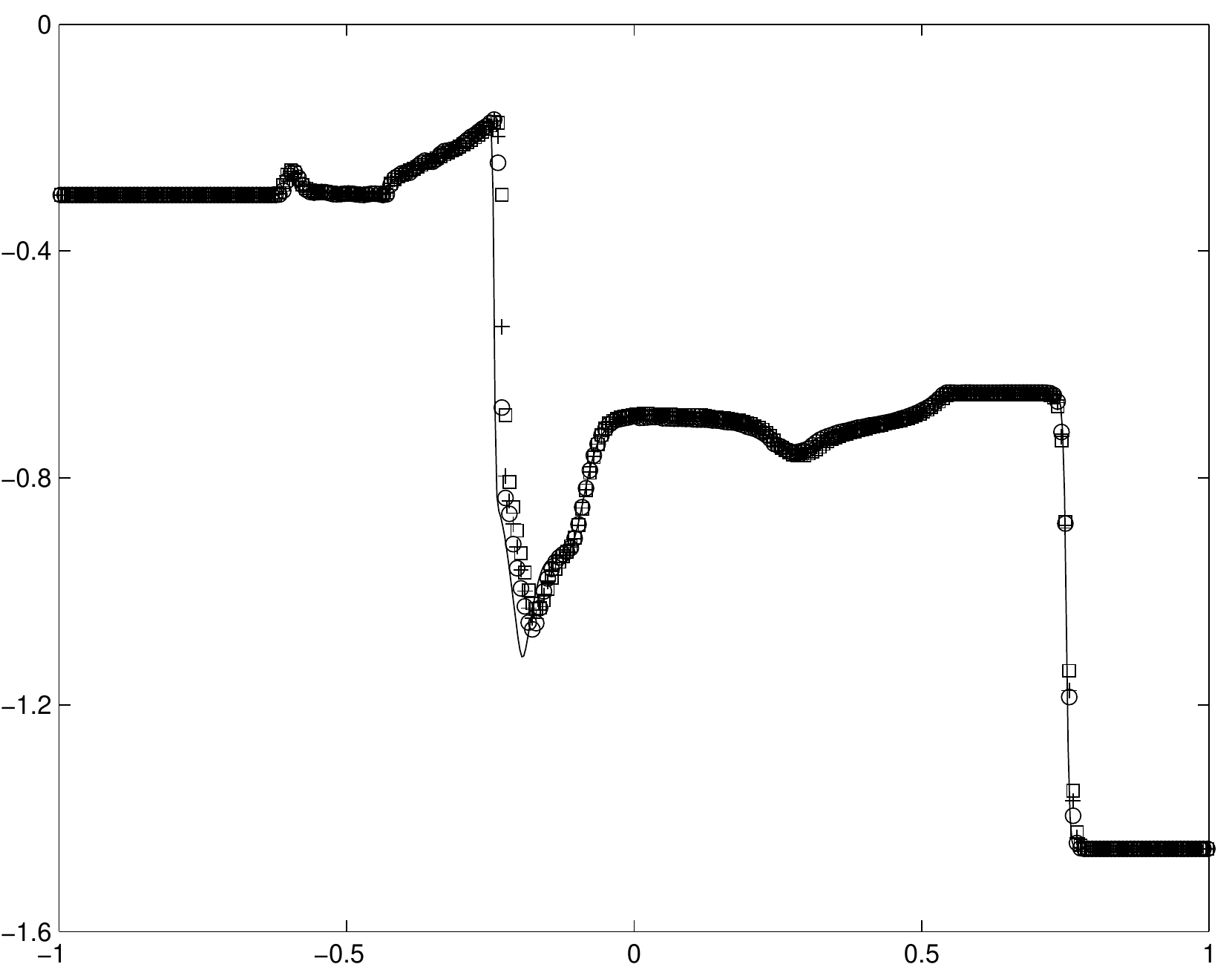}\\
	\end{tabular}
	\caption{Example \ref{exRHDRM2DT1}:
 Density logarithm  $\log\rho$ at
 $t=0.8$ along the line $y=x$.
The solid line denotes the reference solution obtained by using fifth order accurate
 WENO scheme with $600\times 600$ uniform cells,
 while the symbol ¡°$\circ$¡±, ¡°$+$¡±, and ``$\square$''
 are the solutions by the \DG{},  \CDG{} ($\theta=0.5$), \CDG{} ($\theta=1$), respectively.
 From top to bottom, $K=1,~2,~3$.
}
	\label{fig:cmpallRHD2DT1}
\end{figure}

\fi{}

\begin{table}[!htbp]
	\centering
	\caption{Example \ref{exRHDRM2DT1}: The percentage of ¡°troubled¡± cells  at $t=0.8$. }
	\begin{tabular}{|c|c|c|c|}
		\hline
		& \multirow{2}{80pt}{non-central DG} & \multicolumn{2}{|c|}{CDG}\\
		\cline{3-4}
		& & $\theta=1$ & $\theta=0.5$\\
		\hline
		$P^1$& 0.19 & 0.04  & 1.12 \\
		\hline
		$P^2$& 5.51 &  4.78 & 3.68\\
		\hline
		$P^3$& 10.68 & 8.14 & 7.95\\
		\hline
	\end{tabular}
	\label{tab:cellperRM2DT1}
\end{table}

\begin{table}[!htbp]
	\centering
	\caption{Example \ref{exRHDRM2DT1}: CPU times (second) of the CDG and non-central DG methods. }
	\begin{tabular}{|c|c|c|}
		\hline
		& CDG & non-central DG\\
		\hline
		$P^1$& 5728.7 & 2491.5\\
		\hline
		$P^2$& 16722.4 & 6119.8\\
		\hline
		$P^3$& 54490.8 & 19287.8\\
		\hline
	\end{tabular}
	\label{tab:RHDRM2DT1cmpRC}
\end{table}

\begin{figure}[!htbp]
	\centering{}
	\begin{tabular}{cc}
		\includegraphics[width=0.35\textwidth]{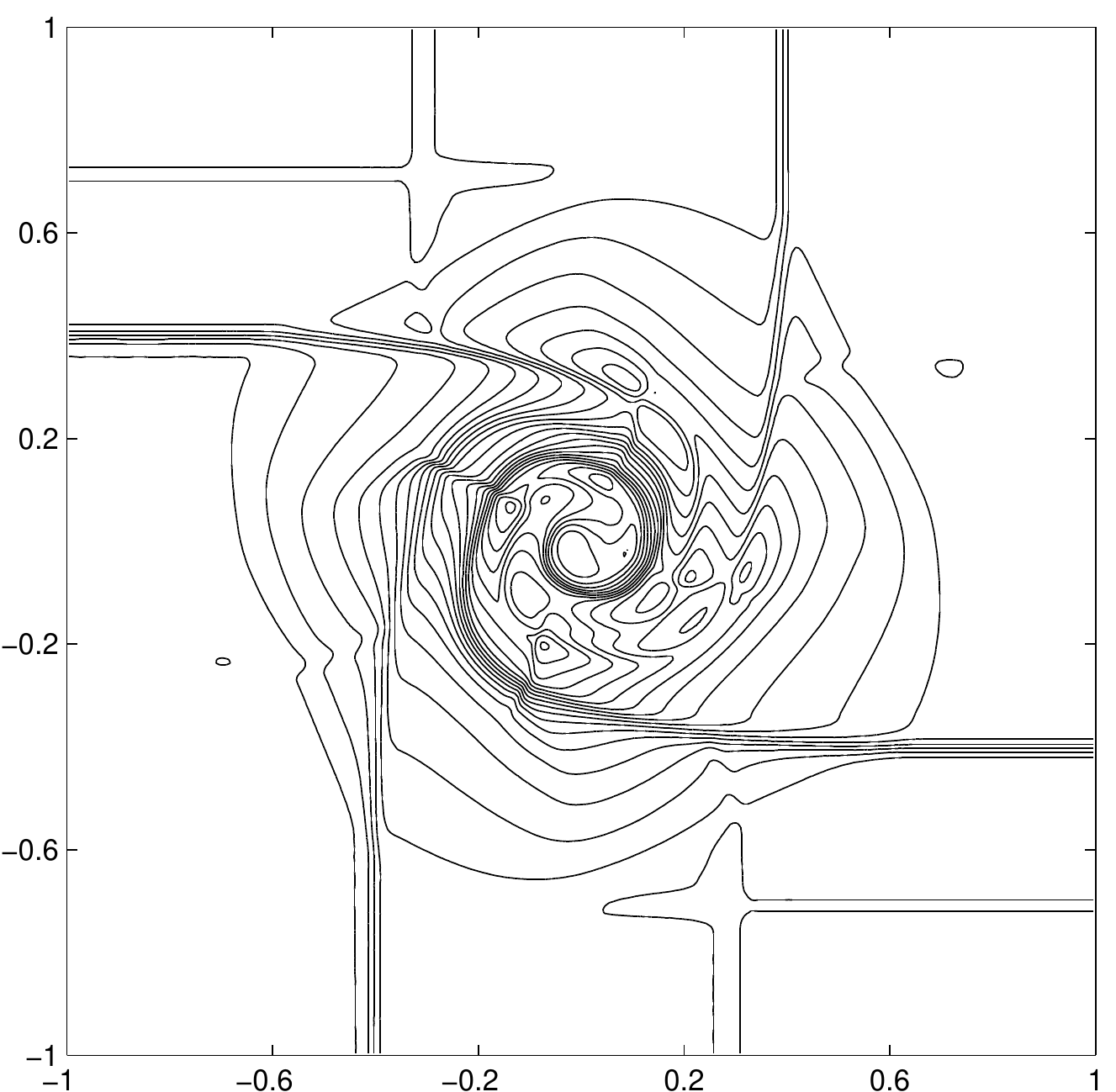}&
		\includegraphics[width=0.35\textwidth]{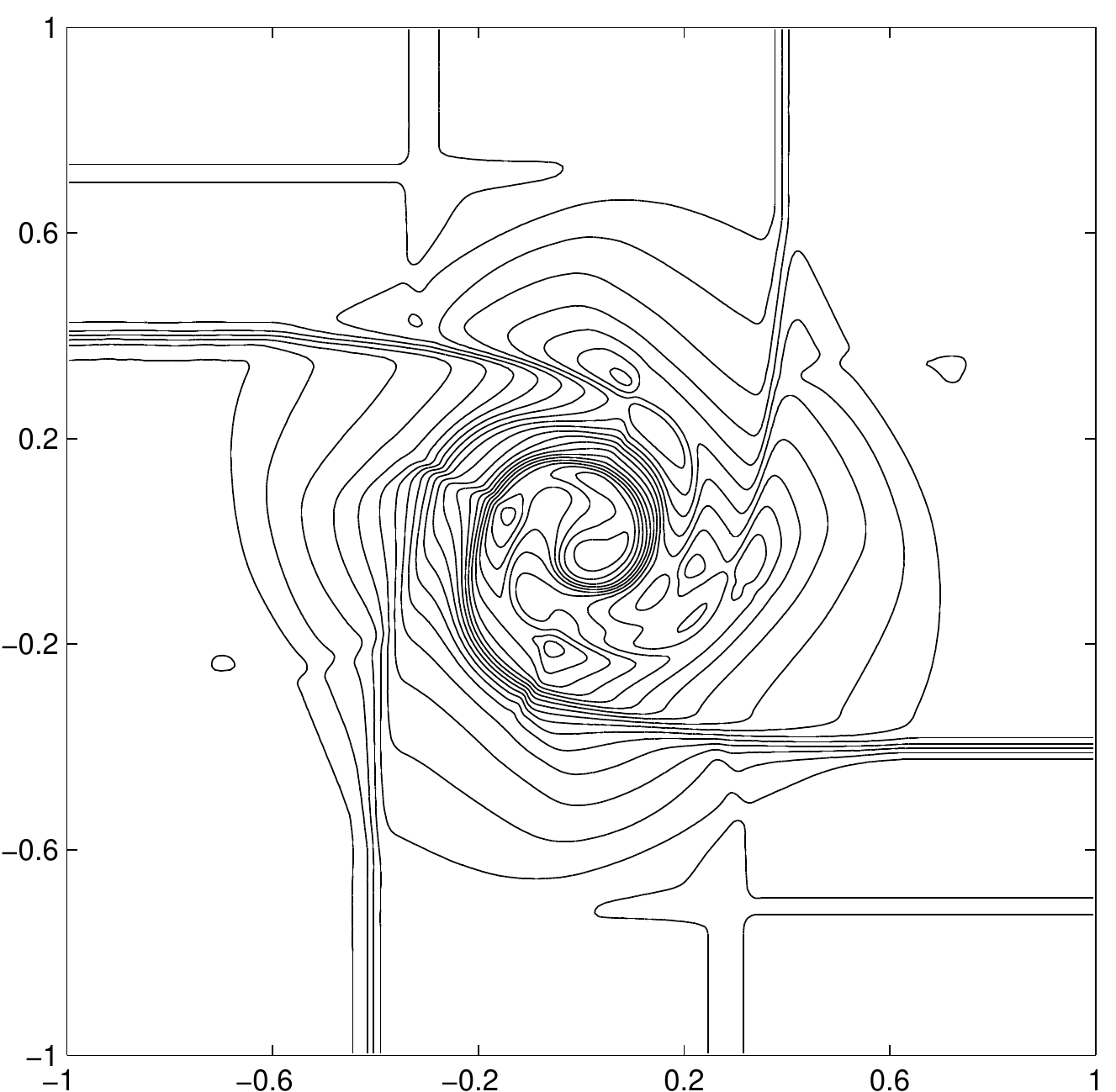}\\
		\includegraphics[width=0.35\textwidth]{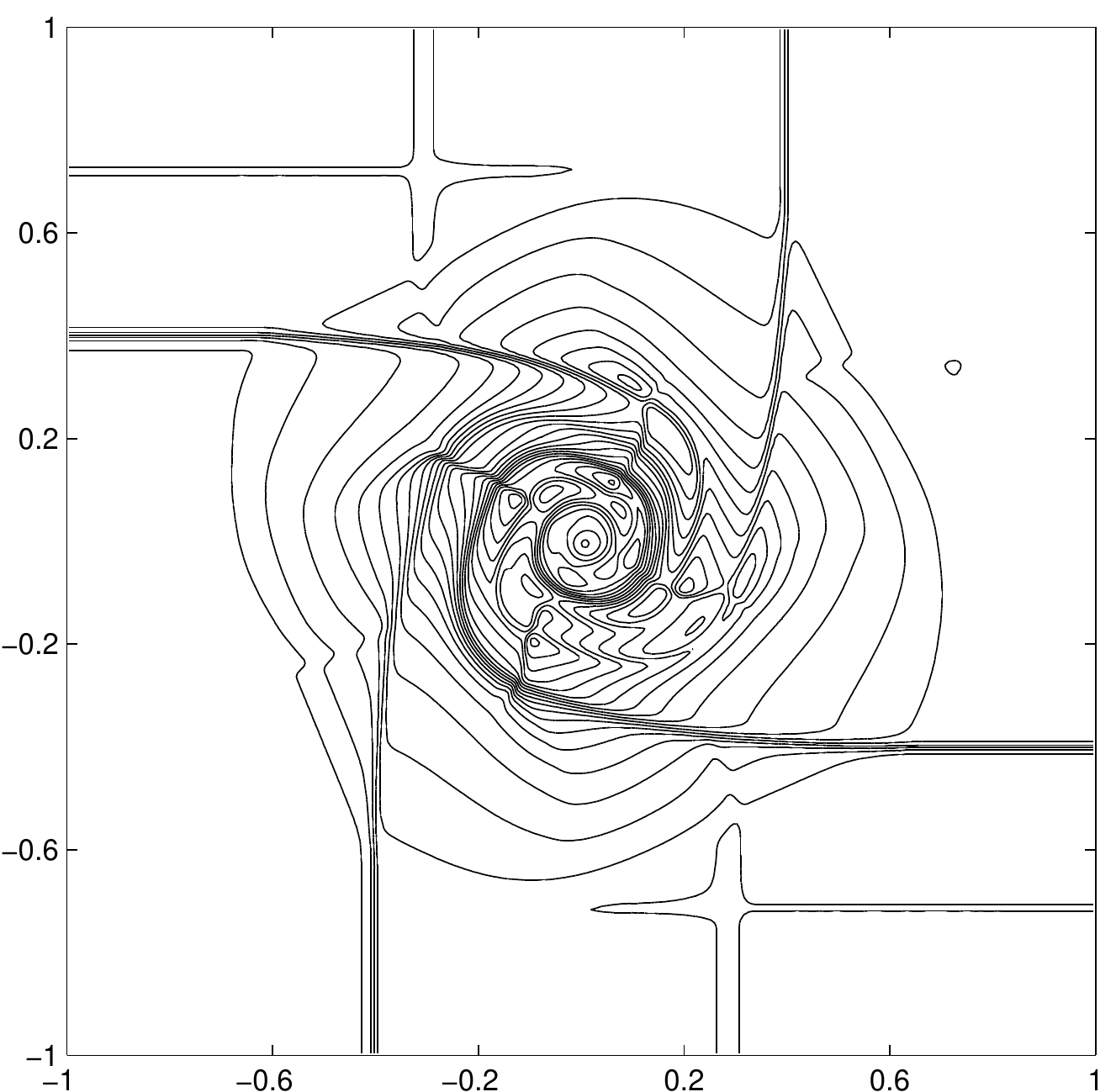}&
		\includegraphics[width=0.35\textwidth]{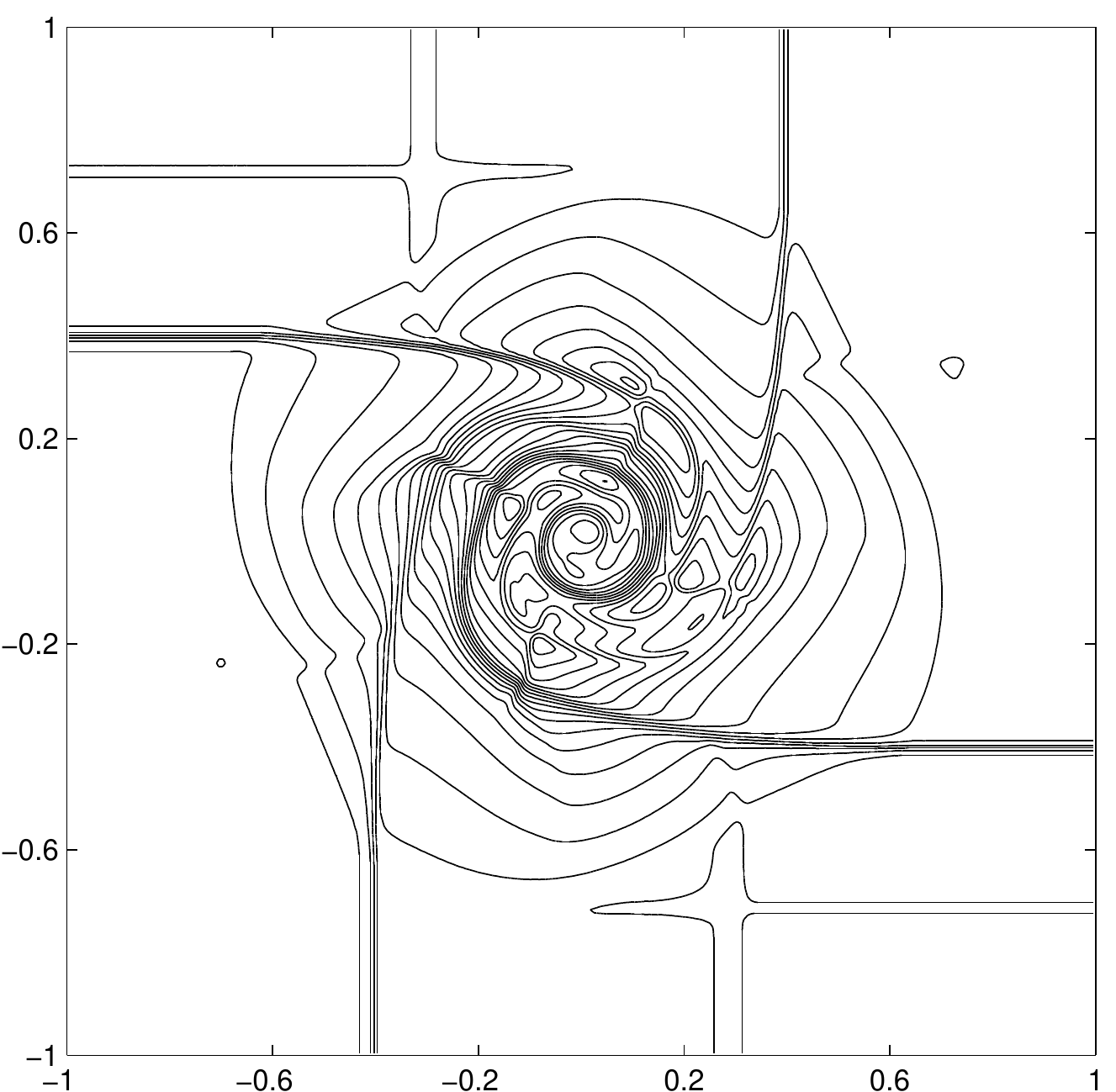}\\
		\includegraphics[width=0.35\textwidth]{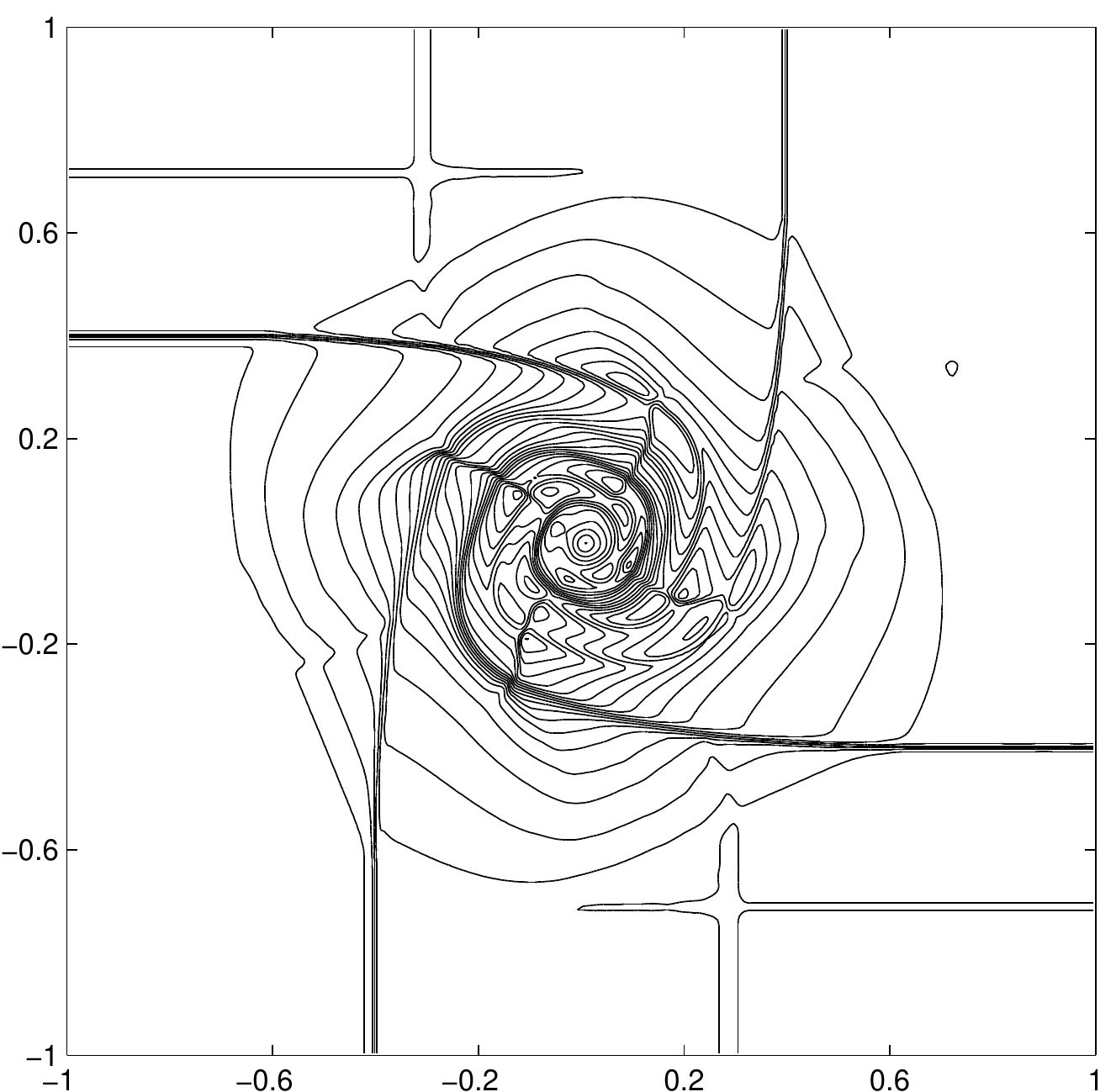}&
		\includegraphics[width=0.35\textwidth]{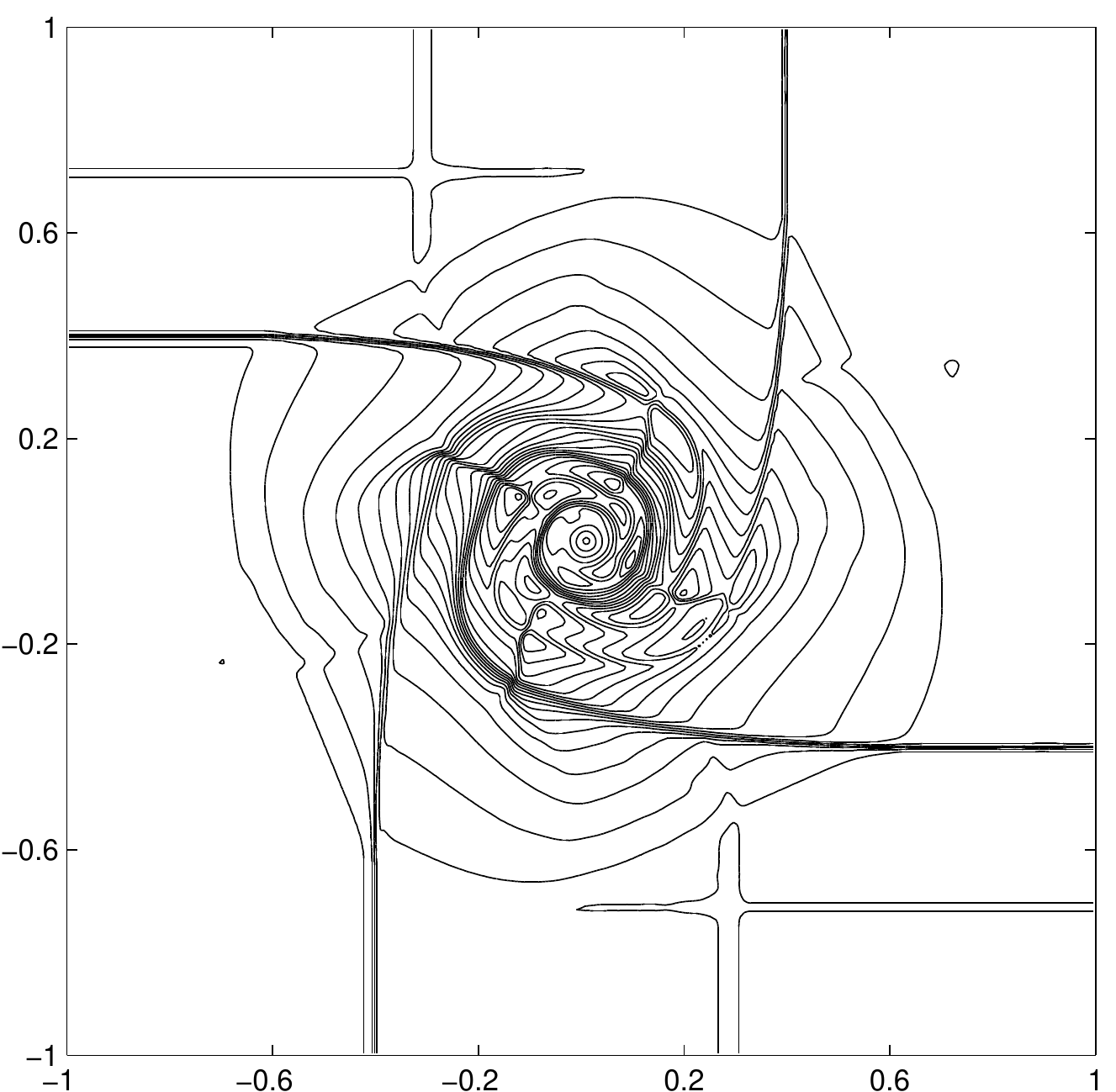}\\
	\end{tabular}
	\caption{Example \ref{exRHDRM2DT2Cdg}:
The contour plots of density logarithm  $\log\rho$ at
 $t=0.8$ (30 equally spaced contour lines from
$-1.98$ to 0.56) obtained by using $P^K$-based  \CDG{} with $300\times 300$ cells.
Left: the old, right: the new.
From top to bottom, $K=1,~2,~3$. }
	\label{fig:RHDT2cmpNOrho}
\end{figure}

\begin{Example}[Riemann problem 2]\label{exRHDRM2DT2Cdg}\rm
 It is about the interaction of four contact discontinuities (vortex sheets) with the same sign (the negative sign) for the ideal relativistic fluid. The initial data may be found in
\cite{ZhaoTang2013}.   
The new and old \CDG{} are used to solve this problem.
In order to compare the CPU times, the maximum CFL numbers in
Table~\ref{tab:cfl} and $\theta=1$ are  considered for them.
Fig.~\ref{fig:RHDT2cmpNOrho} displays the contours of density at $t=0.8$.
It is seen that the resolutions of  the new
$P^1$- and $P^2$-based \CDG{} are slightly worse  than the old, but
the difference between two $P^3$-based methods are not obvious.
Their CPU times listed in Table~\ref{tab:RHDT2cmpNO} show that
for the $P^2$ case, the new method is slower than the old,
but  the new $P^3$-based method has a great advantage.
It is mainly because the CFL number of the new $P^2$-based method is
about half of that of the old, while the CFL number of the new $P^3$-based method
is about two-thirds of the old. 
\end{Example}

{}

\begin{table}[!htbp]
	\centering
	\caption{Example~\ref{exRHDRM2DT2Cdg}: CPU times (second) of the new and old CDG methods. }
	\begin{tabular}{|c|c|c|}
		\hline
		& new method & old method\\
		\hline
		$P^1$& 2324.5 & 2460.3\\
		\hline
		$P^2$& 14210.2 & 13611.9\\
		\hline
		$P^3$& 30702.3 & 48538.7\\
		\hline
	\end{tabular}
	\label{tab:RHDT2cmpNO}
\end{table}


\begin{Example}[Shock and light bubble interaction]\label{exRHDSBT1cdg}\rm
	This example describes the  interaction between the  shock wave and a light bubble,
The setup of the problem is as follows.
Initially, within the computational domain $[0,325]\times [-45,45]$, there is a left-moving shock wave at
	$x=265$ with the left and right states
	$$(\rho,v_1,v_2,p)=\begin{cases}(1,0,0,0.05),&x<265,\\
	(1.865225080631180,-0.196781107378299,0,0.15),&x>265.
\end{cases}$$
A cylindrical bubble is centered at $(215,0)$ with the radius of 25 in front of the initial shock wave.
 The fluid within the bubble  is in a mechanical equilibrium with the surrounding
fluid and lighter than the ambient fluid. The detailed state of the
fluid within the bubble is taken as $(\rho,v_1,v_2,p)=(0.1358, 0,0,0.05)$.
In our computations, the domain is divided into $500\times 140$ uniform cells,
the CFL numbers $\mu$ for the
$P^1$-, $P^2$-, $P^3$-based \CDG{}  are taken as $0.6,0.5,0.4$, respectively, and
$\theta=0.5$. Moreover, the parameter $M$ in the modified TVB minmod function is taken as $0.005$,
the reflective boundaries are specified at $y=\pm 45$,
and the fluid states on two boundaries in $x$-direction
are set to the left and right shock wave states, respectively.

Figs.~\ref{fig:RHDSBT1Ite3rho} and \ref{fig:RHDSBT1Ite3cell}
show the schlieren images of  density and the distributions of ``troubled'' cells at $t=270$,
while Figs.~\ref{fig:RHDSBT1Ite5rho} and \ref{fig:RHDSBT1Ite5cell}~ give corresponding results at
$t=450$.
It is seen from those plots
that the $P^K$-based \DG{} resolve the complex wave structure better than the $P^K$-based \CDG{},
$P^3$-based \DG{} gives relatively fine wave structure, the distribution of ``troubled'' cells is very consistent with  the solutions,
and  the ``troubled'' cell proportions are almost similar
for the same order of the two types of DG methods.
Table~\ref{tab:cellperSB1} gives the percentage of    ``troubled'' cells at two different times.
The CPU times are estimated in Table~\ref{tab:RHDSB1cmpRC} and
show that the $P^1$- and $P^3$-based \CDG{} are faster than corresponding \DG{},  but
 the $P^2$-based method is an exception.
\end{Example}

\ifx\outnofig\undefine
\begin{figure}[!htbp]
	\centering{}
	\begin{tabular}{cc}
		
		\includegraphics[width=0.35\textwidth]{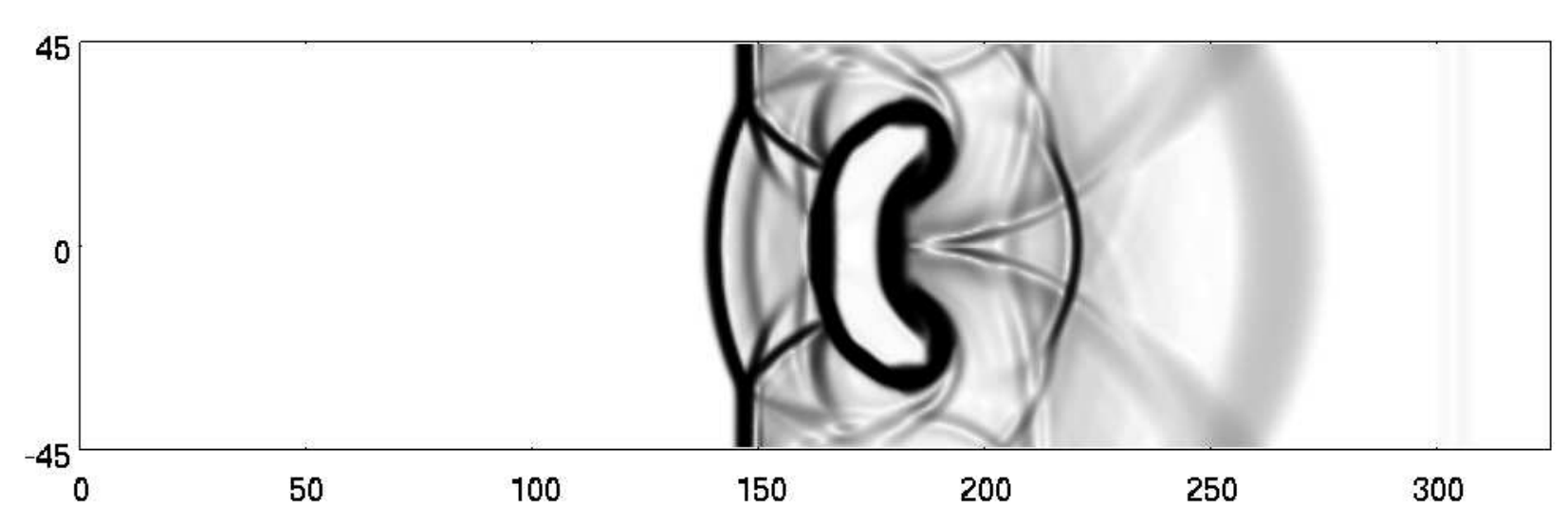}&
		\includegraphics[width=0.35\textwidth]{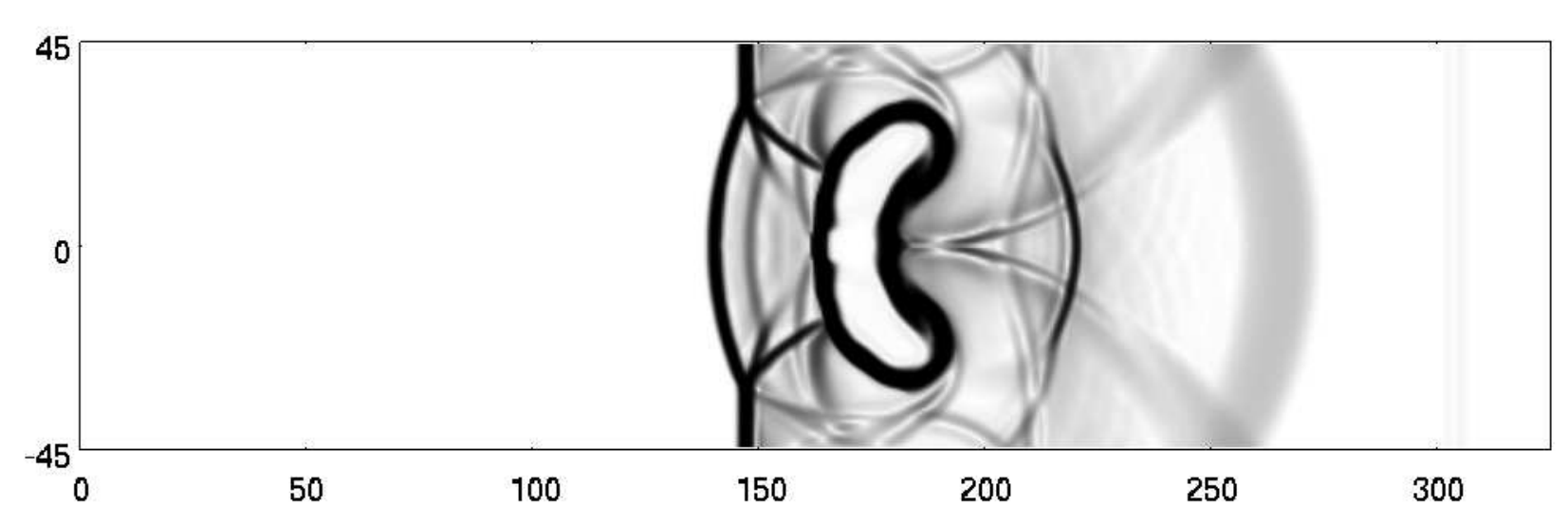}\\
		\includegraphics[width=0.35\textwidth]{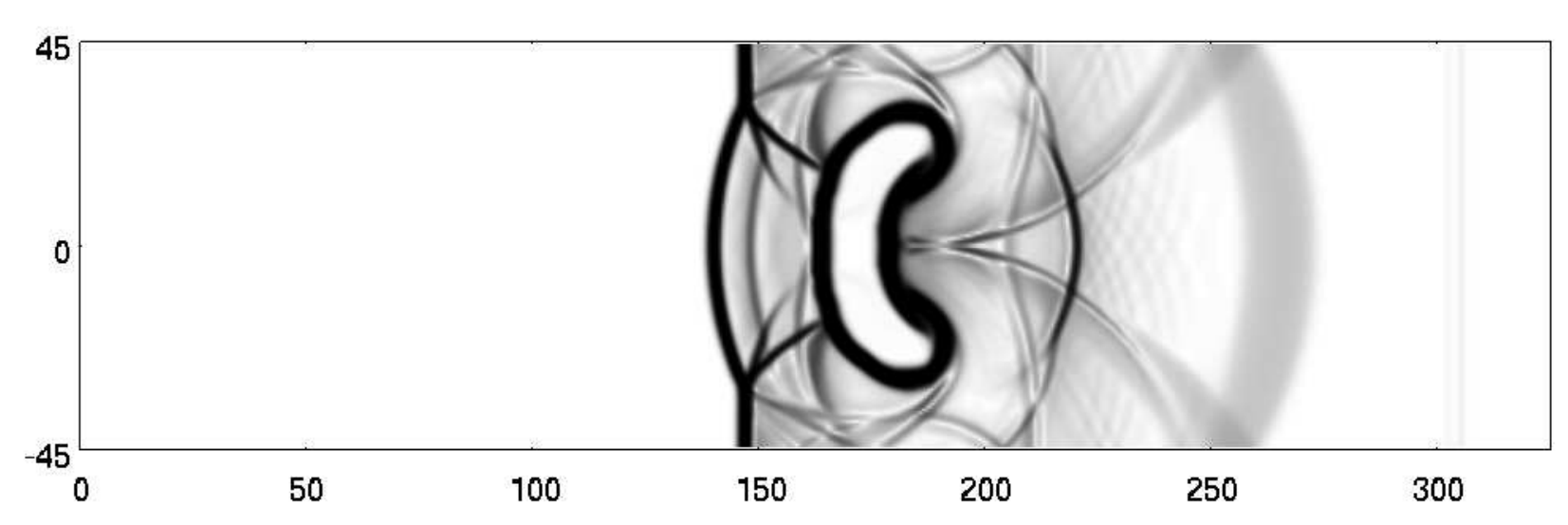}&
		\includegraphics[width=0.35\textwidth]{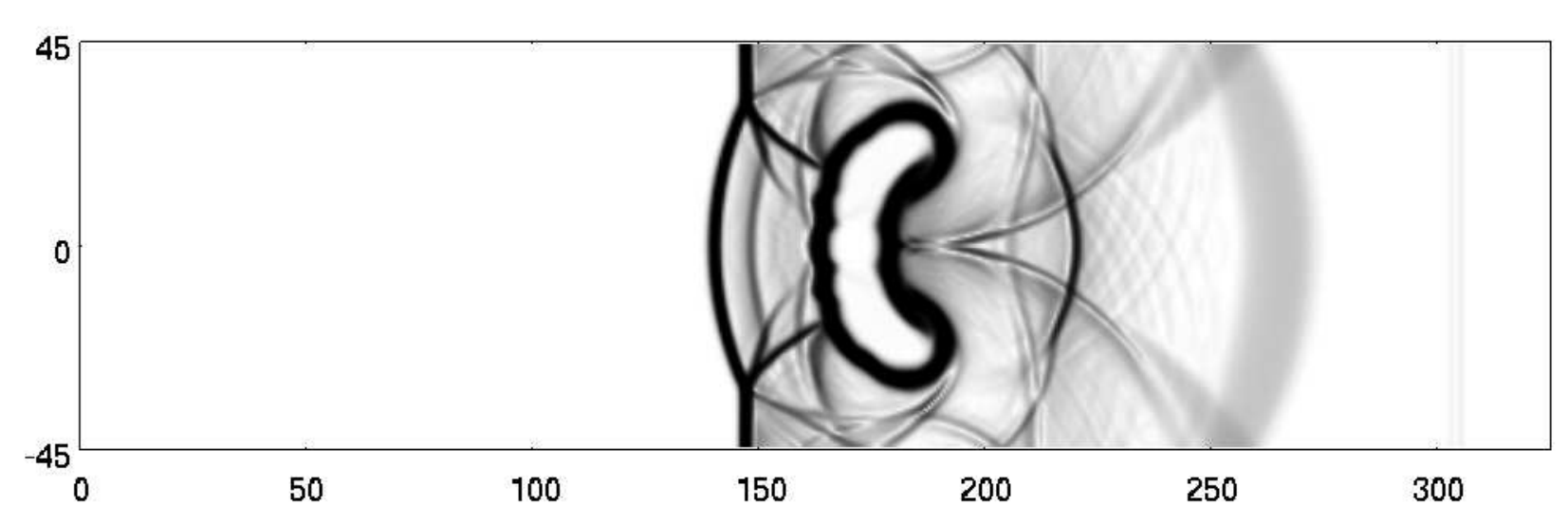}\\
		\includegraphics[width=0.35\textwidth]{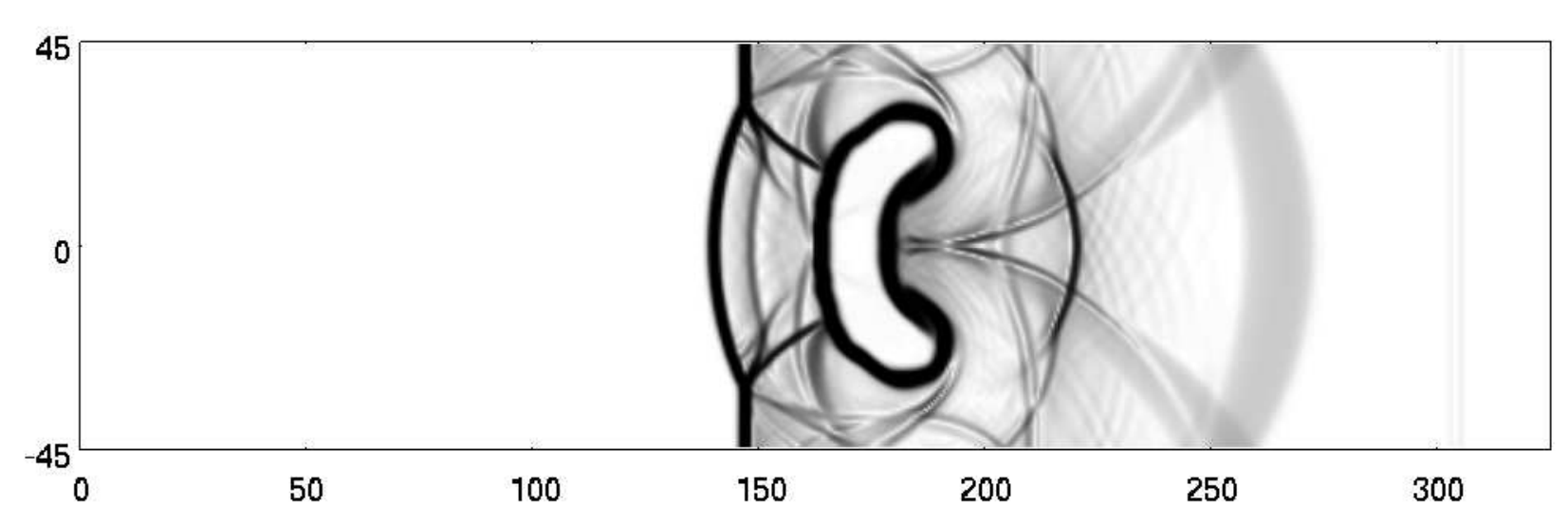}&
		\includegraphics[width=0.35\textwidth]{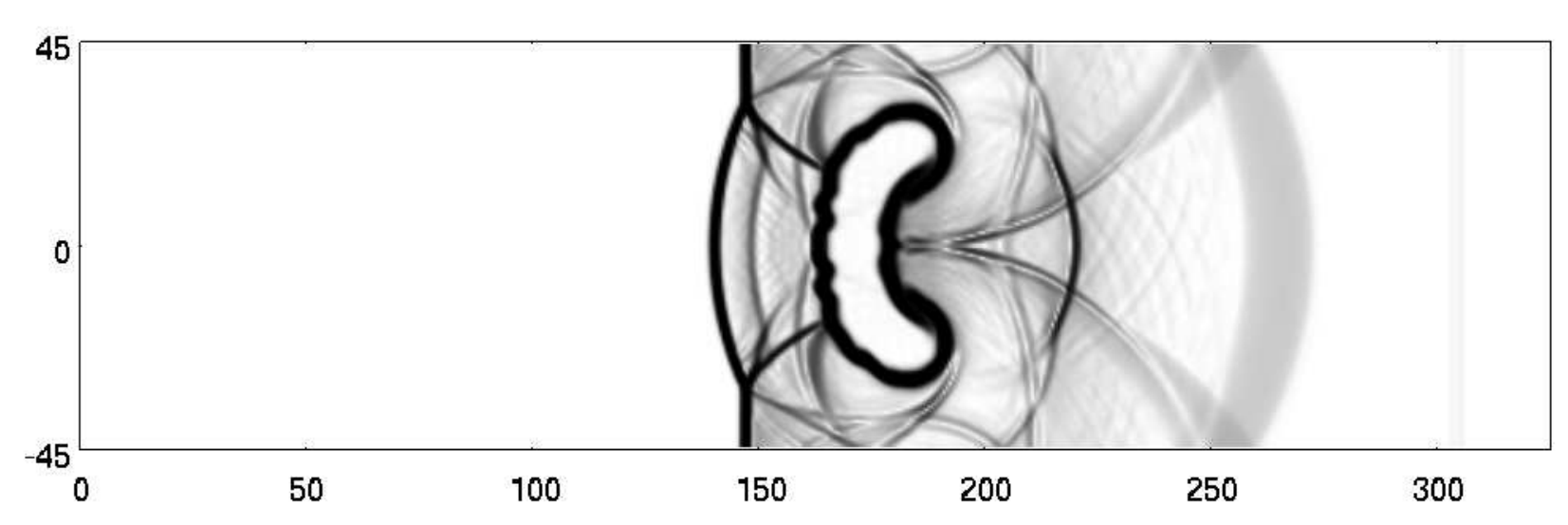}
	\end{tabular}
	\caption{Example~\ref{exRHDSBT1cdg}:
The  schlieren images of density at $t=270$ with $500\times
		140$ cells. Left:
		$P^K$-based \CDG{}; right: $P^K$-based \DG{}. From top to bottom, $K=1,2,3$.}
	\label{fig:RHDSBT1Ite3rho}
\end{figure}

\begin{figure}[!htbp]
	\centering{}
	\begin{tabular}{cc}
		\includegraphics[width=0.35\textwidth]{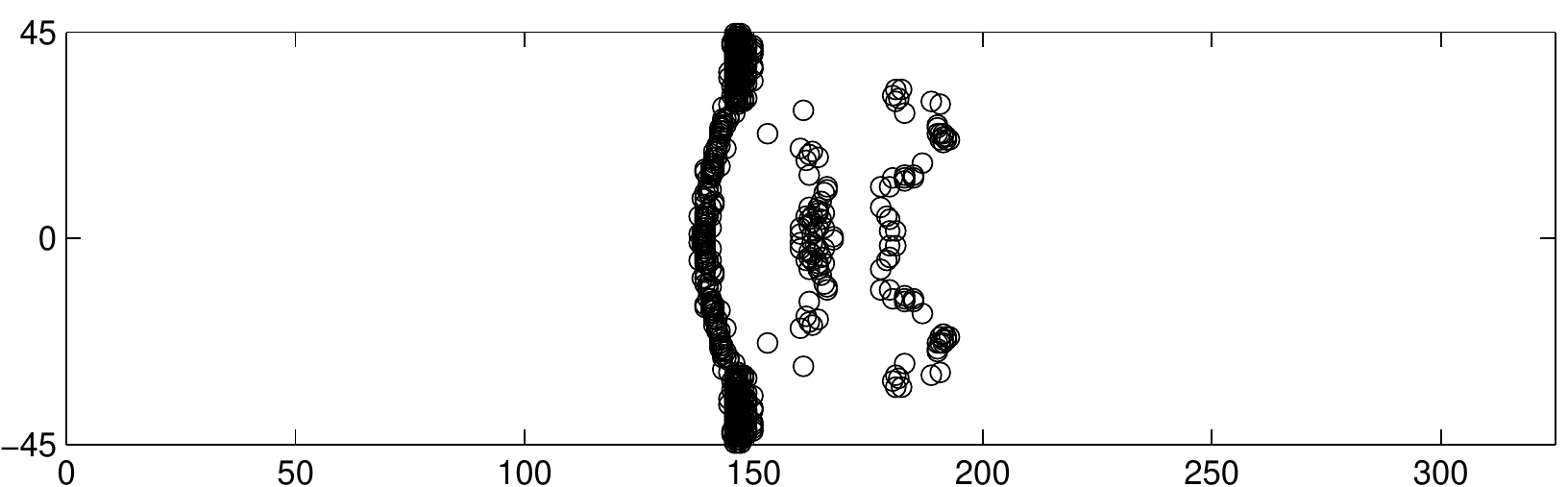}&
		\includegraphics[width=0.35\textwidth]{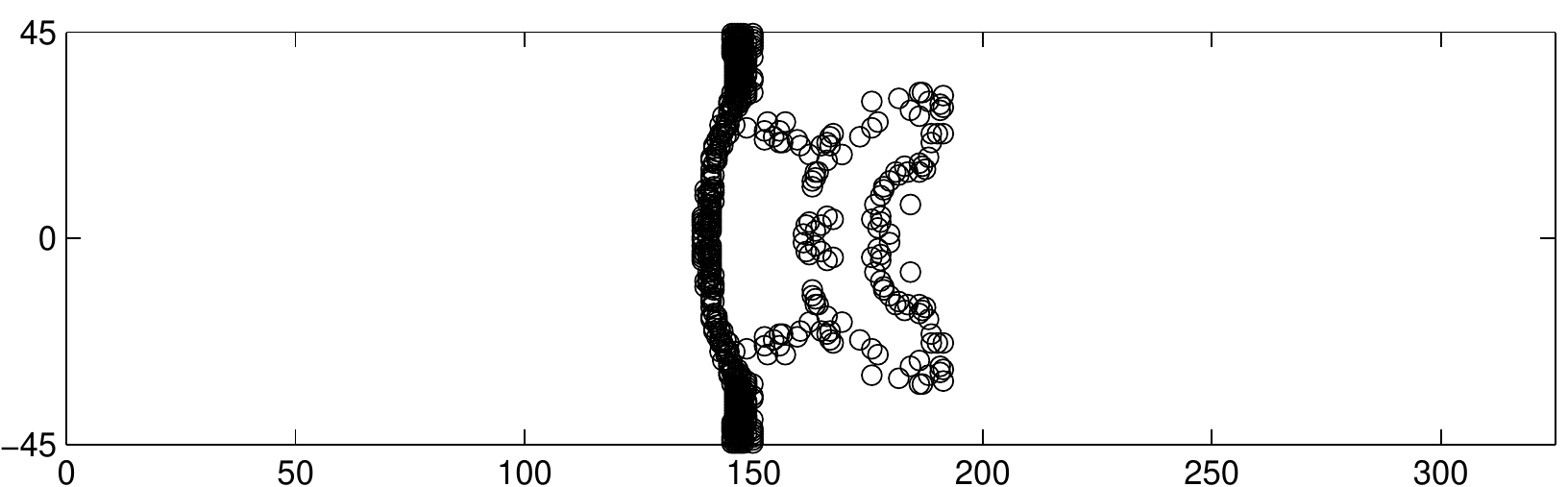}\\
		\includegraphics[width=0.35\textwidth]{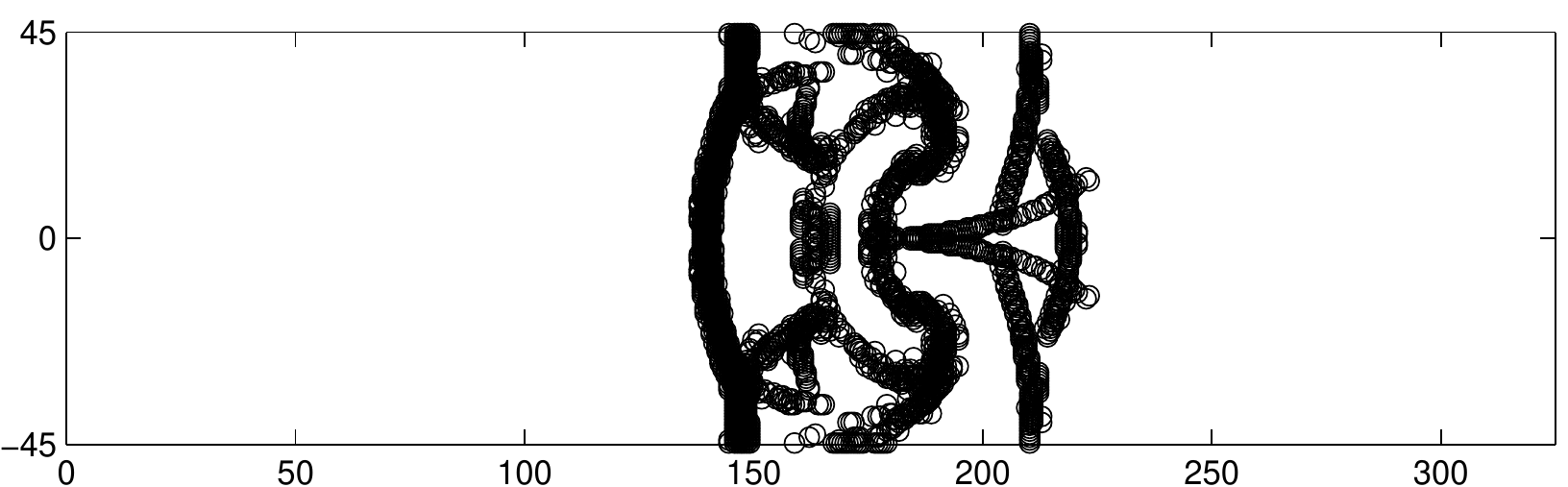}&
		\includegraphics[width=0.35\textwidth]{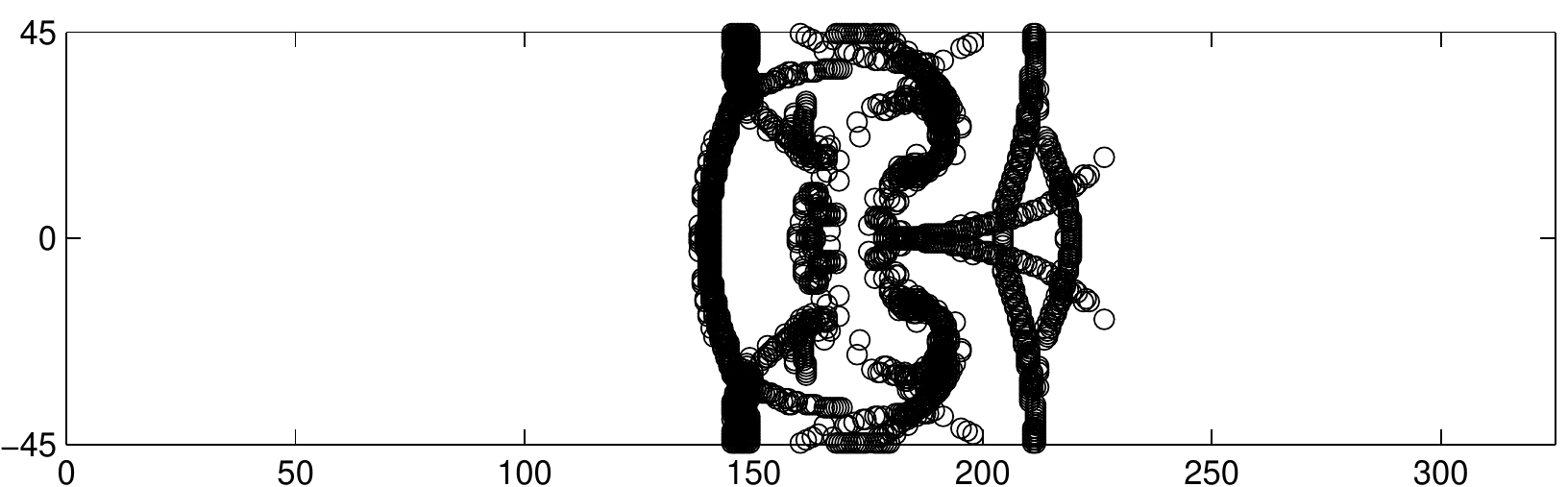}\\
		\includegraphics[width=0.35\textwidth]{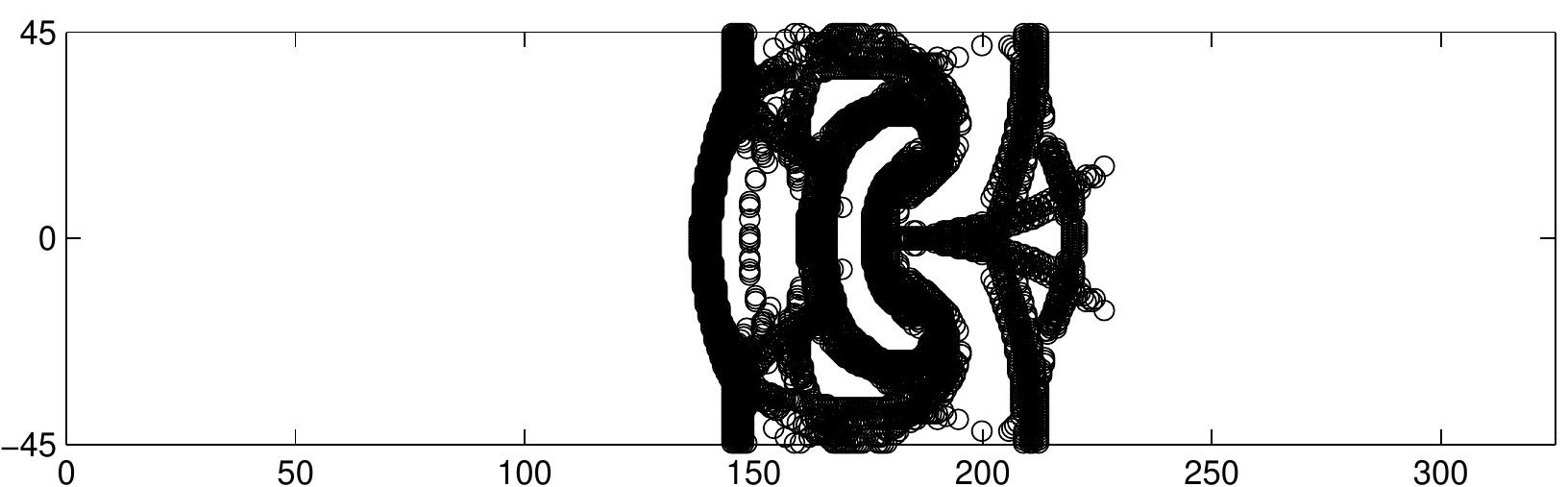}&
		\includegraphics[width=0.35\textwidth]{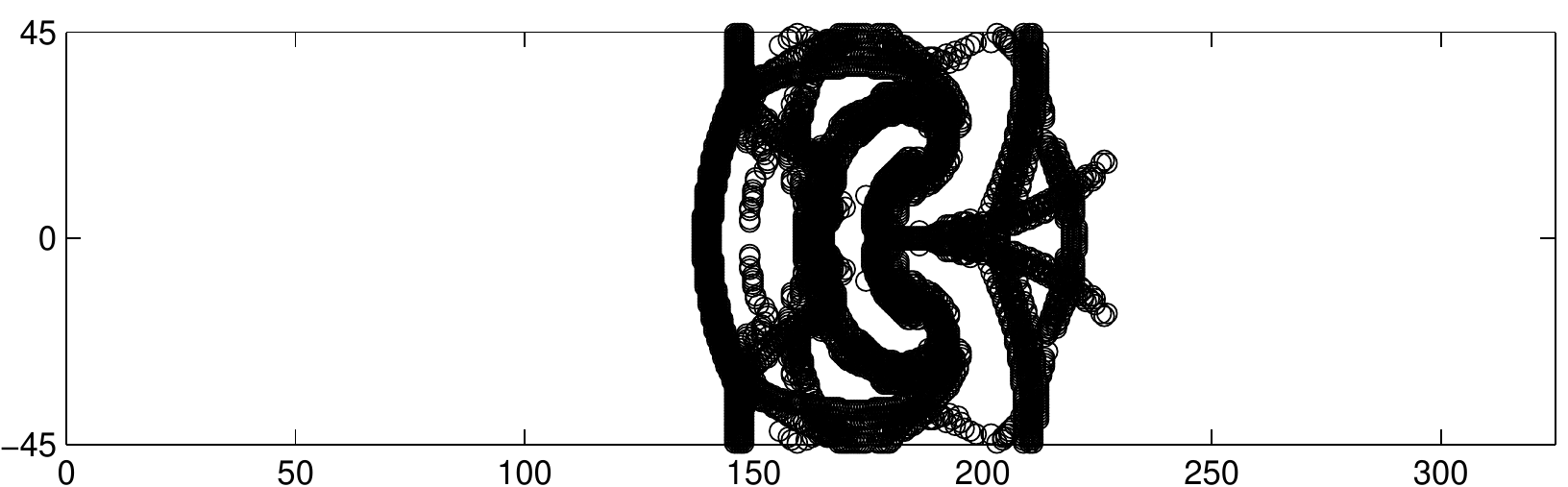}
	\end{tabular}
	\caption{Same as Fig.~\ref{fig:RHDSBT1Ite3rho} except for the ``troubled'' cells at $t=270$. }
	\label{fig:RHDSBT1Ite3cell}
\end{figure}

\begin{figure}[!htbp]
	\centering{}
	\begin{tabular}{cc}
		
		\includegraphics[width=0.35\textwidth]{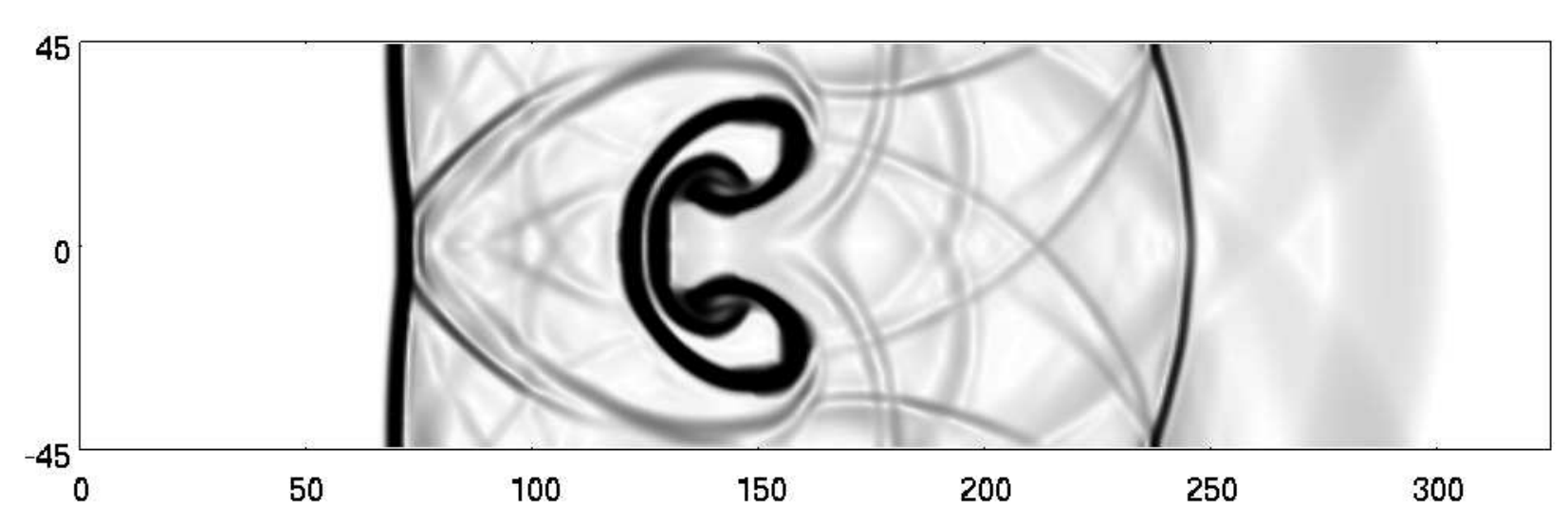}&
		\includegraphics[width=0.35\textwidth]{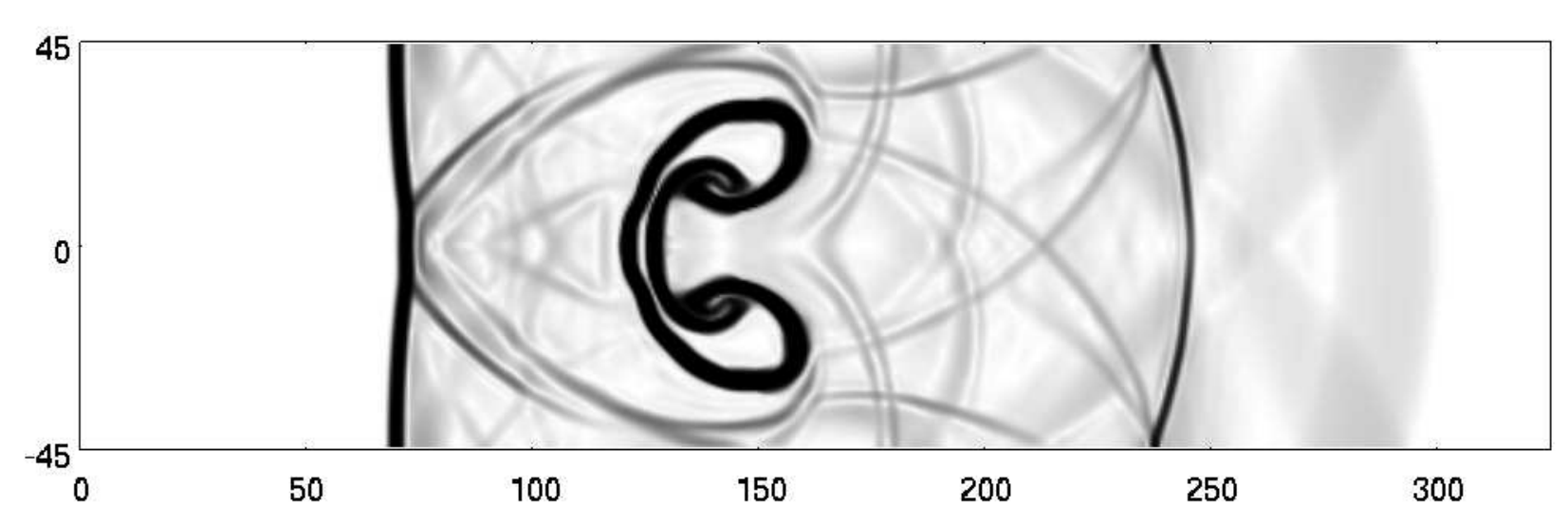}\\
		\includegraphics[width=0.35\textwidth]{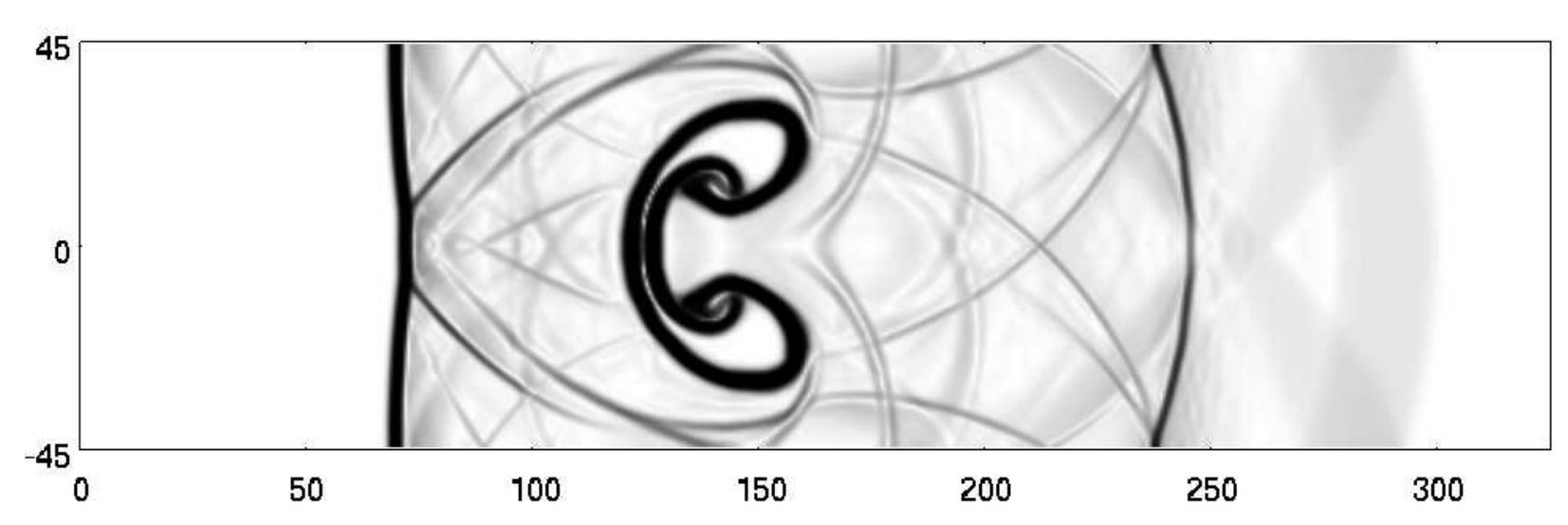}&
		\includegraphics[width=0.35\textwidth]{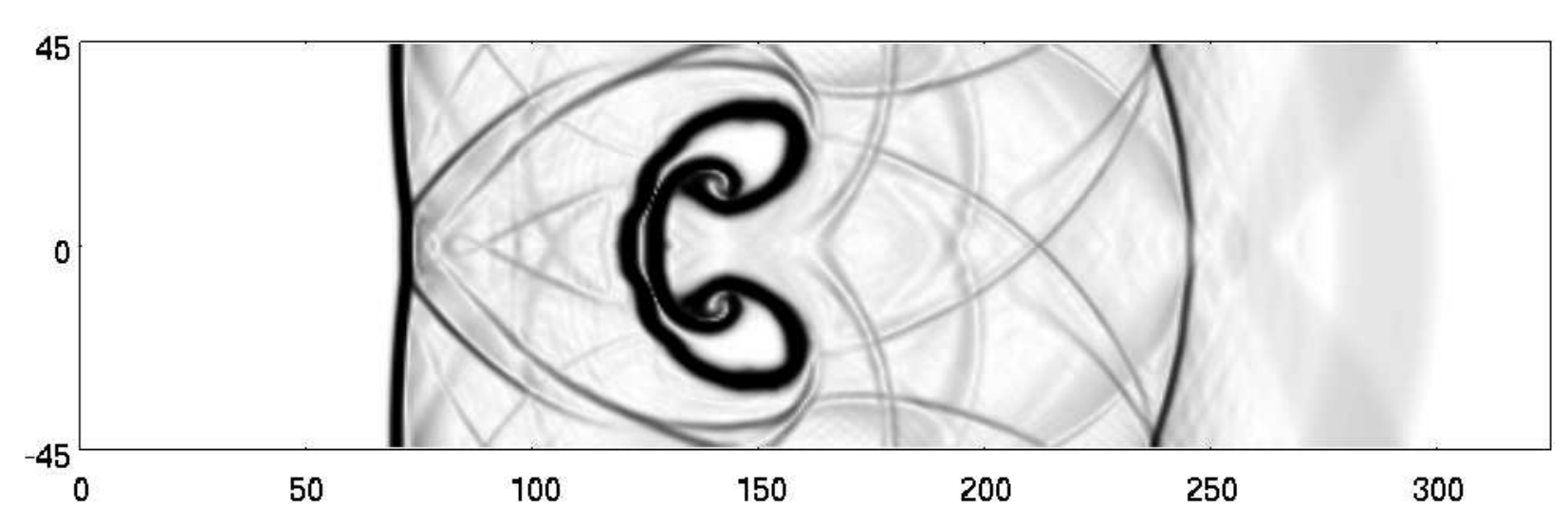}\\
		\includegraphics[width=0.35\textwidth]{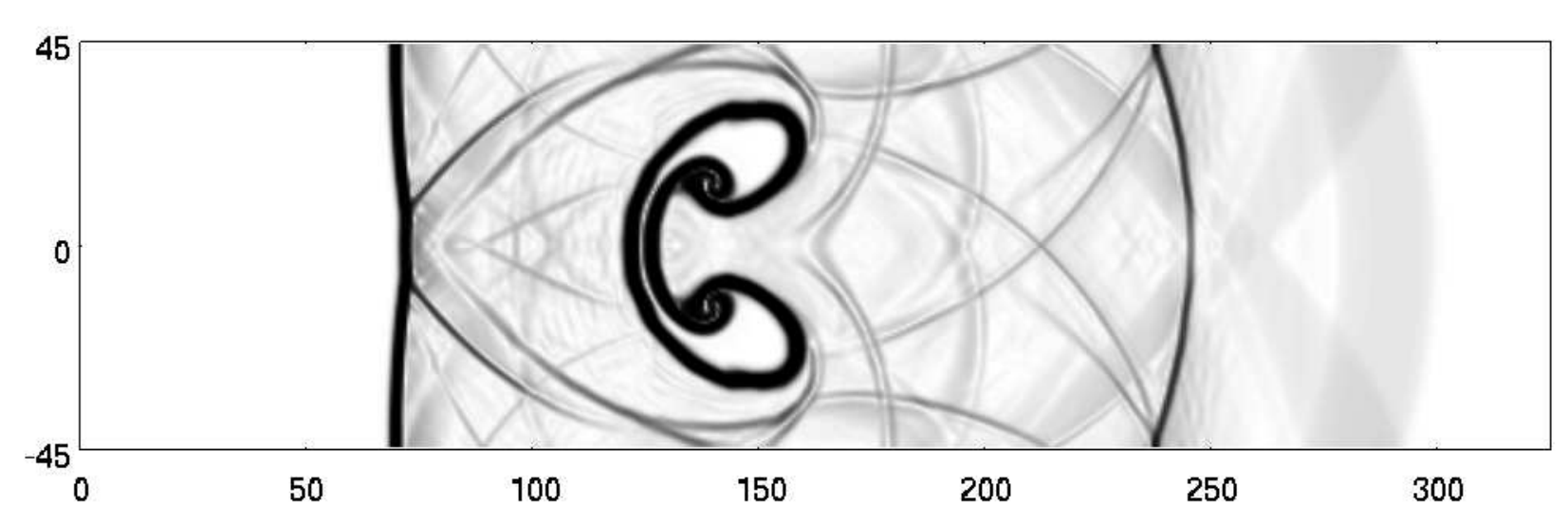}&
		\includegraphics[width=0.35\textwidth]{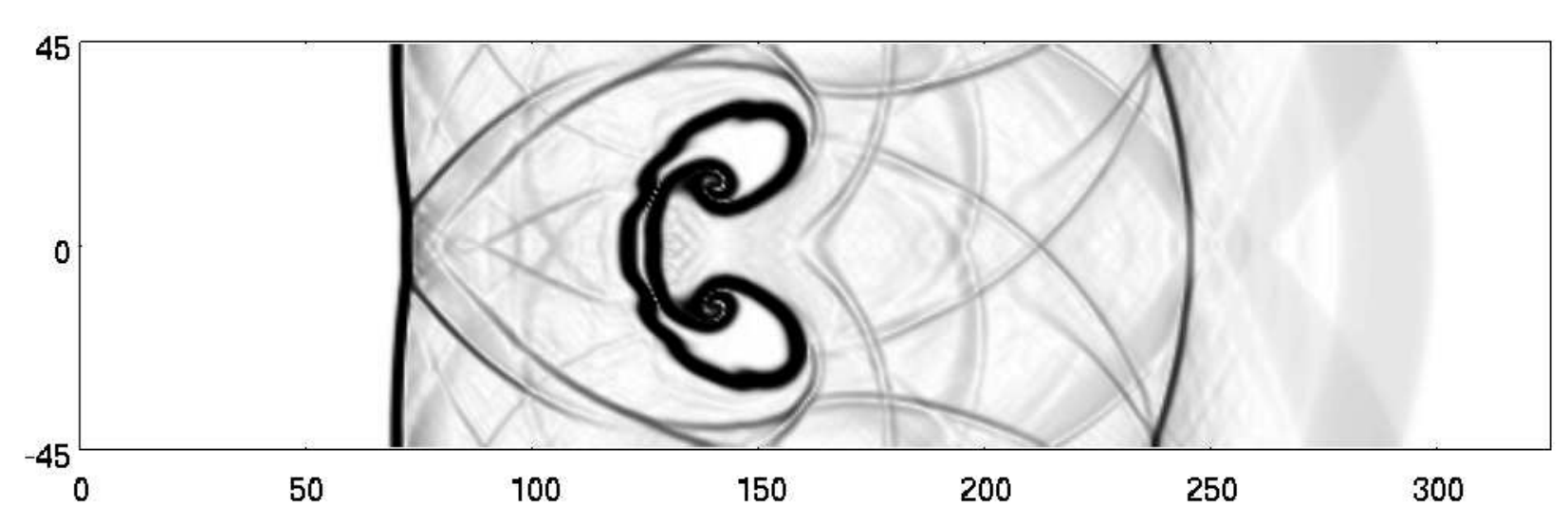}
	\end{tabular}
	\caption{Same as Fig.~\ref{fig:RHDSBT1Ite3rho} except for the  schlieren image of density at $t=450$. }
	\label{fig:RHDSBT1Ite5rho}
\end{figure}

\begin{figure}[!htbp]
	\centering{}
	\begin{tabular}{cc}
		\includegraphics[width=0.35\textwidth]{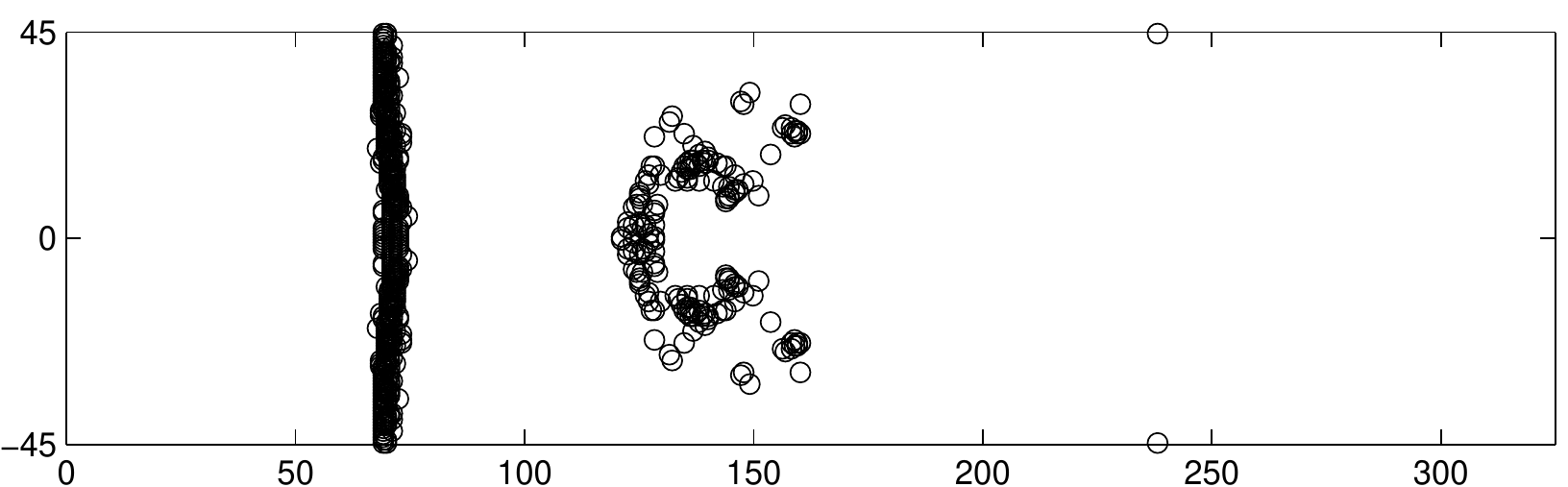}&
		\includegraphics[width=0.35\textwidth]{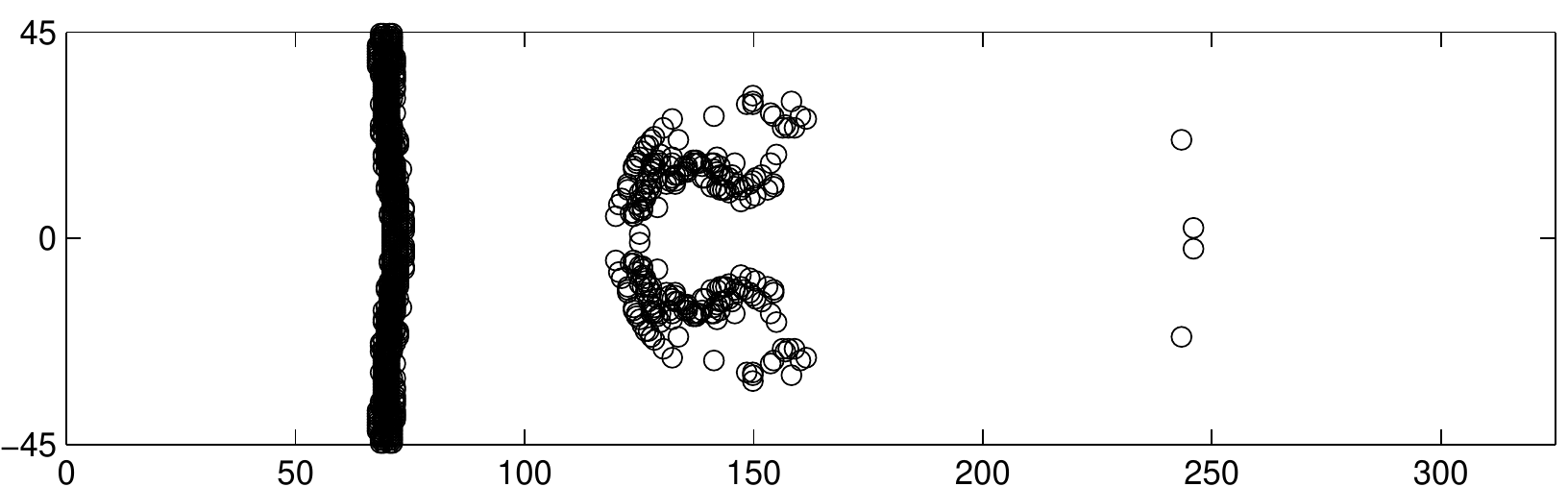}\\
		\includegraphics[width=0.35\textwidth]{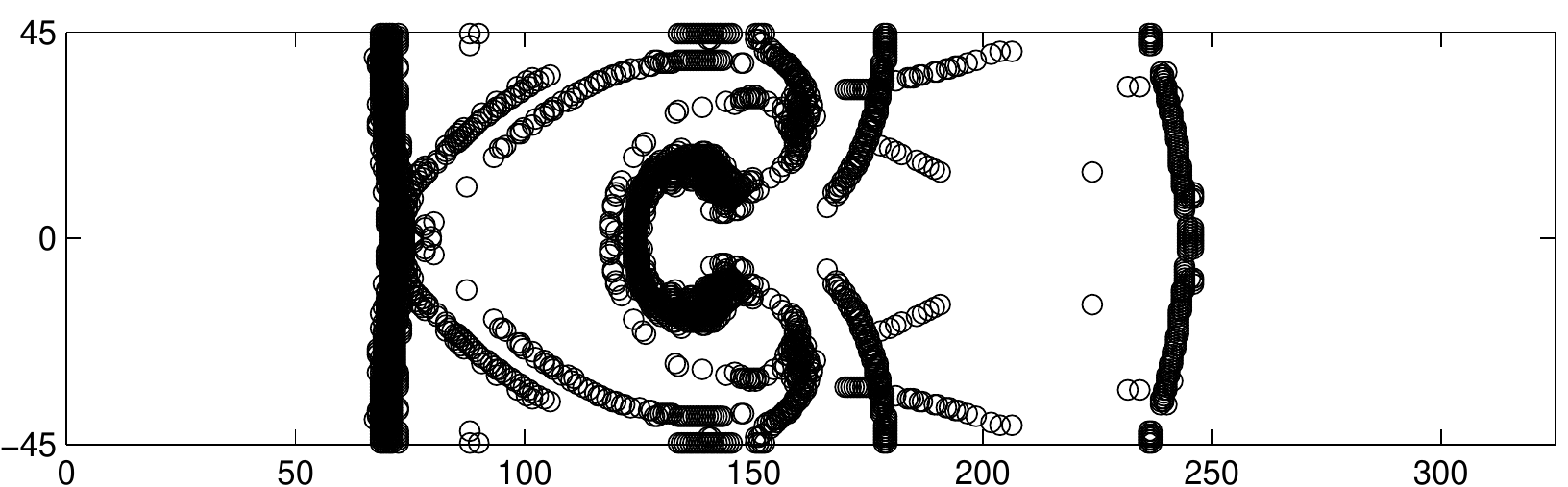}&
		\includegraphics[width=0.35\textwidth]{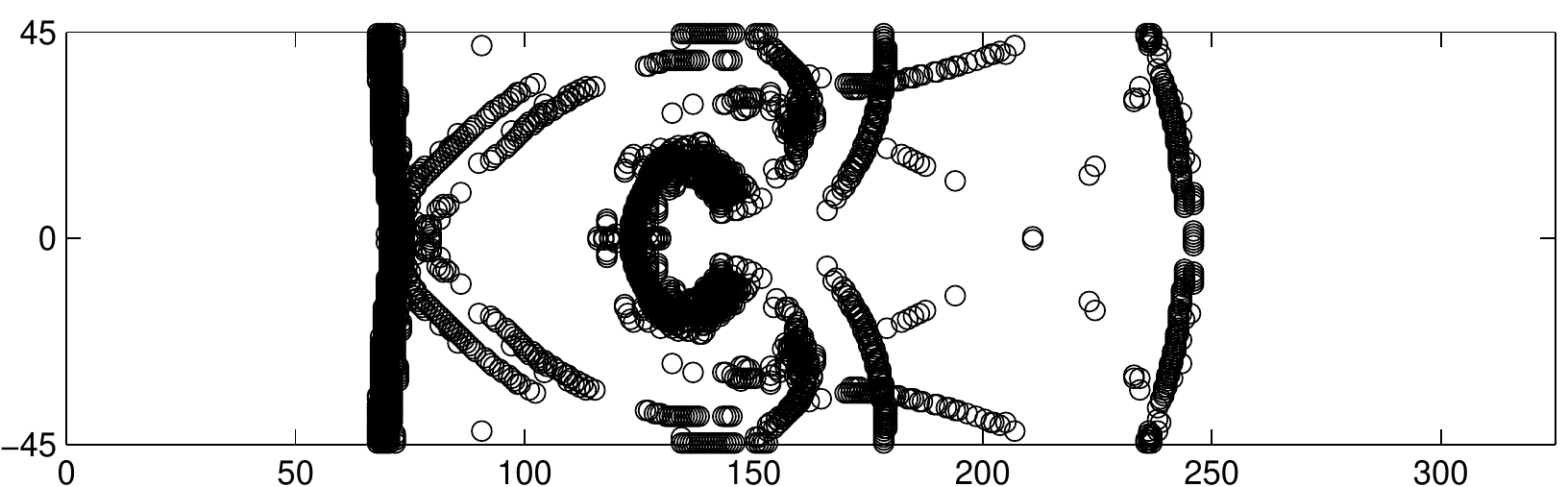}\\
		\includegraphics[width=0.35\textwidth]{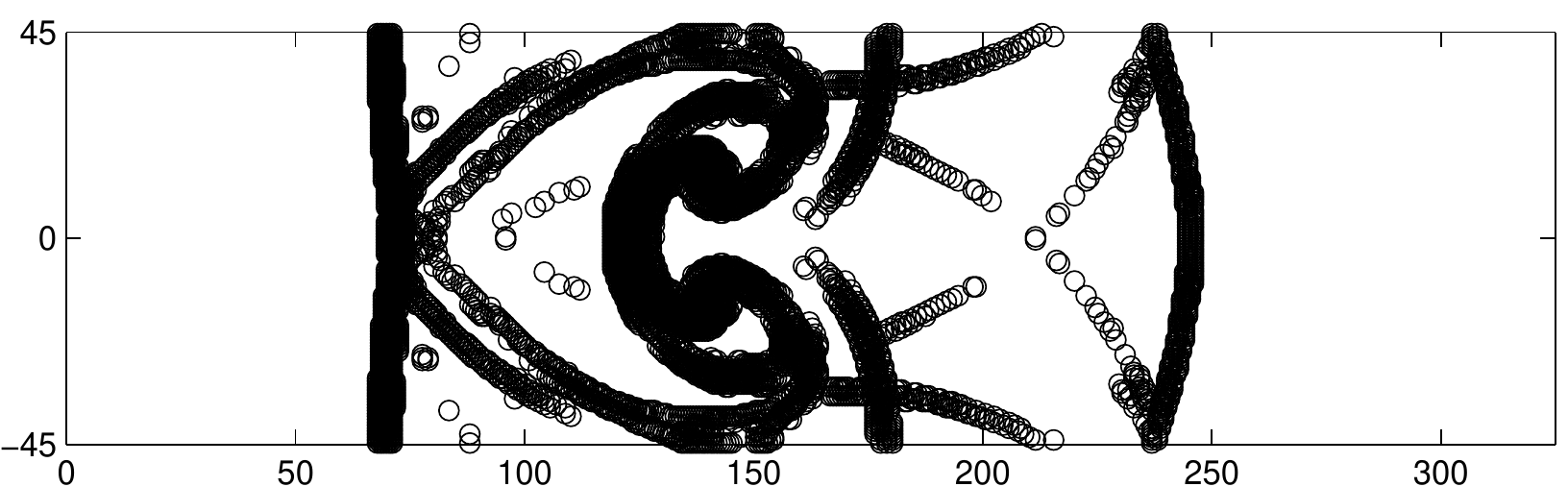}&
		\includegraphics[width=0.35\textwidth]{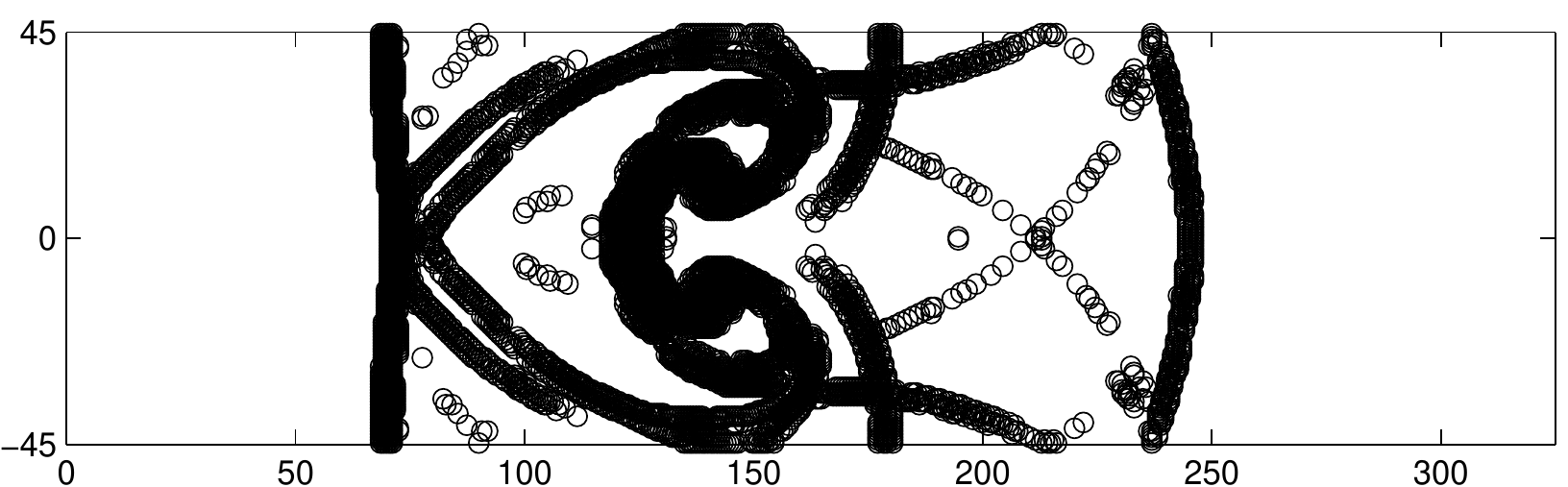}
	\end{tabular}
	\caption{Same as Fig.~\ref{fig:RHDSBT1Ite3rho} except for the ``troubled'' cells at  $t=450$. }
	\label{fig:RHDSBT1Ite5cell}
\end{figure}
\fi{}

\begin{table}[!htbp]
	\centering
	\caption{Example~\ref{exRHDSBT1cdg}: The percentage of ``troubled'' cells at two different times.}
	\begin{tabular}{|c|c|c|c|c|}
		\hline
		\multirow{2}{35pt}{} & \multicolumn{2}{|c|}{$t=270$} & \multicolumn{2}{|c|}{$t=450$}
		\\
		\cline{2-5}
		& CDG & non-central DG & CDG & non-central DG\\
		\hline
		$P^1$& 0.68 & 0.22  & 0.71 & 0.33 \\
		\hline
		$P^2$& 3.12 & 3.58 & 3.40 & 3.56 \\
		\hline
		$P^3$& 5.61 & 5.94 & 7.11 & 7.39\\
		\hline
	\end{tabular}
	\label{tab:cellperSB1}
\end{table}

\begin{table}[!htbp]
	\centering
	\caption{Example \ref{exRHDSBT1cdg}: CPU times (second) of
the CDG and non-central DG methods. }
	\begin{tabular}{|c|c|c|}
		\hline
		& CDG &  non-central DG \\
		\hline
		$P^1$& 1.93e4 & 4.03e3 \\
		\hline
		$P^2$& 6.17e4 & 1.14e4 \\
		\hline
		$P^3$& 1.73e5 & 3.47e4  \\
		\hline
	\end{tabular}
	\label{tab:RHDSB1cmpRC}
\end{table}


\section{Conclusions}
\label{Section-conclusion}

It is much more difficult to solve the relativistic hydrodynamical (RHD) equations
than the non-relativistic case. The appearance of Lorentz factor
enhances the nonlinearity of the RHD equations,  the fluxes
can not be formulated in an explicit
form of the conservative vector, and there are some inherent physical constraints
on the physical state. In practical computations of the RHD system,
 the primitive variable vector has to be first recovered   from the known
conservative vector   by iteratively solving a nonlinear
pressure equation and then the fluxes are  evaluated at each time step.
%

The paper developed the $P^K$-based \CDG{} for the one- and two- dimensional
special RHD equations,  $K=1,2,3$.
In comparison to \DG{}, the \CDG{} found two approximate solutions defined
on  mutually dual meshes.
For  each mesh,
the CDG approximate solutions on its dual mesh were used to calculate the flux values in the cell and
on the cell boundary so that
    the approximate solutions on two  mutually dual meshes were coupled with each other,
    and  the use of numerical flux might be avoided.
In addition, the \CDG{} allowed the use of a larger  CFL number.

The WENO limiter was adaptively implemented via two steps:  the
``troubled'' cells were first identified by using a modified TVB minmod function,
and then  the WENO technique is used to locally reconstruct  new polynomials of degree $(2K + 1)$
replacing the CDG solutions inside the ``troubled'' cells by using the cell average values of the CDG solutions in the neighboring cells
as well as the original cell averages of the ``troubled'' cells.

%
The accuracy of the CDG without the numerical
dissipation was analyzed and calculation of the flux integrals over the cells was also discussed.
Because the DG approximate solutions were discontinuous at the cell interface in general,
the integrals over each cell of the DG approximate solutions defined on corresponding
dual meshes became a sum of several integrals over subcell of the DG polynomial solutions
 on corresponding dual meshes,
which would lead to that more numerical integration  points are needed
to ensure  the accuracy of \CDG.
An attempt was made that
 only approximate solution on the dual mesh was
 used to evaluate the flux on the cell boundary
 for the DG approximations of RHD system on the mesh,
 thus  the integrals of  DG solutions over the dual cell might
 be avoided and  the cost of numerical integration was hopefully reduced.
%
For the linear scalar equation,  the Fourier method and numerical experiments were used to analyze
the stability of such new method, and estimate  the CFL numbers for the stability.

Several numerical experiments  demonstrated the accuracy, robustness,  and
discontinuity resolution of our methods.
The  results showed that
 the \CDG{} with WENO limiter were robust and could  capture the contact discontinuities,
shock waves, and other complex wave structures well,
 the WENO limiter was only implemented in a few ``troubled'' cells,
the \CDG{} had obvious advantages in simulating the propagation of slow shock wave in comparison to the \DG{}.
Moreover, the new two-dimensional $P^3$-based \CDG{} could more significantly improve the computational efficiency
than the old. 
%
In solving RHD problems with large Lorentz factor, or strong discontinuities,
or low rest-mass density or pressure etc., it is still possible for the $P^K$-based \CDG{}
to give nonphysical solutions. To cure such difficulty, the $P^0$-based method
may be locally used to replace  the $P^K$-based. The genuinely effective way
is to employ the physical-constraints preserving methods, see e.g.
\cite{Wu-Tand-JCP2015,Wu-Tang-RHD2016}.


\section*{Acknowledgements}
This work was partially supported by the National Natural Science
Foundation of China (Nos. 91330205 \& 11421101).

\bibliography{ref/journalname,ref/pkuth}
           \bibliographystyle{plain}

\end{document}